\documentclass[12pt,a4paper,fleqn]{article}
\usepackage{tikz}
\usepackage{caption}
\usepackage{lscape}

\usepackage{fullpage}
\usepackage{amsmath}
\usepackage{amsthm}
\usepackage{amssymb}
\usepackage{amscd}
\usepackage{plain}
\usepackage{graphicx}
\usepackage{t1enc,epsfig}
\usepackage[mathscr]{eucal}
\usepackage{epstopdf}
\usepackage[german,english]{babel}
\usepackage{stmaryrd}

\usepackage{xcolor}
\usepackage{float}
\newtheorem{theorem}{Theorem}[section]

\newtheorem{corollary}[theorem]{Corollary}
\newtheorem{definition}[theorem]{Definition}
\newtheorem{example}[theorem]{Example}
\newtheorem{lemma}[theorem]{Lemma}

\newtheorem{proposition}[theorem]{Proposition}
\newtheorem{remark}[theorem]{Remark}
\newtheorem{conjecture}[theorem]{Conjecture}
\newtheorem*{thma}{Theorem 1}
\newtheorem*{thmb}{Theorem 2}
\newtheorem*{thmc}{Theorem 3}
\newtheorem*{thmd}{Theorem 4}
\newtheorem*{thme}{Theorem 5}
\newtheorem*{thmf}{Theorem 6}
\newtheorem*{thmg}{Theorem 7}
\newtheorem*{thmh}{Theorem 8}
\newtheorem*{thmi}{Theorem 9}
\newtheorem*{thmj}{Theorem 10}
\newtheorem*{thms}{Theorem S}
\newtheorem*{lemI}{Lemma I}
\newtheorem*{lemII}{Lemma II}
\newtheorem*{lema}{Lemma 3.4.1}
\newtheorem*{lemb}{Lemma 3.4.2}
\newtheorem*{lemc}{Lemma 3.4.3}
\newtheorem*{lemd}{Lemma 3.4.4}
\newtheorem*{leme}{Lemma 3.4.5}
\newtheorem*{lemf}{Lemma 3.5}
\newtheorem*{lemg}{Lemma 3.6}
\newtheorem*{lemh}{Lemma 3.7}
\newtheorem*{lemi}{Lemma 3.8}
\newtheorem*{lemj}{Lemma 3.9}
\newtheorem*{remarka}{Remark 3.1}

\numberwithin{equation}{section}

\renewenvironment{proof}{{\bf Proof. }}{\hfill$\rule{1ex}{1ex}$\par\vskip 5 truemm}

\newenvironment{prooflemI}{{\bf Proof of Lemma I. }}{\hfill$\rule{1ex}{1ex}$\par\vskip 5 truemm}
\newenvironment{prooflemII}{{\bf Proof of Lemma II. }}{\hfill$\rule{1ex}{1ex}$\par\vskip 5 truemm}

\begin{document}

\newcommand{\bthm}{\begin{theorem}}
\newcommand{\ethm}{\end{theorem}}
\newcommand{\bthma}{\begin{thma}}
\newcommand{\ethma}{\end{thma}}
\newcommand{\bthmb}{\begin{thmb}}
\newcommand{\ethmb}{\end{thmb}}
\newcommand{\bthmc}{\begin{thmc}}
\newcommand{\ethmc}{\end{thmc}}
\newcommand{\bthmd}{\begin{thmd}}
\newcommand{\ethmd}{\end{thmd}}
\newcommand{\bthme}{\begin{thme}}
\newcommand{\ethme}{\end{thme}}
\newcommand{\bthmf}{\begin{thmf}}
\newcommand{\bthms}{\begin{thms}}
\newcommand{\ethmf}{\end{thmf}}
\newcommand{\ethms}{\end{thms}}
\newcommand{\bthmg}{\begin{thmg}}
\newcommand{\ethmg}{\end{thmg}}
\newcommand{\bthmh}{\begin{thmh}}
\newcommand{\ethmh}{\end{thmh}}
\newcommand{\bthmi}{\begin{thmi}}
\newcommand{\ethmi}{\end{thmi}}
\newcommand{\bthmj}{\begin{thmj}}
\newcommand{\ethmj}{\end{thmj}}
\newcommand{\bd}{\begin{definition}}
\newcommand{\ed}{\end{definition}}
\newcommand{\bs}{\begin{proposition}}
\newcommand{\es}{\end{proposition}}
\newcommand{\bp}{\begin{proof}}
\newcommand{\ep}{\end{proof}}
\newcommand{\bpthmc}{\begin{proofc}}
\newcommand{\epthmc}{\end{proofc}}
\newcommand{\bplemI}{\begin{prooflemI}}
\newcommand{\eplemI}{\end{prooflemI}}
\newcommand{\bplemII}{\begin{prooflemII}}
\newcommand{\eplemII}{\end{prooflemII}}
\newcommand{\be}{\begin{equation}}
\newcommand{\ee}{\end{equation}}
\newcommand{\ul}{\underline}
\newcommand{\br}{\begin{remark}}
\newcommand{\er}{\end{remark}}
\newcommand{\bra}{\begin{remarka}}
\newcommand{\era}{\end{remarka}}
\newcommand{\bex}{\begin{example}}
\newcommand{\eex}{\end{example}}
\newcommand{\bc}{\begin{corollary}}
\newcommand{\ec}{\end{corollary}}
\newcommand{\bl}{\begin{lemma}}
\newcommand{\el}{\end{lemma}}
\newcommand{\blI}{\begin{lemI}}
\newcommand{\elI}{\end{lemI}}
\newcommand{\blII}{\begin{lemII}}
\newcommand{\elII}{\end{lemII}}
\newcommand{\bla}{\begin{lema}}
\newcommand{\ela}{\end{lema}}
\newcommand{\blb}{\begin{lemb}}
\newcommand{\elb}{\end{lemb}}
\newcommand{\blc}{\begin{lemc}}
\newcommand{\elc}{\end{lemc}}
\newcommand{\bld}{\begin{lemd}}
\newcommand{\eld}{\end{lemd}}
\newcommand{\ble}{\begin{leme}}
\newcommand{\ele}{\end{leme}}
\newcommand{\blf}{\begin{lemf}}
\newcommand{\elf}{\end{lemf}}
\newcommand{\blg}{\begin{lemg}}
\newcommand{\elg}{\end{lemg}}
\newcommand{\blh}{\begin{lemh}}
\newcommand{\elh}{\end{lemh}}
\newcommand{\bli}{\begin{lemi}}
\newcommand{\eli}{\end{lemi}}
\newcommand{\blj}{\begin{lemj}}
\newcommand{\elj}{\end{lemj}}
\newcommand{\bj}{\begin{conjecture}}
\newcommand{\ej}{\end{conjecture}}

\newcommand*{\ldblbrace}{\{\mskip-5mu\{}
\newcommand*{\rdblbrace}{\}\mskip-5mu\}}

\def\diy{\displaystyle}

\def\a{{\alpha}}
\def\b{{\beta}}
\def\Gam{{\Gamma}}
\def\gam{{\gamma}}
\def\del{{\delta}}
\def\eps{{\epsillon}}
\def\veps{{\varepsilon}}
\def\del{{\delta}}
\def\veps {\varepsilon}
\def\vphi{{\varphi}}
\def\vrho{\varrho}
\def\ups{{\upsilon}}
\def\Om {\Omega}

\def\pl {\partial}

\def\rd{{\rm d}} \def\re{{\rm e}} \def\ri{{\rm i}}
\def\rt{{\rm t}} \def\rv{{\rm v}} \def\rw{{\rm w}}
\def\rx{{\rm x}} \def\ry{{\rm y}} \def\rz{{\rm z}}
\def\rD{{\rm D}}
\def\rO{{\rm O}} \def\rP{{\rm P}} \def\rV{{\rm V}}

\def\bD{{\mathbf D}}

\def\ovphi{\ov\vphi}
\def\ophi{\ov\phi}

\def\bbA{{\mathbb A}} \def\bbB{{\mathbb B}} \def\bbC{{\mathbb C}}
\def\bbD{{\mathbb D}} \def\bbE{{\mathbb E}} \def\bbF{{\mathbb F}}
\def\bbG{{\mathbb G}} \def\bbH{{\mathbb H}}
\def\bbL{{\mathbb L}} \def\bbM{{\mathbb M}} \def\bbN{{\mathbb N}}
\def\bbP{{\mathbb P}} \def\bbQ{{\mathbb Q}}
\def\bbR{{\mathbb R}} \def\bbS{{\mathbb S}} \def\bbT{{\mathbb T}}
\def\bbU{{\mathbb U}}
\def\bbV{{\mathbb V}} \def\bbW{{\mathbb W}} \def\bbZ{{\mathbb Z}}

\def\bbZd{\bbZ^d}

\def\cA{{\mathcal A}} \def\cB{{\mathcal B}} \def\cC{{\mathcal C}}
\def\cD{{\mathcal D}} \def\cE{{\mathcal E}} \def\cF{{\mathcal F}}
\def\cG{{\mathcal G}} \def\cH{{\mathcal H}} \def\cJ{{\mathcal J}}
\def\cP{{\mathcal P}} \def\cT{{\mathcal T}} \def\cW{{\mathcal W}}
\def\cX{{\mathcal X}} \def\cY{{\mathcal Y}} \def\cZ{{\mathcal Z}}

\def\tA{{\tt A}} \def\tB{{\tt B}}
\def\tC{{\tt C}} \def\tF{{\tt F}} \def\tH{{\tt H}}
\def\tO{{\tt O}}
\def\tP{{\tt P}} \def\tQ{{\tt Q}} \def\tR{{\tt R}}
\def\tS{{\tt S}} \def\tT{{\tt T}}

\def\tx{{\tt x}} \def\ty{{\tt y}}

\def\t0{{\tt 0}} \def\t1{{\tt 1}}

\def\be{{\mathbf e}} \def\bh{{\mathbf h}}
\def\bn{{\mathbf n}} \def\bu{{\mathbf u}}
\def\bx{{x}} \def\by{{y}} \def\bz{{z}}
\def\B1{{\mathbf 1}} \def\co{\complement}

\def\bmu{{\mbox{\boldmath${\mu}$}}}
\def\bnu{{\mbox{\boldmath${\nu}$}}}
\def\bPhi{{\mbox{\boldmath${\Phi}$}}}

\def\fA{{\mathfrak A}} \def\fB{{\mathfrak B}} \def\fC{{\mathfrak C}}
\def\fD{{\mathfrak D}} \def\fE{{\mathfrak E}} \def\fF{{\mathfrak F}}
\def\fW{{\mathfrak W}} \def\fX{{\mathfrak X}} \def\fY{{\mathfrak Y}}
\def\fZ{{\mathfrak Z}}

\def\rA{{\rm A}}  \def\rB{{\rm B}}  \def\rC{{\rm C}}
\def\rF{{\rm  F}} \def\rM{{\rm  M}}
\def\rS{{\rm S}}
\def\rT{{\rm T}}  \def\rW{{\rm W}}

\def\ov{\overline}  \def\wh{\widehat}  \def\wt{\widetilde}

\def\es {{\varnothing}}

\def\bt {\circ}
\def\cc {\circ}

\def\wt{\widetilde}

\def\wtm{{\wt m}} \def\wtn{{\wt n}} \def\wtk{{\wt k}}

\def\be{\begin{equation}}
\def\ee{\end{equation}}

\def\beal{\begin{array}{l}}
\def\beac{\begin{array}{c}}
\def\bear{\begin{array}{r}}
\def\beacl{\begin{array}{cl}}
\def\beall{\begin{array}{ll}}
\def\bealll{\begin{array}{lll}}
\def\beallll{\begin{array}{llll}}
\def\bealllll{\begin{array}{lllll}}
\def\beacr{\begin{array}{cr}}
\def\ena{\end{array}}

\def\bma{\begin{matrix}}
\def\ema{\end{matrix}}

\def\bpma{\begin{pmatrix}}
\def\epma{\end{pmatrix}}

\def\bcs{\begin{cases}}
\def\ecs{\end{cases}}

\def\diy{\displaystyle}

\def\sA{\mathscr A} \def\sB{\mathscr B} \def\sC{\mathscr C}
\def\sD{\mathscr D} \def\sE{\mathscr E}
\def\sF{\mathscr F} \def\sG{\mathscr G} \def\sI{\mathscr I}
\def\sL{\mathscr L} \def\sM{\mathscr M}
\def\sN{\mathscr N} \def\sO{\mathscr O}
\def\sP{\mathscr P} \def\sR{\mathscr R}
\def\sS{\mathscr S} \def\sS{\mathscr T}
\def\sU{\mathscr U} \def\sV{\mathscr V} \def\sW{\mathscr W}
\def\sX{\mathscr X} \def\sY{\mathscr Y} \def\sZ{\mathscr Z}

\def\BZ{{\mathbf Z}}

\def\D{D}

\def\bs {\overline \phi}
\def\bbLv {{\cal E}}

\def\iy{\infty}
\def\ct{\cdot}
\def\cl{\centerline}

\def\bbL{{\mathbb L}} \def\Pf{{\mathbf Z}} \def\f{{\vphi}} \def\g{{\Gam}} \def\s{{\phi}}
\def\boeta{{\mbox{\boldmath$\eta$}}}
\def\sq{\square} \def\tr{\triangle} \def\lz{\lozenge}

\def\upharp{\hskip-3pt\upharpoonright}
%%%%%%%%%%%%%%%%%%%%%%%%%%%%%%%%%%%%%
%%%%%%%%%%%%%%%%%%%%%%%%%%%%%%%%%%%%%

\makeatletter
 \def\fps@figure{htbp}
\makeatother

%%%%%%%%%%%%%%%%%%%%%%%%%%%%%%%%%%%%%%

\title{\bf % High-density hard-core model on a square lattice
High-density hard-core model on $\bbZ^2$\\
and norm equations in ring $\bbZ [{\sqrt[6]{-1}}]$}  % quadratic rings}

\author{\bf A. Mazel$^1$, I. Stuhl$^2$, Y. Suhov$^{2,3}$}

\date{}
\footnotetext{2010 {\em Mathematics Subject Classification:\; primary 60A10,
60B05, 60C05, 60J20, secondary 68P20, 68P30, 94A17}}
\footnotetext{{\em Key words and phrases:} square lattice, hard-core exclusion
distance, admissible hard-core configuration, Gibbs measure, high density,
dense-packing, periodic ground state, phase transition, contour representation of the partition
function, Peierls bound, dominance, Pirogov-Sinai theory, Zahradnik's argument, norm equation, quadratic integer ring, coset by the unit group, leading solution

\noindent
$^1$ \footnotesize{AMC Health, New York, NY, USA;} $^2$ \footnotesize{Math Department, Penn State University, PA, USA;}
$^3$ \footnotesize{DPMMS, University of Cambridge and St John's College, Cambridge, UK.}}

\maketitle

\begin{abstract}
We study the Gibbs statistics of high-density hard-core configurations on a unit
square lattice $\bbZ^2$, for a general Euclidean exclusion distance $D$. %, {\color{blue}with $D^2=a^2+b^2$,
%$a,b\in\bbZ$ do we need it here???}. %This work is a continuation of \cite{MSS1} where the hard-core model 
%was solved for large fugacities on a unit triangular lattice $\bbA_2$. 
As a by-product, we solve the disk-packing problem on $\bbZ^2$ for disks of diameter $D$. To this end we 
exploit connections between the structure of the ground states and the algebraic number theory. 
Here the key point is an analysis of solutions to norm equations in the ring
$\bbZ[{\sqrt[6]{-1}}]$. We use
Delaunay triangulations to describe the {\it ground state} configurations in terms of {\it M-triangles}, i.e., non-obtuse
$\bbZ^2$-triangles of a minimal area with the side-lengths 
$\geq D$. Further, we identify $\bbZ^2$-triangles and elements of $\bbZ[{\sqrt[6]{-1}}]$; 
this yields a convenient classification of values $D^2$ and the corresponding periodic 
ground states. 

First, there is a {\it finite} class (Class S) formed by values $D^2$ generating
{\it sliding}, a phenomenon leading to countable families of periodic ground states. Following an idea from \cite{K}, we identify all $D^2$ with sliding. Each of the remaining classes is proven to be {\it infinite}; they are characterized by uniqueness 
or non-uniqueness of a minimal triangle for a given $D^2$, up to $\bbZ^2$-congruencies. 
For values of $D^2$ with uniqueness (Class A) we describe the periodic ground states as
admissible sub-lattices in $\bbZ^2$ of maximum density. By using the 
{\it Pirogov-Sinai} theory \cite{PiS}, \cite{Si}, it allows us
to identify the {\it extreme Gibbs measures} (pure phases) for large values of fugacity
and describe symmetries between them. For this purpose we establish a non-trivial 
{\it Peierls bound} based on an appropriate definition of a {\it contour}.
This yields a complete diagram for the phase coexistence at large fugacities. 

Furthermore, we analyze the values $D^2$ with non-uniqueness. 
For some $D^2$ all M-triangles are  ${\bbR}^2$-congruent
but not $\bbZ^2$-congruent (Class B0). For other values of $D^2$
there exist non-${\bbR}^2$-congruent M-triangles,  with different collections 
of side-lengths (Class B1). Moreover, there are values $D^2$ for which both cases occur 
(Class B2). The large-fugacity phase diagram for Classes B0, B1, B2 is determined by {\it dominant} ground states.  

Thus, for all classes but Class S, we show the non-uniqueness of Gibbs measures and hence, establish the existence of a phase transition. Algebraically, all classes A, B0-B2 are 
described in terms of {\it cosets} in $\bbZ[{\sqrt[6]{-1}}]$ by the group of units.
\end{abstract}

\section{The hard-core model on $\bbZ^2$. Gibbs measures}\label{Sec1}

\subsection{Introduction. A summary of the main results.}\label{SubSec1.1} High-density
hard-core (HC) configurations in a continuous space or on a discrete structure (lattices,
(random) graphs, etc.) are important in a number of theoretical and applied
areas, which includes Mathematical and Theoretical Physics, Computer/Information
Sciences, Theoretical Biology etc. The problem is to study the statistical mechanics of
configurations formed by non-overlapping hard spheres of a given dimension and diameter (or
exclusion distance) $D$.
A survey of different aspects of the hard-core model and its applications can be found, e.g.,
in \cite{CKPUZ}, \cite{PaZ}, \cite{PeS}, \cite{BKZZ}, \cite{KMRTSZ}, \cite{JaL}.

High-density hard-core models are intrinsically related to the sphere-packing problem.
In this regard, one can mention the Thue theorem in $\bbR^2$, the Kepler conjecture in $\bbR^3$ and recent results in $\bbR^8$, $\bbR^{24}$. Cf. \cite{FT, Hs, Ha, Vy, CKMRV}.
Randomized hard-core configurations are interesting as they may describe a non-ideal
situation where dense-packing can be occasionally frustrated. The corresponding
mathematical model is provided by Gibbs probability measures in a high
density/large fugacity regime. In a continuous space $\bbR^d$, the analysis of such measures
for $d\geq 2$ is a well-known challenging open problem. 

We do not think that a discrete version of the model in a high-density/large-fugacity regime is an appropriate approximation for the continuous model. It turns out that a discrete HC model has a variety of properties which cannot be expected in the continuous case. Nevertheless, discrete models are interesting, and in some aspects more challenging than the continuous ones. In particular, discrete models depend heavily on the choice of the underlying lattice and the exclusion distance. 

On a unit cubic lattice $\bbZ^d$, initial results for the HC model have been obtained in 1966-8 by
Dobrushin
(and in a part by Suhov); cf. \cite{Dob, Suhov}. Namely, non-uniqueness of a Gibbs
measure was established for large fugacities for $D=\sqrt{2}$. A detailed study of the
high-density/large-fugacity hard-core model in a unit triangular lattice $\bbA_2$ and a unit honeycomb graph $\mathbb{H}_{2}$ for a general $D$ is presented in \cite{MSS1}. In \cite{MSSz3} a number of results have been established about the HC model on a unit cubic lattice $\bbZ^3$ for $D^2=2, 3, 4, 5, 6, 8, 9, 10, 11, 12, 2\ell^2, \ell\in\bbN$ in the large fugacity regime. Short expositions of these results are provided in \cite{MSSshort} and \cite{MSSz3short}.

H-C models on two-dimensional lattices gained a growing popularity in the recent Physics 
literature: see \cite{AGMS}, \cite{AMMGPS}, \cite{DKBP}, \cite{FAL}, \cite{JMTR}, \cite{JTMSR}, \cite{NR1}, \cite{NR2},  
\cite{TF} and references therein. Most of these papers specify the model as $k$-NN exclusion, indicating that a particle 
excludes all $1$-st to $k$-th nearest-neighboring sites (all in the planar Euclidean metric). The value $D$ used to specify the model in the current paper 
is the Euclidean distance to the $(k+1)$-st nearest-neighbors. Depending on the lattice type and the value of $D$, one expects different forms of phase transitions when the particle 
density/fugacity/chemical potential in the system varies. Simulations and analytical 
techniques have been applied to make predictions for critical points and types of the related phase transitions for some initial 
values of $k$. For HC models on $\bbZ^2$ we refer to \cite{FAL}, \cite{NR1}, \cite{NR2},
on $\bbA_2$ to \cite{AGMS}, \cite{JMTR}, and on $\bbH_2$ to \cite{DKBP}, \cite{TF}. However, a rigorous analysis of criticality and related phase transitions remains an open (and challenging) mathematical problem.

In the present paper we focus on random HC configurations in
$\bbZ^2$ distributed according to a Gibbs measure in a large-fugacity regime.
The HC exclusion distance $D$ is imposed in the Euclidean metric $\rho$ on $\bbR^2$. After the work \cite{Dob} the further
progress  has been blocked, in particular, because of the so-called {\it sliding} phenomenon in dense-packing configurations 
first observed by Dobrushin for the value $D=2$. See Figure \ref{Fig1}.
\begin{figure}
\centering
A\begin{tikzpicture}[scale=0.5]%\label{1FParaD^2=4}
\clip (1.6, -3.2) rectangle (8.6, 3.2);

\filldraw[black] (0, 0) circle (.2); \filldraw[black] (0, 2) circle (.2); \filldraw[black] (0, -2) circle (.2);
\filldraw[black] (2, 0) circle (.2); \filldraw[black] (2, 2) circle (.2); \filldraw[black] (2, -2) circle (.2);
\filldraw[black] (4, 0) circle (.2); \filldraw[black] (4, 2) circle (.2); \filldraw[black] (4, -2) circle (.2);
\filldraw[black] (6, 0) circle (.2); \filldraw[black] (6, 2) circle (.2); \filldraw[black] (6, -2) circle (.2);
\filldraw[black] (8, 0) circle (.2); \filldraw[black] (8, 2) circle (.2); \filldraw[black] (8, -2) circle (.2);

%\clip[yscale=sqrt(3/4), xslant=0.5] (0, -2) rectangle (5, 1);

\draw[yscale=1, xslant=0] (-90,-60) grid (90, 60);
\draw[yscale=1, yslant=0] (-90,-60) grid (90, 60);

%\draw (3.3, -0.7) node;   %{\Large {${\mbox{\boldmath$O$}}$}};
%\draw (6.8, 2.7) node;    %{\Large {${\mbox{\boldmath$V_2$}}$}};
%\draw (3.2, 2.7) node;    %{\Large {${\mbox{\boldmath$V_1$}}$}};
%\draw (6.8, -0.7) node;   %{\Large {${\mbox{\boldmath$V_3$}}$}};

\path [draw=black, line width=0.4mm] (4,0) -- (6,2);
\path [draw=black, line width=0.6mm] (4,2) -- (6,2);
\path [draw=black, line width=0.6mm] (4,2) -- (4,0);
\path [draw=black, line width=0.6mm] (4,0) -- (6,0);
\path [draw=black, line width=0.6mm] (6,0) -- (6,2);

\end{tikzpicture}\quad B\begin{tikzpicture}[scale=0.5]%\label{2FParaD62-4}
\clip (1.6, -3.2) rectangle (8.6, 3.2);

\filldraw[black] (0, 0) circle (.2);
\filldraw[black] (0, 2) circle (.2);
\filldraw[black] (0, -2) circle (.2);
\filldraw[black] (2, 3) circle (.2);
\filldraw[black] (2, 1) circle (.2);
\filldraw[black] (2, -1) circle (.2);
\filldraw[black] (2, -3) circle (.2);
\filldraw[black] (4, 0) circle (.2);
\filldraw[black] (4, 2) circle (.2);
\filldraw[black] (4, -2) circle (.2);
\filldraw[black] (6, 3) circle (.2);
\filldraw[black] (6, 1) circle (.2);
\filldraw[black] (6, -1) circle (.2);
\filldraw[black] (6, -3) circle (.2);
\filldraw[black] (8, 0) circle (.2);
\filldraw[black] (8, 2) circle (.2);
\filldraw[black] (8, -2) circle (.2);

%\clip[yscale=sqrt(3/4), xslant=0.5] (0, -2) rectangle (5, 1);

\draw[yscale=1, xslant=0] (-90,-60) grid (90, 60);
\draw[yscale=1, yslant=0] (-90,-60) grid (90, 60);

%\draw (3.3, -0.7) node   {\Large {${\mbox{\boldmath$O$}}$}};
%\draw (6.9, 1.3) node    {\Large {${\mbox{\boldmath$V_2$}}$}};
%\draw (3.3, 2.5) node    {\Large {${\mbox{\boldmath$V_1$}}$}};
%\draw (6.7, -1.5) node   {\Large {${\mbox{\boldmath$V_3$}}$}};

\path [draw=black, line width=0.4mm] (4,0) -- (6,1);
\path [draw=black, line width=0.4mm] (4,2) -- (6,1);
\path [draw=black, line width=0.4mm] (4,0) -- (6,-1);
%\path [draw=black, line width=0.8mm] (4,2) -- (6,2);
\path [draw=black, line width=0.6mm] (4,2) -- (4,0);
%\path [draw=black, line width=0.8mm] (4,0) -- (6,0);
\path [draw=black, line width=0.6mm] (6,-1) -- (6,1);

\end{tikzpicture}\quad C\begin{tikzpicture}[scale=0.5]%\label{1FParaD^2=4}
\clip (2.6, -3.2) rectangle (7.4, 3.2);

%\clip[yscale=sqrt(3/4), xslant=0.5] (0, -2) rectangle (5, 1);

\draw[yscale=1, xslant=0] (-90,-60) grid (90, 60);
\draw[yscale=1, yslant=0] (-90,-60) grid (90, 60);

\path [draw=black, line width=0.4mm] (4,0) -- (6,2);
\path [draw=black, line width=0.6mm] (4,2) -- (6,2);
\path [draw=black, line width=0.6mm] (4,2) -- (4,0);
\path [draw=black, line width=0.4mm] (4,2) -- (6,1);
\path [draw=black, line width=0.4mm] (4,0) -- (6,1);
\path [draw=black, line width=0.4mm] (4,2) -- (6,0);
\path [draw=black, line width=0.6mm] (4,0) -- (6,0);

\draw (3.5,2.5) node   {\Large {${\mbox{$O$}}$}};
\draw (6.7,2.5) node   {\Large {${\mbox{$A$}}$}};
\draw (6.7,0.9) node    {\Large {${\mbox{$B$}}$}};
\draw (6.7,-0.7) node    {\Large {${\mbox{$C$}}$}};
\draw (3.5,-0.7) node    {\Large {${\mbox{$W$}}$}};
\end{tikzpicture}

\vskip 3 truemm
D\begin{tikzpicture}[scale=0.49]%\label{1FParaD^2=5}
\clip (-1.4, -3.4) rectangle (11.4, 4.5);

\filldraw[black] (3, 4) circle (.2);
\filldraw[black] (8, 4) circle (.2);

\filldraw[black] (0, 3) circle (.2);
\filldraw[black] (5, 3) circle (.2);
\filldraw[black] (10, 3) circle (.2);

\filldraw[black] (2, 2) circle (.2);
\filldraw[black] (7, 2) circle (.2);

\filldraw[black] (-1, 1) circle (.2);
\filldraw[black] (4, 1) circle (.2);
\filldraw[black] (9, 1) circle (.2);

\filldraw[black] (1, 0) circle (.2);
\filldraw[black] (6, 0) circle (.2);
\filldraw[black] (11, 0) circle (.2);

\filldraw[black] (3, -1) circle (.2);
\filldraw[black] (8, -1) circle (.2);

\filldraw[black] (0, -2) circle (.2);
\filldraw[black] (5, -2) circle (.2);
\filldraw[black] (10, -2) circle (.2);

\filldraw[black] (2, -3) circle (.2);
\filldraw[black] (7, -3) circle (.2);

\draw[yscale=1, xslant=0] (-90,-60) grid (90, 60);
\draw[yscale=1, yslant=0] (-90,-60) grid (90, 60);

%\draw (0.3, -0.7) node   {\Large {${\mbox{\boldmath$O$}}$}};
%\draw (4.9, 1.3) node    {\Large {${\mbox{\boldmath$V_2$}}$}};
%\draw (2.4, 2.7) node    {\Large {${\mbox{\boldmath$V_1$}}$}};
%\draw (3.7, -1.8) node   {\Large {${\mbox{\boldmath$V_3$}}$}};

\path [draw=black, line width=0.4mm] (1,0) -- (4,1);
\path [draw=black, line width=0.4mm] (4,1) -- (2,2);
\path [draw=black, line width=0.4mm] (1,0) -- (2,2);
\path [draw=black, line width=0.4mm] (1,0) -- (3,-1);
\path [draw=black, line width=0.4mm] (3,-1) -- (4,1);
%\path [draw=black, line width=0.8mm] (4,0) -- (6,0);
%\path [draw=black, line width=0.8mm] (6,-1) -- (6,1);

\end{tikzpicture}

\caption{\footnotesize{Frames A and B show sliding in $\bbZ^2$ for $D=2$: a 1D array of occupied
sites can be $\mathbb{Z}^2$-shifted
without violating the exclusion distance and decreasing the density. Frame C shows
competing fundamental triangles $OAW$, $OBW$ and $OCW$, which generates sliding. Note 
arizing trapezes $OABW$ and $OBCW$ and rectangle $OACW$. In contrast, frame D shows
the absence of sliding for the next attainable exclusion distance $D=\sqrt{5}$.}}
\label{Fig1} %\label{sliding}
\end{figure}
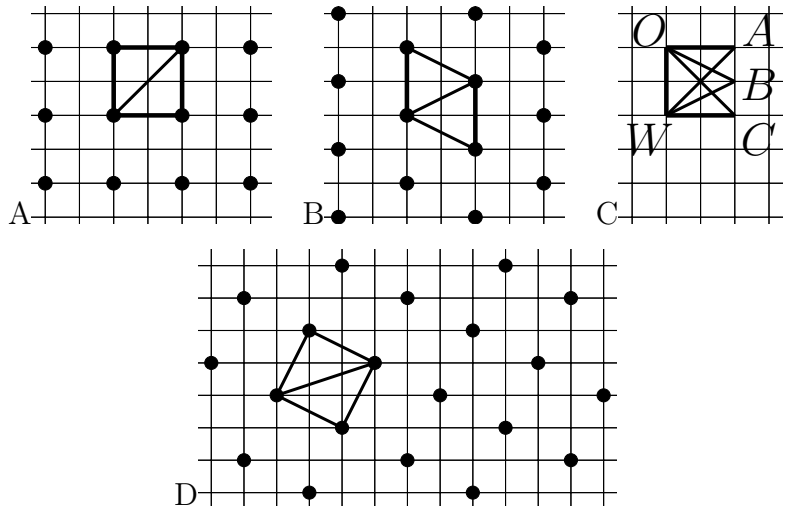
 Here sliding means that the model exhibits a continuum of dense-packing configurations which can be mapped into each other by $\bbZ^2$-shifts of one dimensional arrays of occupied sites. This particular case of sliding has been recently analyzed in \cite{HaPe}, where the existence of a columnar order has been established. In the current paper it is shown that the sliding is a rare phenomenon, and there exist only $39$ values of $D$ with sliding. The behavior of the model in the high-fugacity regime for the $38$ cases of sliding different from $D=2$ remains an open question and in this paper we focus on non-sliding cases of $D$. Our main results can be described as follows.

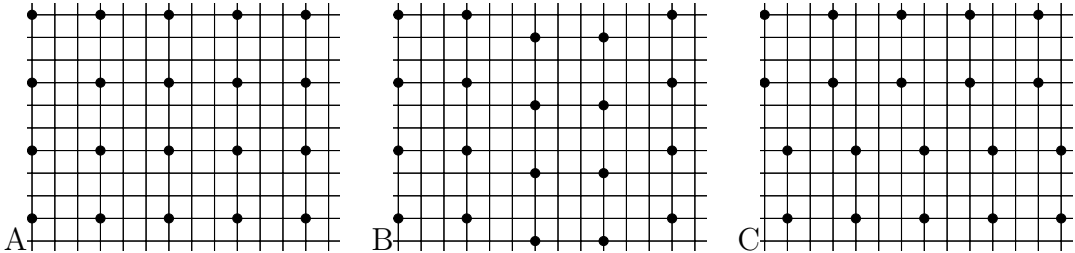
\begin{figure}[H]
\centering
A\begin{tikzpicture}[scale=0.33]%\label{PeriodicLattGSD^2=9}
\clip (-2.2, -3.4) rectangle (11.5, 7.5);

\filldraw[black] (-2, 7) circle (.2);
\filldraw[black] (-2, 4) circle (.2);
\filldraw[black] (-2, 1) circle (.2);
\filldraw[black] (-2, -2) circle (.2);

\filldraw[black] (1, 7) circle (.2);
\filldraw[black] (1, 4) circle (.2);
\filldraw[black] (1, 1) circle (.2);
\filldraw[black] (1, -2) circle (.2);

\filldraw[black] (4, 7) circle (.2);
\filldraw[black] (4, 4) circle (.2);
\filldraw[black] (4, 1) circle (.2);
\filldraw[black] (4, -2) circle (.2);

\filldraw[black] (7, 7) circle (.2);
\filldraw[black] (7, 4) circle (.2);
\filldraw[black] (7, 1) circle (.2);
\filldraw[black] (7, -2) circle (.2);

\filldraw[black] (10, 7) circle (.2);
\filldraw[black] (10, 4) circle (.2);
\filldraw[black] (10, 1) circle (.2);
\filldraw[black] (10, -2) circle (.2);

%\filldraw[black] (4, -1) circle (.2);
%\filldraw[black] (10, -1) circle (.2);

%\filldraw[black] (1, -3) circle (.2);
%%\filldraw[black] (4, -3) circle (.2);
%\filldraw[black] (7, -3) circle (.2);

\draw[yscale=1, xslant=0] (-90,-60) grid (90, 60);
\draw[yscale=1, yslant=0] (-90,-60) grid (90, 60);

%\draw (0.5, -0.7) node    {\Large {${\mbox{\boldmath$O$}}$}};
%\draw (4.5, 2.7) node     {\Large {${\mbox{\boldmath$V_2$}}$}};
%\draw (1, 3.7) node     {\Large {${\mbox{\boldmath$V_1$}}$}};
%\draw (4.4, -1.9) node    {\Large {${\mbox{\boldmath$V_3$}}$}};

%\path [draw=black, line width=0.4mm] (1,0) -- (4,2);
%\path [draw=black, line width=0.4mm] (4,2) -- (1,3);
%\path [draw=black, line width=0.8mm] (1,3) -- (1,0);
%\path [draw=black, line width=0.8mm] (4,2) -- (4,-1);
%\path [draw=black, line width=0.4mm] (1,0) -- (4,-1);
\end{tikzpicture}\quad B\begin{tikzpicture}[scale=0.33]%\label{NonLatticeGSD^2=9}
\clip (-2.2, -3.4) rectangle (11.5, 7.5);

\filldraw[black] (-2, 7) circle (.2);
\filldraw[black] (-2, 4) circle (.2);
\filldraw[black] (-2, 1) circle (.2);
\filldraw[black] (-2, -2) circle (.2);

\filldraw[black] (1, 7) circle (.2);
\filldraw[black] (1, 4) circle (.2);
\filldraw[black] (1, 1) circle (.2);
\filldraw[black] (1, -2) circle (.2);

\filldraw[black] (4, 6) circle (.2);
\filldraw[black] (4, 3) circle (.2);
\filldraw[black] (4, 0) circle (.2);
\filldraw[black] (4, -3) circle (.2);

\filldraw[black] (7, 6) circle (.2);
\filldraw[black] (7, 3) circle (.2);
\filldraw[black] (7, 0) circle (.2);
\filldraw[black] (7, -3) circle (.2);

\filldraw[black] (10, 7) circle (.2);
\filldraw[black] (10, 4) circle (.2);
\filldraw[black] (10, 1) circle (.2);
\filldraw[black] (10, -2) circle (.2);

%\filldraw[black] (4, -1) circle (.2);
%\filldraw[black] (10, -1) circle (.2);

%\filldraw[black] (1, -3) circle (.2);
%%\filldraw[black] (4, -3) circle (.2);
%\filldraw[black] (7, -3) circle (.2);

\draw[yscale=1, xslant=0] (-90,-60) grid (90, 60);
\draw[yscale=1, yslant=0] (-90,-60) grid (90, 60);

%\draw (0.5, -0.7) node    {\Large {${\mbox{\boldmath$O$}}$}};
%\draw (4.5, 2.7) node     {\Large {${\mbox{\boldmath$V_2$}}$}};
%\draw (1, 3.7) node     {\Large {${\mbox{\boldmath$V_1$}}$}};
%\draw (4.4, -1.9) node    {\Large {${\mbox{\boldmath$V_3$}}$}};

%\path [draw=black, line width=0.4mm] (1,0) -- (4,2);
%\path [draw=black, line width=0.4mm] (4,2) -- (1,3);
%\path [draw=black, line width=0.8mm] (1,3) -- (1,0);
%\path [draw=black, line width=0.8mm] (4,2) -- (4,-1);
%\path [draw=black, line width=0.4mm] (1,0) -- (4,-1);
\end{tikzpicture}\quad C\begin{tikzpicture}[scale=0.33]%\label{NonPeriodicGSD^2=9}
\clip (-2.2, -3.4) rectangle (11.5, 7.5);

\filldraw[black] (-2, 7) circle (.2);
\filldraw[black] (-2, 4) circle (.2);
\filldraw[black] (-1, 1) circle (.2);
\filldraw[black] (-1, -2) circle (.2);

\filldraw[black] (1, 7) circle (.2);
\filldraw[black] (1, 4) circle (.2);
\filldraw[black] (2, 1) circle (.2);
\filldraw[black] (2, -2) circle (.2);

\filldraw[black] (4, 7) circle (.2);
\filldraw[black] (4, 4) circle (.2);
\filldraw[black] (5, 1) circle (.2);
\filldraw[black] (5, -2) circle (.2);

\filldraw[black] (7, 7) circle (.2);
\filldraw[black] (7, 4) circle (.2);
\filldraw[black] (8, 1) circle (.2);
\filldraw[black] (8, -2) circle (.2);

\filldraw[black] (10, 7) circle (.2);
\filldraw[black] (10, 4) circle (.2);
\filldraw[black] (11, 1) circle (.2);
\filldraw[black] (11, -2) circle (.2);

%\filldraw[black] (4, -1) circle (.2);
%\filldraw[black] (10, -1) circle (.2);

%\filldraw[black] (1, -3) circle (.2);
%%\filldraw[black] (4, -3) circle (.2);
%\filldraw[black] (7, -3) circle (.2);

\draw[yscale=1, xslant=0] (-90,-60) grid (90, 60);
\draw[yscale=1, yslant=0] (-90,-60) grid (90, 60);

%\draw (0.5, -0.7) node    {\Large {${\mbox{\boldmath$O$}}$}};
%\draw (4.5, 2.7) node     {\Large {${\mbox{\boldmath$V_2$}}$}};
%\draw (1, 3.7) node     {\Large {${\mbox{\boldmath$V_1$}}$}};
%\draw (4.4, -1.9) node    {\Large {${\mbox{\boldmath$V_3$}}$}};

%\path [draw=black, line width=0.4mm] (1,0) -- (4,2);
%\path [draw=black, line width=0.4mm] (4,2) -- (1,3);
%\path [draw=black, line width=0.8mm] (1,3) -- (1,0);
%\path [draw=black, line width=0.8mm] (4,2) -- (4,-1);
%\path [draw=black, line width=0.4mm] (1,0) -- (4,-1);
\end{tikzpicture}
\caption{\footnotesize{%In contrast to a triangular lattice $\bbA_2$,
The hard-core model in $\bbZ^2$ may have lattice (A) and non-lattice (B) periodic
ground states as well as non-periodic ground states (C) of maximal density (shown
for the exclusion distance $D=3$, with sliding). Non-periodic ground states
are not interesting in our context (as follows from a general result by Dobrushin-Shlosman \cite{DoS}).
For a non-sliding $D$ the non-lattice
periodic ground states are ruled out by one of our results (see Lemma I in Section 2.3).}}
\label{Fig2}
\end{figure}

\begin{description}
\item[{\bf{(i)}}] For any non-sliding value of $D$ we prove that
all periodic ground states come from optimal (i.e., max-dense)
admissible {\it sub-lattices} constructed from a minimum-area {\it fundamental
parallelogram} (FP). The latter is determined via a minimizer in a
specific discrete optimization problem; see Eqn (2.1) in Section \ref{SubSec2.1}.

\item[{\bf{(ii)}}] For the values of $D$ of Class A,
where the max-dense sub-lattice is unique up to $\bbZ^2$-symmetries
(in the sense defined later), we use the Pirogov-Sinai theory \cite{PiS,Si} (complemented
by Zahradnik \cite{Z} and Dobrushin-Shlosman \cite{DoS}) and give a complete description
of the extreme Gibbs measures for large fugacities. See Theorem 1 in Section 2.3. In
particular, in Class A the large-fugacity extreme Gibbs measures occur in collections
of cardinality $mS$. Here $m=1$ for $D=1,{\sqrt 2}$; for $D>\sqrt 2$, $m=2$ or $4$,
depending on whether the optimal FP is composed of isosceles or non-isosceles
{\it fundamental triangles} (FTs). (The factor $m$ is related to the aforementioned $\bbZ^2$-symmetries occurring in each of these cases.) Here $S=S(D)$ stands for
the area of the optimal FP (i.e., twice the area of an optimal FT). (See Figures \ref{Fig3}, \ref{Fig4}.)
In this approach, sliding is treated as a specific form of non-uniqueness in the minimization problem Eqn (2.1).

\begin{figure}[H]
\centering

A \begin{tikzpicture}[scale=0.42]\label{NSlidingD2=10}
\clip (-1.6, 2.6) rectangle (13.4, 12.6);

%\filldraw[black] (5,5) circle (.2);

\filldraw[lightgray](0,7)--(1,10)--(4,9)--(3,6)--(0,7);

\filldraw[lightgray] (12,13) circle (1.58);
\filldraw[lightgray] (11,10) circle (1.58);  \filldraw[lightgray] (10,7) circle (1.58);
\filldraw[lightgray] (13,6) circle (1.58); \filldraw[lightgray] (14,9) circle (1.58);
\filldraw[lightgray] (12,3) circle (1.58); \filldraw[lightgray] (9,4) circle (1.58);

\filldraw[black] (2,13) circle (.2); \filldraw[black] (5,12) circle (.2);
\filldraw[black] (8,11) circle (.2); \filldraw[black] (11,10) circle (.2);

\filldraw[black] (1,10) circle (.2); \filldraw[black] (4,9) circle (.2);
\filldraw[black] (7,8) circle (.2);  \filldraw[black] (10,7) circle (.2);
\filldraw[black] (13,6) circle (.2);

\filldraw[black] (0,7) circle (.2); \filldraw[black] (3,6) circle (.2);
\filldraw[black] (6,5) circle (.2); \filldraw[black] (9,4) circle (.2);
\filldraw[black] (12,3) circle (.2); \filldraw[black] (15,2) circle (.2);

\filldraw[black] (-1,4) circle (.2); \filldraw[black] (2,3) circle (.2);
\filldraw[black] (5,2) circle (.2);

\path [draw=black, line width=0.4mm] (0,7) -- (1,10);
\path [draw=black, line width=0.4mm] (0,7) -- (3,6);
\path [draw=black, line width=0.4mm] (0,7) -- (4,9);
\path [draw=black, line width=0.4mm] (1,10) -- (4,9)--(3,6);

\draw[yscale=1, xslant=0] (-90,-60) grid (90, 60);
\draw[yscale=1, yslant=0] (-90,-60) grid (90, 60);
\draw (0,6) node   {\Large {${\mbox{\boldmath$O$}}$}};
\end{tikzpicture}\;\;\;B \begin{tikzpicture}[scale=0.42]\label{NSlidingD2=10}
\clip (-1.6, 2.6) rectangle (13.4, 12.6);

%\filldraw[black] (5,5) circle (.2); \sqrt{10}=3.162

\filldraw[lightgray](0,7)--(1,4)--(4,5)--(3,8)--(0,7);

\filldraw[lightgray] (8,13) circle (1.58); \filldraw[lightgray] (11,14) circle (1.58);
\filldraw[lightgray] (12,11) circle (1.58);  \filldraw[lightgray] (12,11) circle (1.58);
\filldraw[lightgray] (13,8) circle (1.58); \filldraw[lightgray] (9,10) circle (1.58);
\filldraw[lightgray] (10,7) circle (1.58); \filldraw[lightgray] (11,4) circle (1.58);
\filldraw[lightgray] (14,5) circle (1.58); %\filldraw[lightgray] (12,4) circle (2);
%\filldraw[lightgray] (15,7) circle (2);

%\path [draw=gray, line width=1.4mm] (16,9)--(17,12)--(14,13)--(13,10)--(16,9);

\filldraw[black] (-1,10) circle (.2); \filldraw[black] (2,11) circle (.2);
\filldraw[black] (5,12) circle (.2);

\filldraw[black] (0,7) circle (.2); \filldraw[black] (3,8) circle (.2);
\filldraw[black] (6,9) circle (.2); \filldraw[black] (9,10) circle (.2);
\filldraw[black] (12,11) circle (.2);

\filldraw[black] (1,4) circle (.2); \filldraw[black] (4,5) circle (.2);
\filldraw[black] (7,6) circle (.2); \filldraw[black] (10,7) circle (.2);
\filldraw[black] (13,8) circle (.2);

\filldraw[black] (8,3) circle (.2); \filldraw[black] (11,4) circle (.2);

\path [draw=black, line width=0.4mm] (0,7) -- (3,8);
\path [draw=black, line width=0.4mm] (0,7) -- (1,4);
\path [draw=black, line width=0.4mm] (0,7) -- (4,5);
\path [draw=black, line width=0.4mm] (1,4) -- (4,5)--(3,8);

\draw[yscale=1, xslant=0] (-90,-60) grid (90, 60);
\draw[yscale=1, yslant=0] (-90,-60) grid (90, 60);
\draw (0,8) node   {\Large {${\mbox{\boldmath$O$}}$}};
\end{tikzpicture}\vskip 3mm

C\begin{tikzpicture}[scale=0.42]\label{NSlidingD2=16}
\clip (-1.6, 1.6) rectangle (14.4, 13.6);

%\filldraw[black] (5,5) circle (.2);

\filldraw[lightgray](5,12)--(1,11)--(0,7)--(4,8)--(5,12);

\filldraw[lightgray] (14,3) circle (2);
\filldraw[lightgray] (11,6) circle (2); \filldraw[lightgray] (12,10) circle (2);
\filldraw[lightgray] (15,7) circle (2); \filldraw[lightgray] (9,13) circle (2);
\filldraw[lightgray] (13,14) circle (2); %\filldraw[lightgray] (12,4) circle (2);
%\filldraw[lightgray] (15,7) circle (2);

\path [draw=gray, line width=1.4mm] (10,13)--(10,9)--(14,9)--(14,13)--(10,13);

\filldraw[black] (6,1) circle (.2); \filldraw[black] (10,2) circle (.2);
\filldraw[black] (14,3) circle (.2); %F

\filldraw[black] (-1,3) circle (.2);
\filldraw[black] (0,7) circle (.2); \filldraw[black] (4,8) circle (.2);
\filldraw[black] (8,9) circle (.2); \filldraw[black] (12,10) circle (.2);
%\filldraw[black] (16,11) circle (.2); %F

\filldraw[black] (3,4) circle (.2); \filldraw[black] (7,5) circle (.2);
\filldraw[black] (11,6) circle (.2); \filldraw[black] (15,7) circle (.2); %F

\filldraw[black] (1,11) circle (.2);
\filldraw[black] (5,12) circle (.2); \filldraw[black] (9,13) circle (.2);

\path [draw=black, line width=0.4mm] (0,7) -- (1,11);
\path [draw=black, line width=0.4mm] (0,7) -- (4,8);
\path [draw=black, line width=0.4mm] (1,11) -- (4,8)--(5,12)--(1,11);

\draw[yscale=1, xslant=0] (-90,-60) grid (90, 60);
\draw[yscale=1, yslant=0] (-90,-60) grid (90, 60);
\draw (0,6) node   {\Large {${\mbox{\boldmath$O$}}$}};
\end{tikzpicture}\;D \begin{tikzpicture}[scale=0.42]\label{2NSlidingD2=16}
\clip (-1.6, 1.6) rectangle (14.4, 13.6);

%\filldraw[black] (5,5) circle (.2);

\filldraw[lightgray](1,3)--(5,2)--(4,6)--(0,7)--(1,3);

\filldraw[lightgray] (10,12) circle (2); \filldraw[lightgray] (14,11) circle (2);
\filldraw[lightgray] (7,9) circle (2); \filldraw[lightgray] (11,8) circle (2);
\filldraw[lightgray] (8,5) circle (2); \filldraw[lightgray] (12,4) circle (2);
\filldraw[lightgray] (15,7) circle (2);

\path [draw=gray, line width=1.4mm] (9,13)--(9,9)--(13,9)--(13,13)--(9,13);

\filldraw[black] (6,13) circle (.2); \filldraw[black] (10,12) circle (.2);
\filldraw[black] (14,11) circle (.2);

\filldraw[black] (-1,11) circle (.2);
\filldraw[black] (0,7) circle (.2); \filldraw[black] (4,6) circle (.2);
\filldraw[black] (8,5) circle (.2); \filldraw[black] (12,4) circle (.2);
\filldraw[black] (16,3) circle (.2);

\filldraw[black] (3,10) circle (.2); \filldraw[black] (7,9) circle (.2);
\filldraw[black] (11,8) circle (.2); \filldraw[black] (15,7) circle (.2);

\filldraw[black] (1,3) circle (.2);
\filldraw[black] (1,3) circle (.2);
\filldraw[black] (5,2) circle (.2); \filldraw[black] (9,1) circle (.2);
\filldraw[black] (13,0) circle (.2);

\path [draw=black, line width=0.4mm] (1,3) -- (4,6);
\path [draw=black, line width=0.4mm] (1,3) -- (0,7);
\path [draw=black, line width=0.4mm] (1,3) -- (5,2);
\path [draw=black, line width=0.4mm] (0,7) -- (4,6);
\path [draw=black, line width=0.4mm] (5,2) -- (4,6);

\draw[yscale=1, xslant=0] (-90,-60) grid (90, 60);
\draw[yscale=1, yslant=0] (-90,-60) grid (90, 60);
\draw (0,8) node   {\Large {${\mbox{\boldmath$O$}}$}};
\end{tikzpicture}

\caption{\footnotesize{%In contrast to a triangular lattice $\bbA_2$,
Counting the extreme Gibbs measures for exclusion distances $D=\sqrt{10}$ (A, B)
and $D=4,{\sqrt{17}}$ (C,D), of Class A. The minimal area of the fundamental parallelograms (light-gray)
is $S=10$ and $S=15$, and  the fundamental triangles are isosceles, with
side-lengths and $\sqrt{10}$, $\sqrt{10}$, $\sqrt{20}$ and $\sqrt{17}$, $\sqrt{17}$,
$\sqrt{18}$, respectively. In each case, there are two max-dense sub-lattices
obtained from each other by ${\mathbb Z}^2$-symmetries:
$\pm\frac{\pi}{2}$-rotations or reflections. Each sub-lattice
contributes $10$ and $15$ periodic ground states obtained by ${\mathbb Z}^2$-shifts
along the fundamental parallelogram; this explains why the number of periodic
ground states  equals $20$ for $D=\sqrt{10}$ and $30$ for $D=4$. Our Theorem 1
states that the number of extreme Gibbs measures for large fugacities is exactly $20$ and
$30$, respectively. For $D=4$, the $4\times 4$ squares (dark gray) would lead to a
non-optimal sub-lattice.}}
\label{Fig3}
\end{figure}
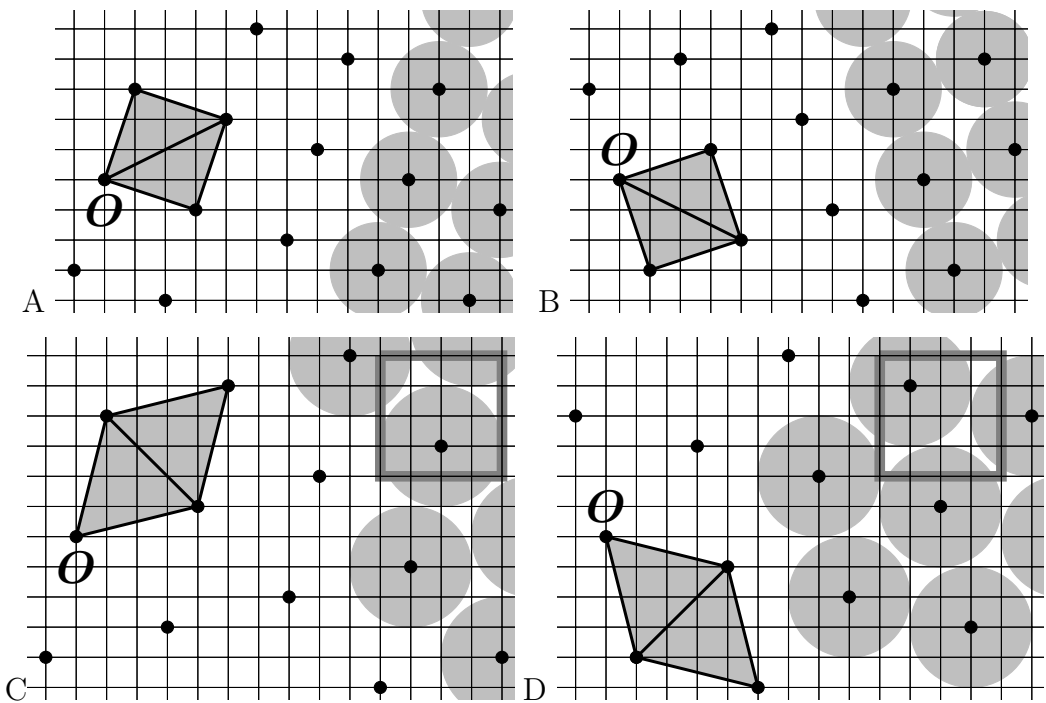

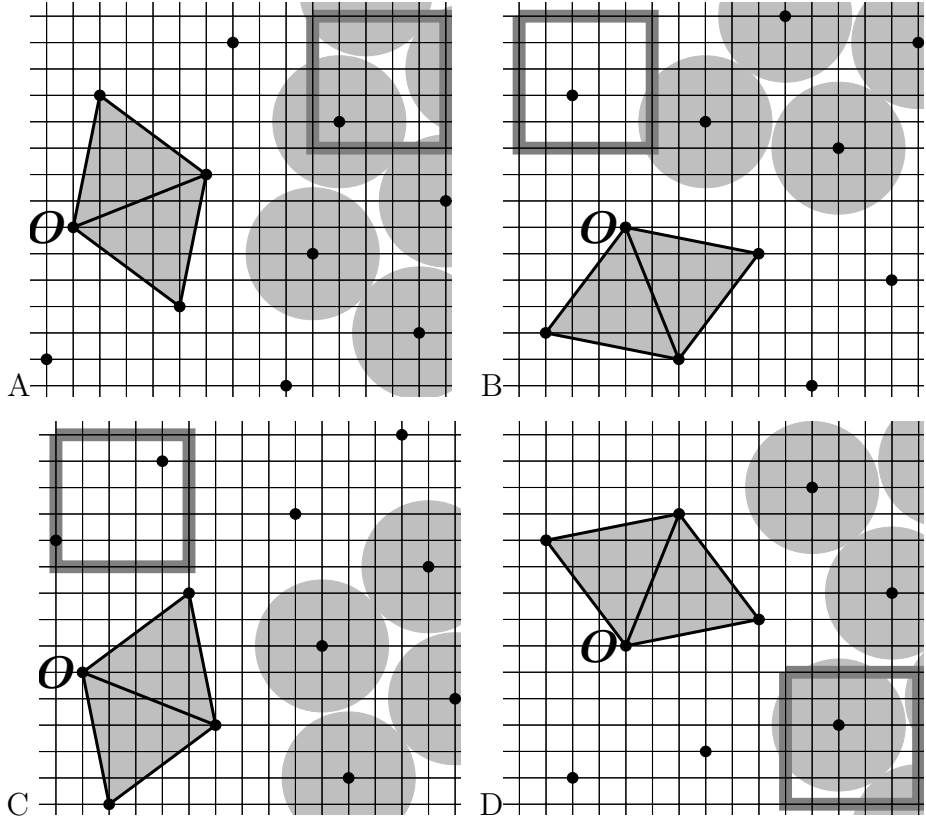
\begin{figure}
\centering

A\begin{tikzpicture}[scale=0.35]\label{1FParaD^2=8}
\clip (-1.6, -9.4) rectangle (14.2, 5.5);

\filldraw[lightgray] (0,-3)--(5,-1)--(1,2)--(0,-3);
\filldraw[lightgray] (0,-3)--(4,-6)--(5,-1)--(0,-3);

\filldraw[lightgray] (10,1) circle (2.5); \filldraw[lightgray] (9,-4) circle (2.5);
\filldraw[lightgray] (14,-2) circle (2.5); \filldraw[lightgray] (13,-7) circle (2.5);
\filldraw[lightgray] (11,6) circle (2.5);  \filldraw[lightgray] (15,3) circle (2.5);

\path [draw=gray, line width=1.7mm] (9,5)--(14,5)--(14,0)--(9,0)--(9,5);

\path [draw=black, line width=0.4mm] (0,-3) -- (1,2);
\path [draw=black, line width=0.4mm] (0,-3) -- (5,-1);
\path [draw=black, line width=0.4mm] (0,-3) -- (4,-6);
\path [draw=black, line width=0.4mm] (1,2) -- (5,-1)--(4,-6);

\filldraw[black] (-3,5) circle (.2);
\filldraw[black] (0,-3) circle (.2);
\filldraw[black] (5,-1) circle (.2);  \filldraw[black] (10,1) circle (.2);
\filldraw[black] (1,2) circle (.2);   \filldraw[black] (6,4) circle (.2);
\filldraw[black] (-4,0) circle (.2);
\filldraw[black] (-5,-5) circle (.2); \filldraw[black] (14,-2) circle (.2);
\filldraw[black] (-1,-8) circle (.2);
\filldraw[black] (4,-6) circle (.2); \filldraw[black] (9, -4) circle (.2);
\filldraw[black] (8,-9) circle (.2); \filldraw[black] (13,-7) circle (.2);

\draw[yscale=1, xslant=0] (-90,-60) grid (90, 60);
\draw[yscale=1, yslant=0] (-90,-60) grid (90, 60);

\draw (-1,-3) node   {\Large {${\mbox{\boldmath$O$}}$}};
%\draw (1,3) node   {\Large {${\mbox{\boldmath$V_1$}}$}};
%\draw (6,-1) node   {\Large {${\mbox{\boldmath$V_2$}}$}};
%\draw (4,-7) node   {\Large {${\mbox{\boldmath$V_3$}}$}};
\end{tikzpicture} \; B\begin{tikzpicture}[scale=0.35]\label{1FParaD^2=8}
\clip (-4.6, -9.4) rectangle (11.2, 5.5);

\filldraw[lightgray] (0,-3)--(-3,-7)--(2,-8)--(5,-4)--(0,-3);

\filldraw[lightgray] (6,5) circle (2.5);
\filldraw[lightgray] (11,4) circle (2.5);
\filldraw[lightgray] (8,0) circle (2.5);
\filldraw[lightgray] (3,1) circle (2.5);
%\path [draw=lightgray, line width=1.7mm] (6,3)--(11,3)--(11,-2)--(6,-2)--(6,3);

\path [draw=gray, line width=1.7mm] (-4,5)--(1,5)--(1,0)--(-4,0)--(-4,5);

\path [draw=black, line width=0.4mm] (0,-3) -- (-3,-7);
\path [draw=black, line width=0.4mm] (0,-3) -- (5,-4);
\path [draw=black, line width=0.4mm] (0,-3) -- (2,-8);
\path [draw=black, line width=0.4mm] (-3,-7) -- (2,-8)--(5,-4);

\filldraw[black] (6,5) circle (.2);
\filldraw[black] (11,4) circle (.2);

\filldraw[black] (-2,2) circle (.2);
\filldraw[black] (3,1) circle (.2); \filldraw[black] (8,0) circle (.2);

\filldraw[black] (0,-3) circle (.2);
\filldraw[black] (5,-4) circle (.2);  \filldraw[black] (10,-5) circle (.2);

\filldraw[black] (-3,-7) circle (.2);
\filldraw[black] (2,-8) circle (.2); \filldraw[black] (7,-9) circle (.2);

\draw[yscale=1, xslant=0] (-90,-60) grid (90, 60);
\draw[yscale=1, yslant=0] (-90,-60) grid (90, 60);

\draw (-1,-3) node   {\Large {${\mbox{\boldmath$O$}}$}};
%\draw (1,3) node   {\Large {${\mbox{\boldmath$V_1$}}$}};
%\draw (6,-1) node   {\Large {${\mbox{\boldmath$V_2$}}$}};
%\draw (4,-7) node   {\Large {${\mbox{\boldmath$V_3$}}$}};
\end{tikzpicture}
\vskip 3mm

C \begin{tikzpicture}[scale=0.35]\label{1FParaD^2=8}
\clip (-1.6, -8.4) rectangle (14.2, 6.5);

\filldraw[lightgray] (14,-4) circle (2.5);
\filldraw[lightgray] (13,1) circle (2.5);
\filldraw[lightgray] (10,-7) circle (2.5);
\filldraw[lightgray] (9,-2) circle (2.5);

\path [draw=gray, line width=1.7mm] (-1,6)--(4,6)--(4,1)--(-1,1)--(-1,6);

\filldraw[lightgray] (0,-3)--(1,-8)--(5,-5)--(4,0)--(0,-3);

\path [draw=black, line width=0.4mm] (0,-3) -- (1,-8);
\path [draw=black, line width=0.4mm] (0,-3) -- (5,-5);
\path [draw=black, line width=0.4mm] (0,-3) -- (4,0);
\path [draw=black, line width=0.4mm] (1,-8) -- (5,-5)--(4,0);

\filldraw[black] (-1,2) circle (.2);
\filldraw[black] (0,-3) circle (.2);
\filldraw[black] (1,-8) circle (.2);

\filldraw[black] (3,5) circle (.2);
\filldraw[black] (4,0) circle (.2);
\filldraw[black] (5,-5) circle (.2);

\filldraw[black] (8,3) circle (.2);
\filldraw[black] (9,-2) circle (.2);
\filldraw[black] (10,-7) circle (.2);

 \filldraw[black] (12, 6) circle (.2);
\filldraw[black] (13,1) circle (.2);
\filldraw[black] (14,-4) circle (.2);

%\filldraw[black] (4,0) circle (.2);
%\filldraw[black] (8,-9) circle (.2);

\draw[yscale=1, xslant=0] (-90,-60) grid (90, 60);
\draw[yscale=1, yslant=0] (-90,-60) grid (90, 60);

\draw (-1,-3) node   {\Large {${\mbox{\boldmath$O$}}$}};
%\draw (1,3) node   {\Large {${\mbox{\boldmath$V_1$}}$}};
%\draw (6,-1) node   {\Large {${\mbox{\boldmath$V_2$}}$}};
%\draw (4,-7) node   {\Large {${\mbox{\boldmath$V_3$}}$}};
\end{tikzpicture}\; D\begin{tikzpicture}[scale=0.35]\label{1FParaD^2=8}
\clip (-4.6, -9.4) rectangle (11.2, 5.5);

\filldraw[lightgray] (13,-5) circle (2.5);
\filldraw[lightgray] (12,4) circle (2.5);
\filldraw[lightgray] (11,-10) circle (2.5);
\filldraw[lightgray] (10,-1) circle (2.5);
\filldraw[lightgray] (8,-6) circle (2.5);
\filldraw[lightgray] (7,3) circle (2.5);

\filldraw[lightgray] (0,-3)--(5,-2)--(2,2)--(-3,1)--(0,-3);

\path [draw=gray, line width=1.7mm] (6,-4)--(11,-4)--(11,-9)--(6,-9)--(6,-4);

\path [draw=black, line width=0.4mm] (0,-3) -- (5,-2);
\path [draw=black, line width=0.4mm] (0,-3) -- (2,2);
\path [draw=black, line width=0.4mm] (0,-3) -- (-3,1);
\path [draw=black, line width=0.4mm] (-3,1)--(2,2)--(5,-2);

\filldraw[black] (-2,-8) circle (.2);

\filldraw[black] (-3,1) circle (.2);
\filldraw[black] (0,-3) circle (.2); \filldraw[black] (3,-7) circle (.2);

\filldraw[black] (2,2) circle (.2);
\filldraw[black] (5,-2) circle (.2); \filldraw[black] (8,-6) circle (.2);

\filldraw[black] (7,3) circle (.2);
\filldraw[black] (10,-1) circle (.2);

%\filldraw[black] (-1,-8) circle (.2);

\draw[yscale=1, xslant=0] (-90,-60) grid (90, 60);
\draw[yscale=1, yslant=0] (-90,-60) grid (90, 60);
\draw (-1,-3) node   {\Large {${\mbox{\boldmath$O$}}$}};
\end{tikzpicture}

\caption{\footnotesize{%In contrast to a triangular lattice $\bbA_2$,
Counting the extreme Gibbs measures for exclusion distance $D=5$ of Class A.
The minimal area of the fundamental parallelograms (light-gray)
is $S=23$. The corresponding fundamental triangles are non-isosceles, with
side-lengths $5$, $\sqrt{26}$ and $\sqrt{29}$. There are four max-dense sub-lattices A--C
obtained from each other by ${\mathbb Z}^2$-symmetries:
$\pm\frac{\pi}{2}$-rotations (A, B and C, D) or reflections (A, C and B, D).
Each sub-lattice contribute $23$ periodic ground states obtained by ${\mathbb Z}^2$-shifts
along the fundamental parallelogram;
this explains why the number of periodic ground states for $D=5$ equals $4S=92$. Our Theorem 1
is that the number of extreme Gibbs measures for large fugacities is exactly $92$. The $5\times 5$ squares would lead to a (non-optimal) sub-lattice. \\
$\qquad$ Observe that Figures 3 and 4 demonstrate partial touching or absence of touching of
disks in ground states.}}
\label{Fig4}
\end{figure}

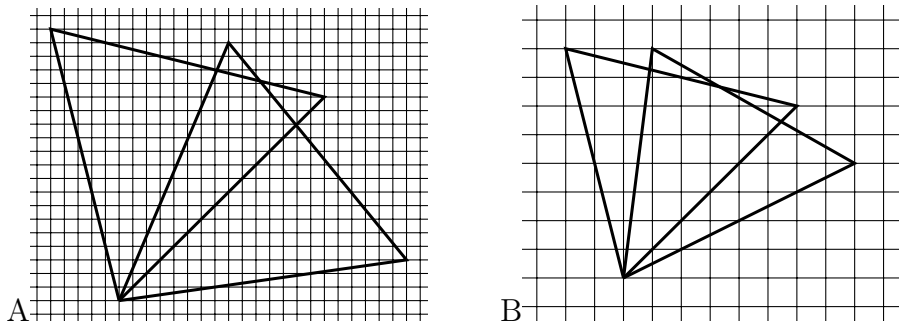
\begin{figure}
\centering
A\begin{tikzpicture}[scale=0.18]
\clip (-6.5, -1.5) rectangle (22.5, 21.5);

%\clip[yscale=sqrt(3/4), xslant=0.5] (0, -2) rectangle (5, 1);

\draw[yscale=1, xscale=1] (-100,-100) grid (100, 100);
%\draw[yscale=1, yslant=0] (-90,-60) grid (90, 60);

\path [draw=black, line width=0.4mm] (0,0) -- (15,15) -- (-5,20) -- (0,0);
\path [draw=black, line width=0.4mm] (0,0) -- (21,3) -- (8,19) -- (0,0);

\end{tikzpicture}
\qquad B\begin{tikzpicture}[scale=0.38]
\clip (-3.5, -1.5) rectangle (9.5, 9.5);

%\clip[yscale=sqrt(3/4), xslant=0.5] (0, -2) rectangle (5, 1);

\draw[yscale=1, xscale=1] (-100,-100) grid (100, 100);
%\draw[yscale=1, yslant=0] (-90,-60) grid (90, 60);

\path [draw=black, line width=0.4mm] (0,0) -- (8,4) -- (1,8) -- (0,0);
\path [draw=black, line width=0.4mm] (0,0) -- (6,6) -- (-2,8) -- (0,0);

\end{tikzpicture}

\caption{\footnotesize{Examples of non-uniqueness of optimal fundamental
triangles: $D^2=425$, $S=375$ (Frame A) and $D^2=65$, $S=60$ (Frame B). The
optimal squared side lengths for the case  $D^2=425$ are $425,425,450$, with
two non-equivalent $\mathbb{Z}^2$ implementations, which falls in Class B0.
The optimal
squared side-lengths for $D^2=65$ are $65, 65, 80$ (right) and $68, 68, 72$ (left);
both triangles admit a unique implementation up to $\bbZ^2$-symmetries.
This case belongs to Class B1, and to determine the number of extreme Gibbs measures
one has to use the dominance formalism in full. Cf. Section 4.2.}}
\label{Fig5}
\end{figure}

\item[{\bf{(iii)}}] For the values of $D$ with non-uniqueness but without sliding (Class B), we
offer in our Theorem 2, Section \ref{SubSec2.3}, a less explicit result, that
extreme high-density Gibbs measures come from a particular `dominant'/`stable' group (or
groups) of optimal FPs. Cf. \cite{BrS, Z}. An identification of a dominant group requires a
separate analysis and is outside the scope of this paper. However, we provide a guidance
in this direction, via the so-called {\it dominance analysis} based on counting {\it local
excitations} of periodic ground states. Such a technique was introduced in \cite{MSS1} for the HC
model in $\mathbb{A}_2$.
  \item[{\bf{(iv)}}] The uniqueness and non-uniqueness of an optimal FP and FT can be examined
both analytically and numerically. We prove that uniqueness (Class A) occurs for countably many
exclusion distances by constructing natural infinite sequences of $D$ for which uniqueness of
an optimal FP/FT is valid. The construction has a number-theoretical character and is based upon
{\it norm equations} in the ring $\bbZ[{\sqrt[6]{-1}}]$. (A particular case here is the famous Pell equation.)
  \item[{\bf{(v)}}] Next, we establish that Class B is also infinite. Moreover, we describe all
possible patterns of non-uniqueness which can occur along natural infinite sequences of values of $D$. 
As before, it follows from the study of norm equations in $\bbZ[{\sqrt[6]{-1}}]$.

\item[{\bf{(vi)}}] For all non-sliding values of $D$ we show the existence of a phase transition due to the non-uniqueness of extreme Gibbs measures \cite{Ge}.

  \item[{\bf{(vii)}}] Finally, by developing the construction from paper \cite{K}, we identify all values $D^2$ 
with sliding. See Theorem S in Section 2.2.
\end{description}

Therefore, we provide a complete structure of the phase diagram
of the HC model on $\bbZ^2$ for large values of fugacity and for all non-sliding exclusion
distances. In a sense, it gives a full description of
phenomena that can occur in this model in the high-fugacity regime and in the absence of sliding. %Of course, a number of important questions
%remain open; see comments in the main body of the text.

% analyze numerically min-area FTs for all values $D\leq 2000$
% (i.e., $D^2\leq 4\cdot 10^6$). The obtained numerical results suggest that 
% the value $D^2=37970$ is the largest one with sliding, and all values of $D^2$ with sliding are revealed 
% already by our numeric calculations, see 

In \cite{MSS1} we obtained a number of similar results about the HC model
on a unit triangular lattice $\bbA_2$ and the honeycomb graph $\bbH_2$; the comparisons between the two models will be instructive for understanding the present work and show
lattice-dependence of the obtained results. We conduct such comparisons in Remarks 1.1--4.1 located at appropriate places throughout the text.

%{\bf Remark 1.1.}
\br\label{rem1.1} {\rm As was shown in \cite{MSS1},
the periodic max-density HC configurations in $\bbA_2$ are produced from equilateral
triangular sub-lattices and can be treated as direct analogs of a max-dense
disk-packing in the Euclidean plane $\bbR^2$. The situation
in $\bbZ^2$ is more complex since the max-dense sub-lattices are formed by
non-equilateral triangles. Furthermore, the max-dense sub-lattices for a given $D$ are
not necessarily congruent. Consequently, the results of group  {\bf (i)} (see
above) require a technically involved proof.}  \quad $\blacktriangle$
\er

\subsection{Hard-core Gibbs/DLR measures.}\label{SubSec1.2} As was mentioned, the HC
exclusion is imposed in the Euclidean metric $\rho$ and is described via the HC exclusion
distance $D>0$ (aka the HC {\it  diameter}).
We will suppose that the value $D$ is {\it attainable},  i.e., $D^2$ is a
{\it Gaussian} number admitting the
{\it Fermat representation} $D^2 = m^2 + n^2$, where  $m$ and $n$ are
integers. Physically, attainability means that
there exists a pair of sites $x,y\in\bbZ^2$ with $\rho (x,y)=D$.

For convenience, an initial list of attainable values of $\D$ is given in Table 1.
The assumption of attainability is justified, since if $D_1<D<D_2$ where $D_{1,2}$
are subsequent Gaussian numbers then the HC model with the exclusion distance $D$ is
reduced to that with distance $D_2$.

\begin{table}[ht]
\centering
\boxed{\begin{matrix}%\label{Gaussians}
1,2,4,5,8,9,10,13,16,17,18,20,25,26,29,32,34,36,37,40,41,
45,49,50,52,53,58,61,64,65,\\68,72,73,74,80,81,82,85,89,90,97,98,100,101,104,106,109,113,116,117,121,122,125,\\
128,130,136,137,144,145,146,148,149,153,157,160,162,164,169,170,173,178,180,181,\\
185,193,194,196,197,200,202,205,208,212,218,221,225,226,229,232,233,234,241,242,\\
244,245,250,256,257,260,261,265,269,272,274,277,281,288,289,290,292,293,296,298\\
\end{matrix}}
%\vskip .5cm
\caption{\footnotesize{Gaussian numbers $D^2\leq 300$.}}\label{Gaussians}
\end{table}

%The list of Gaussian integers is widely available in the literature and on the Internet.
%Cf. http://nonagon.org/ExLibris/fermat-sum-two-squares-calculator.
%In Table 1 below we provide attainable values with $\D^2\leq 300$.

We work with $D$-{\it admissible}
configurations ($D$-ACs) in $\bbZ^2$, represented
by maps $\phi:\bbZ^2\to\{0,1\}$ such that $\rho (x,y)\geq D$ for every
pair of distinct sites $x,y\in\bbZ^2$ with $\phi (x)=\phi (y)=1$. A site $x$ with
$\phi (x)=1$ is interpreted as occupied by a `particle' while site $y$ with $\phi (y)=0$
as vacant. Particles are treated as disjoint open disks of diameter $D$ with
the centers placed at lattice sites. We write $x\in\phi$ if $\phi (x)=1$ and identify $\phi$
with the subset in $\bbZ^2$ where $\phi (x)=1$. 

The collection $\cA=\cA(D)$ of $D$-ACs forms a closed subset in the Cartesian
product $\sX:=\{0,1\}^{\bbZ^2}$ in the
Tykhonov topology. 
%The notion of an $D$-AC can be given for a subset $\bbW\subset{\bbZ^2}$. For $D=1$, $\cA =\sX$.
The notion of an admissible configuration can be defined for any $\bbV\subset\bbZ^2$;
accordingly, one can use the notation $\cA (\bbV)=\cA (D,\bbV)$.

We analyze probability measures $\bmu$ on $(\sX,\fB(\sX))$ sitting on $\cA$
(i.e., such that $\bmu (\cA )=\bmu (\sX)=1$) where $\fB(\sX)$ is the Borel
$\sigma$-algebra in $\sX$. %The measures of interest are identified as {\it extreme} hard-core {\it Gibbs}, or {\it DLR}, probability measures for high densities/large
%fugacities. The precise definitions follow.  We will use the term a $D$-HC Gibbs/DLR
%measure or simply a Gibbs measure when the dependence upon $D$ can be omitted.
Let $\bbV\subset{\bbZ}^2$ be a finite set and $\phi\in{\mathcal A}$ a
$D$-AC. We say that a finite configuration $\psi^\bbV\in\{0,1\}^\bbV$
is $(\phi ,\bbV)$-compatible if the concatenated configuration
$\psi^\bbV\vee (\phi \upharp_{\bbZ^2\setminus\bbV})\in{\mathcal A}$
(which requires that $\psi^\bbV\in\cA (\bbV)$). Here and below, symbol $\upharpoonright$ stands for the restriction. 
The set of $(\phi ,\bbV)$-compatible configurations is denoted by $\cA(\bbV\|\phi )$.

Given $u>0$, consider a probability measure
$\mu_{\bbV}(\;\cdot\; \|\phi )$ on $\{0,1\}^\bbV$ given by
$$ \mu_{\bbV}(\psi^\bbV\,\|\phi )=\begin{cases}\diy\frac{u^{\sharp
(\psi^\bbV)}}{\BZ(\bbV\|\phi )},&\hbox{if $\psi^\bbV\in\cA(\bbV\|\phi )$,}\\
0,&\hbox{if $\psi^\bbV\not\in\cA(\bbV\|\phi )$.}
\end{cases}
\eqno (1.1)$$
Here $\sharp (\psi^\bbV)$ stands for the number of particles in $\psi^\bbV$:
$\sharp (\psi^\bbV):=\#\big\{x\in\bbV:\;\phi (x)=1\big\}$. Next, $\BZ(\bbV\;\|\phi )$
is the {\it partition function} in $\bbV$ with the boundary condition $\phi$\,:
$$\BZ(\bbV\|\phi )=\sum\limits_{\psi^\bbV\in\cA (\bbV\|\phi )}u^{\sharp
(\psi^\bbV)}.\eqno (1.2)$$
Parameter $u>0$ is called {\it fugacity}, or {\it activity} (of an occupied site).

A probability measure $\bmu$ on $(\sX,\fB(\sX))$ is called a $D$-HC Gibbs/DLR measure \ if (i) $\bmu ({\mathcal A})=1$, (ii)
$\forall$ \ finite $\bbV\subset\bbZ^2$ and a function $f:\phi\in
\sX\mapsto f(\phi )\in{\mathbb C}$ depending only on the
restriction $\phi\upharp_\bbV$, the integral $\bmu (f)=\int_\sX f(\phi
)\rd\bmu (\phi )$ has the form
$$\begin{array}{c}\bmu (f)=\diy\int\limits_\sX\int\limits_{\{0,1\}^\bbV}
f(\psi^\bbV\vee\phi\upharp_{\bbZ^2\setminus\bbV})\rd\mu_\bbV(\psi^\bbV\,
\|\phi )\rd\bmu (\phi ).\end{array}\eqno (1.3)$$
One can say that under such measure $\bmu$, the probability of a configuration
$\psi^\bbV$ in a finite `volume' $\bbV\subset{\mathbb Z}^2$, conditional on
a configuration
$\phi\upharp_{\bbZ^2\setminus\bbV}$, coincides with $\mu_{\bbV}(\psi^\bbV\,\|\phi)$,
for $\bmu$-a.a. $\phi\in\{0,1\}^{\bbZ^2}$.

In the literature, equality (1.3) is referred to as the DLR (Dobrushin-Lanford-Ruelle) equation for a measure
$\bmu$ (it represents a system of equations labeled by $\bbV$ and $f$).

The HC Gibbs measures form a {\it Choquet simplex} (in the weak-convergence
topology on the set of probability measures on $(\sX,\fB(\sX))$), which we denote by
$\sG=\sG(D,u)$. We are interested in {\it extreme Gibbs measures}, i.e., extreme points of $\sG$. 
%$\bmu$ which do not admit a non-trivial decomposition $\bmu =\alpha\bmu^{(1)}
%+(1-\alpha )\bmu^{(2)}$ in terms of other HC Gibbs measures $\bmu^{(i)}$, $i=1,2$.
In Physics the extreme Gibbs measures are interpreted as {\it pure phases}.
The collection of the extreme Gibbs measures (EGMs) is denoted by $\sE =\sE(D,u)$ and its cardinality is denoted by $\sharp \sE$. Any HC DLR-measure $\bmu$ is the barycenter for some unit mass distribution over
$\sE$. Cf. \cite{Ge}.

%One of the main questions is to verify, for given $D$ and $u$, whether there is
%a unique or multiple EGM.
%It is known that for $u$ small ($u\in (0,u^0)$ where $u^0=u^0(D)$), the EGM
%(and hence the HC Gibbs measure) $\bmu$ is unique. As was stressed, we
%work with $u$ large: $u\geq u_*=u_*(D)$ where $u_* \geq u^0$. Applying the Pirogov-Sinai theory (see \cite{PiS, Si}, with additions from \cite{Z} and \cite{DoS}), we prove that, in
%absence of sliding (precise definitions are given in Section 2.2), for $u\geq u_*$
%there are multiple (but finitely many) EGMs. %; consequently, {\color{blue}the collection of HC DLR-measures $\sG$ is a polytop whose vertices form set $\sE$. ???}
%{\color{blue}When the value $D$ is of Class A (see Section \ref{SubSec2.1}), we provide a complete 
%description of EGMs $\bmu\in\sE$. In particular, we specify the cardinality 
%$\sharp\sE$, establish symmetries between the EGMs and clarify connections 
%between the PGSs and the EGMs. For a complementary Class B
%we offer a less explicit result. Cf. Theorems 1 - 3 in Section 2.3 and Theorem 4 in Section
%4.1. Nevertheless, for both Classes A and B we establish the coexistence of EGMs (pure phases), which implies the existence of a phase transition \cite{Ge}. DO WEE NEED THIS ?}

The current paper consists of 8 Sections. In Section 2 we state our main results: see Theorem S, Theorems 1 - 3 and Lemma I.
Section 3 derives the proof of these assertions from existing results and methods. In particular,
Section 3.3 contains the definition of a contour and the Peierls bound (see Lemma II). In Section 4
we provide some additional statements; cf. Theorem 4. We also discuss the issue of dominance 
in Section 4.2.
Sections 5, 6 establish connections with algebraic number theory, namely, with
{\it norm equations} in the ring $\bbZ\left[{\sqrt[6]{-1}}\right]$.
More precisely, uniqueness and non-uniqueness in problem (2.1) is related to cosets in
$\bbZ\left[{\sqrt[6]{-1}}\right]$ by the group of units. In particular, in Section 5.2  we prove an
important property of (eventual) minimality along a coset. Next, in Section 6 we discuss
simplest examples of infinite sequences of values $D$ of Classes A and B, respectively.
The proof of Theorem 3 (in a more detailed form - Theorems 5-10) is carried out in Section 7. Section 8 contains the proof of Theorem S claiming that all instances of sliding are identified in Table 2; see Section 2.2. 

% We also provide a supplement to the current paper which includes computer programs
% referenced in the main body of the text.

%Section 8 focus on examples of infinite sequences of values $D$ from Class B.
%Cf. Theorems 4--9.

\section{Periodic ground states in $\bbZ^2$. Sliding. Main results}\label{Sec2}

\subsection{Periodic ground sates and an optimization problem for sub-lattices.}\label{SubSec2.1} The Pirogov-Sinai (PS) theory is based on two
(mutually related) key properties: (I) a finite number of periodic ground states (PGSs) and (II)
a {\it Peierls bound}: a lower estimate for an energy increment in a local perturbation
of a PGS. An admissible configuration $\vphi$ is called a HC
{\it ground state} if one cannot remove
finitely many occupied sites from $\vphi$ and replace them by a larger number of
particles without breaking $D$-admissibility. In other words, one cannot find a finite subset
$\bbV\subset\bbZ^2$ and a configuration $\psi^\bbV\in\cA (\bbV\|\vphi )$ compatible with $\vphi$
such that $\sharp\psi^\bbV>\sharp\vphi\upharp_\bbV$. Next, we say that a
configuration $\phi$ is periodic if there exist two linearly independent
vectors $e_1$, $e_2$ such that $\phi (x)=\phi (x+e_1)=\phi (x+e_2)$. %{\color{blue} If a periodic
%configuration $\phi$ contains the origin, we say that $\phi$ is a {\it sub-lattice}. IS PERIODICITY IS ENOUGH? DO NOT WE NEED A PERIODIC LATTICE CONFIGURATION ???} 
The
parallelogram with vertices $0$, $e_1$, $e_2$, $e_1+e_2$ is a {\it fundamental
parallelogram} (FP) for $\phi$. We always assume that $e_{1,2}$ are chosen so that the
shorter diagonal of the FP divides it into two triangles without obtuse angles;
one of these triangles, with a vertex at the origin, is referred to as a {\it fundamental triangle}
(FT). We say that a sub-lattice is isosceles or non-isosceles if the FT is isosceles or not.

Let $\cP=\cP(D)$ denote the set of PGSs for a given attainable $D$; the cardinality of $\cP$ is denoted by $\sharp \cP$. Our first result
(Lemma I in Section 2.3) is that, in absence of sliding, the set $\cP$ is finite and all PGSs are produced from
maximally-dense admissible sub-lattices by means of
$\bbZ^2$-shifts, the rotation
by $\pi/2$ and the reflection about the x-axis in $\bbR^2$. We will use the term $\bbZ^2$-{\it congruence} for each of these transformations. %or, simply congruence for short. 
Since each maximally-dense sub-lattice is generated by a FP/FT of a minimal area we refer to it as an (admissible) min-area sub-lattice. 
%sequences of attainable values $D^2_n$ along which both $S(D_n)$ and $D^2_n-S(D_n)$
%increase monotonically and indefinitely with $n$ (in fact, with exponentially growing minorants).
%Also, along these sequences, the quantity $S(D_n)/(D^2_n{\sqrt 3}/2)-1$ is positive
%and monotonically decreases to $0$.
%Since a min-area admissible FP generates a maximally-dense
%admissible sub-lattice, we often use the term an (admissible) min-area sub-lattice.

% The Peierls condition is a technically more involved property; however,
% it also holds for a non-sliding $D$. It is discussed in detail in Lemma II and
% Eqns (3.4) and (3.5) from Section 3.

\medskip
A natural way to identify a min-area FP for a given $D$ is to solve
the following discrete optimization problem:
$$\beac\hbox{\bf miminize the area of a ${\bbZ}^2$-triangle}\\
\hbox{\bf with one vertex at the origin,}\\
\hbox{\bf with side-lengths $\ell_i\geq D$ and angles $\alpha_i\leq\pi/2$.}
\ena\eqno (2.1)$$
The term $\bbZ^2$-triangle means a triangle with vertices in $\bbZ^2$. The approach with minimizing the area of $\bbZ^2$-triangles is natural because the (formal) amount of triangles in a Delaunay triangulation corresponding to an admissible configuration is always twice the (formal) number of particles in it. An additional constraint in (2.1) is that the optimization is performed over non-obtuse triangles only. The justification of this constraint is the subject of Section 3.1.

In what follows, $S=S(D)$ refers to the doubled minimal area in (2.1); it is
instructive to note that $S(D)\in\bbN$. (Here and below, $\bbN$ stands 
for the set of positive integers.) Also, for all values $\ell_i$ under consideration we have 
$\ell_i^2=m^2+n^2\in\bbN$, with integer $m,n$. Problem (2.1) always 
has a solution but a minimizing $\bbZ^2$-triangle may be non-unique
(in one sense or another). Any minimizing triangle is referred to as an {\it M-triangle}. A notable part of this work attempts at clarifying
algebraic aspects of non-uniqueness in (2.1). From now on we 
suppose that the triple of side-lengths \ $\ell_i$, $i=0,1,2$, obeys
$$D\leq\ell_0\leq\ell_1\leq \ell_2,\eqno (2.2)$$
and $\alpha_i$ stands for the opposite angles: $\pi/2\geq\alpha_0\geq
\alpha_1\geq\alpha_2$. %It will be sometimes convenient to refer to a 
%$\bbZ^2$-triangle $\triangle$ as a triple $(\ell_0,\ell_1,\ell_2)$ (or an integer triple 
%$(\ell^2_0,\ell^2_1,\ell^2_2)$) satisfying (2.2). 

Of course, one cannot inscribe an equilateral triangle in $\bbZ^2$. However, a pair of integers $a= \ell^2_1 - \ell^2_0$ and $b = \ell^2_2 - \ell^2_1$, called  
a {\it signature} of a $\bbZ^2$-triangle, characterizes the closeness of the triangle to an equilateral one.

We say 
that a given value $D$ (or $D^2$) is of Class S if there exist two M-triangles with a common side but different third vertices belonging to the same half-plane with respect to the common side. Assuming that a dense-packing configuration consists of M-triangles only, it is not hard to 
see that for all $D$ belonging to Class S the HC model exhibits sliding. Theorem S proves that all sliding values of $D$ belong to Class S, 
and there are only 39 of them.

We say that the non-sliding value $D$ (or $D^2$) is of Class A if the M-triangle 
in (2.1), (2.2)
is unique modulo $\bbZ^2$-congruences. That is, all M-triangles have
the same triple $(\ell_0,\ell_1,\ell_2)$ and are $\bbZ^2$-congruent. We show that Class A contains infinitely
many values of $D$ (or $D^2$). An initial list of Class A values of $D^2$ is given in Table 3, together with the corresponding 
values of $S(D)$ and $\sharp\,\cE (D)$.

For example, for $D^2=328$, the doubled area of an M-triangle is $S=292$,
with $\ell^2_0=328$, $\ell^2_1=338$, $\ell^2_2=346$. An example of such an M-triangle
has vertices at sites $(0,0)$, $(13,13)$ and $(-5,18)$. 
%; it generates a min-area sub-lattice
%which we denote by $\bbE_{(13,13),(-5,18)}$. In addition, the minimizing triangles generate 
%three min-area sub-lattices that are $\bbZ^2$-congruent to $\bbE_{(13,13),(-5,18)}$: in the
%above line of notation they are $\bbE_{(-13,13),(-18,-5)}$, $\bbE_{(13,-13),(-5,-18)}$, $\bbE_{(13,13),(18,-5)}$. 
As a result (see Theorem 1), the number of PGSs and EGMs equals $1168=4S$. (Here the
factor $4$ indicates four min-area sub-lattices for $D^2=328$.)

All values of $D$ not belonging to Classes A and S necessarily exhibit non-uniqueness of a solution to (2.1), (2.2).
One form of such non-uniqueness is where the corresponding triple $(\ell_0,\ell_1,\ell_2)$ is unique but the M-triangles form more than one $\bbZ^2$-congruence class. % (i.e., more than one equivalence class modulo $\bbZ^2$-congruences). 
In other 
words, the minimizing triangle for a given $D$ (or $D^2$)  is unique up to an $\bbR^2$-congruence 
but has more than one $\bbZ^2$-{\it implementation}. In this case we say that a value 
$D$ (or $D^2$) belongs to Class B0 and  speak of different 
$\bbZ^2$-implementations for a given minimizing
triple (or minimizing type) $(\ell_0,\ell_1,\ell_2)$.
%Class B0: the minimizing triangle for a given $D$ (or $D^2$)  is unique up to 
%$\bbR^2$-congruence but has more than one $\bbZ^2$-implementation.
An initial example of Class B0 is $D^2=425$; see Table 4. Here the M-triangle has the doubled area $S=750$. The squared side-lengths in an 
M-triangle are $\ell_0^2=\ell_1^2=425$, $\ell_2^2=450$ (so it is isosceles). 
However, such a triangle can be $\bbZ^2$-implemented in two ways: (i) with vertices at $(0,0)$, $(15,15)$
and  $(-5,20)$, and (ii) with vertices at  $(0,0)$, $(21,3)$ and $(8,19)$. These $\bbZ^2$-implementations 
are $\bbR^2$- but not $\bbZ^2$-congruent. Therefore, they will generate not $\bbZ^2$-congruent 
min-area sub-lattices. %Consequently, only one of them generates a dominant family of PGSs. Thus, {\color{blue}we know for sure DO WE ???} that the number of s for $D^2=425$ is 
%$1500=2S$ (although the question which implementation, (i) or (ii) above, generates a 
%dominant family requires an additional work). Again, cf. Table 4. 

In fact, each equivalence class includes $1$, $2$ or $4$ $\bbZ^2$-congruent min-area sub-lattices, 
depending on whether the M-triangle is isosceles straight,
isosceles non-straight or non-isosceles. (An isosceles straight triangle occurs only for
$D^2\leq 20$.) This explains why the number of PGSs for $D^2=328$ is $4S$ while 
for $D^2=425$ it is $2S$. 

%In our computer-assisted enumeration,
%among $D\leq 2000$ we found examples of attainable values of $D$ from Class B0
%with $4$ non-equivalent implementations, although theoretically the number of
%non-equivalent implementation is unbounded. 
We establish that Class B0 includes infinitely many values of $D$. Moreover, 
the degree of degeneracy (the number of different $\bbZ^2$-implementations) can be 
arbitrarily large. We present an initial list of Class
B0 values $D$ in Table 4, together with the specified values of $S(D)$ and $\sharp\,\cP (D)$.

Another form of non-uniqueness is where the triple $(\ell_0,\ell_1,\ell_2)$ for an M-triangle is non-unique. Here
we say that a value $D$ (or $D^2$) belongs to Class B1 if (I) it generates at least $2$
different minimizing triples $(\ell_1,\ell_2,\ell_3)$ in problem (2.1), (2.2)
and (II) each triple has a single implementation modulo $\bbZ^2$-congruences. 

An example of Class B1 is $D^2=65$; again see Table 4. Here we have two 
minimizing triples $(\ell_0,\ell_1,\ell_2)$: (a) with $\ell_0^2=\ell_1^2=65$, $\ell_2^2=80$ and 
(b) with $\ell_0^2=\ell_1^2=68$, $\ell_2^2=72$; both isosceles and of the doubled area $S=120$.
The M-triangles of type (a) are obviously not $\bbR^2$-congruent to those
of type (b). But each of triples (a), (b) has a unique $\bbZ^2$-implementation (again, 
up to $\bbZ^2$-congruences). The vertices can be selected as $(0,0)$, $(8,4)$ and $(-1,8)$
for (a)  and $(0,0)$, $(8,2)$ and $(2,8)$ for (b). %As a result, we obtain {\color{blue} two non-$\bbR^2$-congruent 
%pairs} of min-area-sub-lattices that are $\bbZ^2$-congruent within each pair.% Consequently, 
%there will be two families of PGSs but only one of them will generate EGMs. This example is 
%relatively easy to calculate: here we assert that it is the $(68,68,72)$-type that generates 
%%%dominant PGSs. Consequently, the number of EGMs for $D^2=65$ equals $240=2S$; cf. Table 4. 
%{\color{blue}At present, this statement has not been fully proven but we are certain that our answer is correct. ???}
%Cf. Section 4.2. As above, we prove that Class B1 have infinite occurrences 
%(and admits unlimited degrees of degeneracy of the minimizing triple $(\ell_0,\ell_1,\ell_2)$). 
%The computer-assisted enumeration
%detects attainable values $D\leq 2000$ from Class B1 with $2$, $3$, $4$ or $5$ different triples.
%Again, we are able to check in Sections 6--8 that Class B1 also includes infinitely many
%values of $D$.

We give an initial list of Class
B1 values $D$ in Table 4, together with the specified values of $S(D)$ and $\sharp\,\cP (D)$.

In addition, if, for a given $D$, condition (I) holds true  while condition (II) is not true then 
this $D$ (or $D^2$) belongs to Class B2; as before it also occurs for infinitely many values 
of $D$. That is, for a given value $D^2$ of Class B2 we have (i) at least two 
non-$\bbR^2$-congruent minimizing triangles, and (ii) at least one of them has 
at least two non-$\bbZ^2$-congruent implementations. The first example of Class
B2 is $D^2=15005$. Here we have three minimizing triples $(\ell_0,\ell_1,\ell_2)$: (a) with $\ell_0^2=15005,  \ell_1^2=15076$, $\ell_2^2=15133$, (b) with $\ell_0^2= 15028, \ell_1^2=15041$, $\ell_2^2=15145$ and (c) with $\ell_0^2=15013,  \ell_1^2=15061$, $\ell_2^2=15140$; all of them are non-isosceles and of the doubled area $S=13052$. The M-triangles of types (a), (b) and (c) are obviously pair-wise not $\bbR^2$-congruent. Moreover, each of triples (a) and (c) has a unique $\bbZ^2$-implementation (again, 
up to $\bbZ^2$-congruences). The vertices for triple (a) 
can be selected as $(0,0)$, $(91,82)$ and $(-26,120)$. The vertices for triple (c) 
can be selected as $(0,0)$, $(118,33)$ and $(30,119)$. On the contrary, triple (b) has two implementations: the vertices can be selected as $(0,0)$, $(108,58)$ and $(4,123)$ for one of them and $(0,0)$, $(122,12)$, and $(51,112)$ for the other one.

%All these facts  can be related to
%some number-theoretic properties of value~$D$.

In what  follows, we collectively refer to the set of the values $D$ with non-uniqueness and
non-sliding as Class B. By definition, B0, B1 and B2 form a partition of Class~B. Moreover, all 
attainable values of $D^2$ are partitioned into Class S, Class A and Class B.

%{\bf Remark 2.1.}
\br {\rm For the HC model in $\bbA_2$ only the non-uniqueness type B0 is possible; cf.
\cite{MSS1}. Also note that, contrary to \cite{MSS1}, for lattice $\bbZ^2$ we do not partition
 Class~A into Classes A0 and A1 as an acute triangle cannot be non-inclined with respect to $\bbZ^2$.}
 \quad $\blacktriangle$ \er

Observe that an M-triangle corresponding to a given $D$ can also be minimal for smaller attainable exclusion distances. Moreover, for each $\bbZ^2$-triangle which is minimal for some $D$ there exists a minimal 
attainable exclusion distance, denoted by $D^*$, for which this triangle is still minimal. Consequently, this triangle is an M-triangle for a range of values of $D$ starting with $D^*$ and ending with the shortest side-length in the triangle. This observation allows us to speak about M-triangles without specifying the corresponding $D$ (or $D^2$) unless it is required by the context. Cf. Section 5.

%-----------------------%
%{\color{blue}We would like to note that, for all attainable values $D^2>20$,
%the minimal area $S(D)$ of the FP for a $D$-admissible sub-lattice is $<D^2$
%(and the FP itself is not a square). !!!! }%In Section 6 we  produce
%
%------------------------------

Every min-area sub-lattice $\vphi$ gives rise to a collection of PGSs
obtained by $\bbZ^2$-shifts of $\vphi$; the number of such PGSs is equal to $S(D)$.
Consequently, a given congruence class of sub-lattices produces $m\cdot S(D)$ PGSs where
$m=1,2,4$ as stated in paragraph {\bf{(ii)} }in Section 1.1. Thus, we can speak of {\rm{PGS}} equivalence classes where
every {\rm{PGS}} is obtained from another {\rm{PGS}} by means of a $\bbZ^2$-shift,
rotation by $k\pi/2$ and reflections about co-ordinate axes and bisectrices.
This defines ${\bbZ}^2$-{\it congruences} (or ${\bbZ}^2$-{\it symmetries}) within a  
{\rm{PGS}} equivalence class.

In the rest of this section we discuss the relationship between PGSs and EGMs. We say that a PGS $\vphi\in \cP$ generates an EGM $\bmu_{\vphi} \in \cE$ if 
$$\bmu_\vphi =\lim\limits_{\bbV\nearrow\bbZ^2}\mu_\bbV(\;\cdot\;||\vphi )\eqno (2.3)$$
(see (1.1)--(1.3)).

Note that the {\rm{PGS}} equivalence class is unique iff $D$ is of Class A. Our Theorem 1 (see Section 2.3) states that in this case the structure of the set $\cE$ of the
extreme Gibbs measures for large $u$ is relatively simple: each EGM
$\bmu\in\cE$ is generated from a {\rm{PGS}} $\vphi\in\cP$ and each {\rm{PGS}} $\vphi\in\cP$ generates
an EGM $\bmu=\bmu_\vphi$.

On the opposite end, for values $D$ of Class B not all {\rm{PGS}} classes will generate EGMs,
only dominant ones possessing broader varieties of local excitations of
a larger statistical weight. (For definitions and details see \cite{Z} where the term stable is used in place of dominant.) Our Theorem 2 (see Section 2.3) contains the corresponding results. Due to the underlying symmetries, if a PGS $\vphi$ from a dominant PGS equivalence class generates an EGM $\bmu_{\vphi}$ then every PGS from this class generates a corresponding EGM. The identification of dominant PGSs is a non-trivial task, see the discussion in Section 4.2. 

%Physically speaking, a particular perturbation theory
%emerges here (for measures $\mu_\bbV (\;\cdot\;||\vphi )$, where a perturbation of a given
%order is related to the contribution from `locally excited' configurations of a given
%statistical weight). In simple examples of values $D$ from Class B (with $D^2=65,130,324$),
%the issue of dominance is  expected to be resolved at the level of local excitations of
%statistical weight $u^{-2}$ (relative
%to the PGS) where we vacate 3, 4, 5 or 6 (specifically selected) occupied sites
%from a {\rm{PGS}} $\vphi$ and add $1$, $2$, $3$ or $4$ particles placed at (again, specifically
%chosen) positions that were previously empty in $\vphi$.

To summarize, we may say that $\sharp \cP(D)$ is a sum, over all PGS equivalence classes, of the numbers of PGSs in each class (which may be different for different classes). Accordingly, for $u$ large enough $\sharp \cE(D)$ is the sum over all dominant PGS equivalence classes of the number of PGSs in each class. Note that for such $u$ non-periodic EGMs do not occur (see \cite{DoS}).

%Informally, the description of the PGSs and EGMs  can be related to (2.1) via formulas
%(2.4.1,2) below. For a general non-sliding $D$ we can write
%$$\sharp\cP(D)=\sum\limits_{\substack{{\bf all\; minimizing}\\
%{\bf equivalence\; classes}}}\sharp \;{\rm PGSs}\eqno (2.4.1)$$
%and -- the for $u$ large enough --
%$$\sharp {\mathcal E}(D)=\sum\limits_{\substack{{\bf all\; dominant}\\
%{\bf equivalence\;classes}}}\sharp \;{\rm PGSs}.\eqno (2.4.2)$$
%We want to stress that dominance is a class property as it is preserved under
%$\bbZ^2$-symmetries.

%Speaking of  minimizing triangle, we use the acronym M-triangles.
%Recall, they (a) have angles $\alpha_i\leq\pi/2$,
%(b) have the half-integer area $S(D)/2$ and (c) give rise to max-density/min-area sub-lattices.

\subsection{Sliding as a pattern of non-uniqueness in Eqn (2.1).}\label{SubSec2.2} We say 
that a given value $D$
exhibits \ {\it sliding} \ if exist  two or more M-triangles %${\tt T}^{(1)}$, ${\tt T}^{(2)}$, $\ldots$, 
with (i) a common side called a {\it sliding base}, with two common vertices,
and (ii) distinct third vertices lying in the same half-plane relative to the shared side. In what 
follows we refer to the sliding base as $OW$ and the third vertices as $A,B,\ldots $. Cf. Figures 
\ref{Fig1}C and \ref{Fig6}. % The whole collection of values of $D$ (or $D^2$) with sliding is referred to as {\it Class S}.

Thus, sliding can be viewed as a specific form of non-uniqueness in problem (2.1). It leads
to a multitude of periodic and non-periodic ground states characterized by layered or
staggered patterns. The point is that under sliding there are countably many periodic and
continuum of non-periodic ground states. In fact, here both assumptions I and II of the PS theory
are violated (a finite collection of PGSs and the Peierls bound).

The cases of sliding are listed in Table 2 below. The data format in Table 2 is
$$D^2=m^2+n^2\;\;\;\,W\;\;\;\,
[A,B,\ldots ]\;\;\;\,S(D)$$%\eqno (2.5)$$
where $W=(w_1,w_2)$ and $A=(a_1,a_2)$, $B=(b_1,b_2)$, $\ldots$.
We would like to note that in Table 2 we omitted values of $D^2$
for which the M-triangles and areas $S(D)$ are re-produced for a larger value of
$D^2$ viz., $D^2=1513$ and $D^2=1514$ have the same $W=(0,-39)$ and
the same $A=(34,-19)$, $B=(34,-20)$ and $S(D)=1326$ as for $D^2=1517$.
In other words, we only list the values $D^2$ that are maximal for given
M-triangles.

Examples of bases $OW$ and vertices $A,B,C,\ldots$ are shown in Figure \ref{Fig6}.

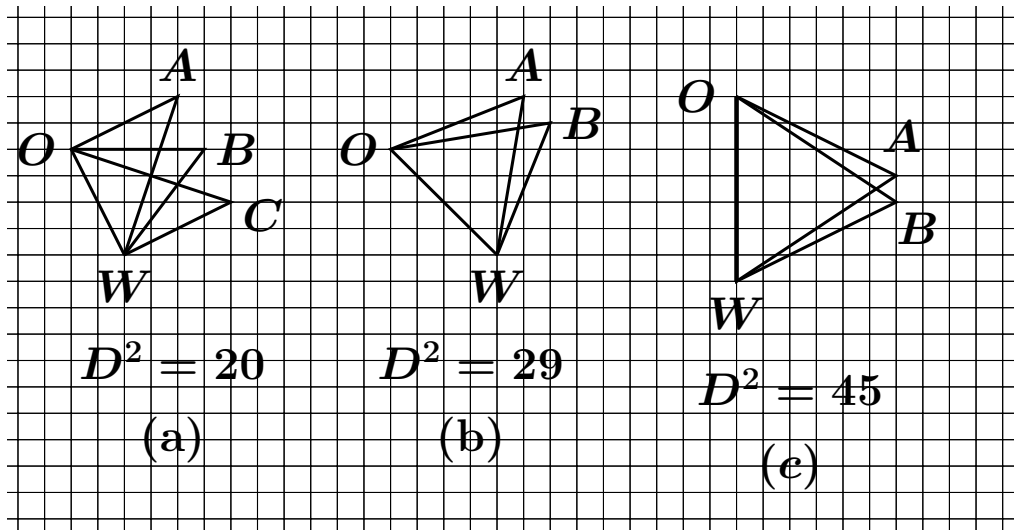
\begin{figure} 
\centering

\begin{tikzpicture}[scale=0.35]\label{SlidingD^2=20}
\clip (-5.4, -6.4) rectangle (32.4, 13.4);

\path [draw=black, line width=0.4mm] (-1,4)--(-3,8)--(1,10)
--(-1,4);
\path [draw=black, line width=0.4mm] (-3,8)--(2,8)
--(-1,4);
\path [draw=black, line width=0.4mm] (-3,8)--(3,6)
--(-1,4);

\draw[yscale=1, xslant=0] (-90,-60) grid (90, 60);
\draw[yscale=1, yslant=0] (-90,-60) grid (90, 60);

\draw (-4.3,8) node  {\Large {${\mbox{\boldmath$O$}}$}};
\draw (1,11.2) node  {\Large {${\mbox{\boldmath$A$}}$}};
\draw (3.2,8) node  {\Large {${\mbox{\boldmath$B$}}$}};
\draw (4.2,5.5) node  {\Large {${\mbox{\boldmath$C$}}$}};
\draw (-1,2.8) node  {\Large {${\mbox{\boldmath$W$}}$}};
\draw (0.8,0) node  {\Large {${\mbox{\boldmath$D^2=20$}}$}};
\draw (0.8,-3) node  {\Large {${\mbox{\boldmath$\rm (a)$}}$}};

\path [draw=black, line width=0.4mm] (13,4)--(9,8)--(14,10)
--(13,4);
\path [draw=black, line width=0.4mm] (9,8)--(15,9)--(13,4);

\draw (7.8,8) node  {\Large {${\mbox{\boldmath$O$}}$}};
\draw (14,11.2) node  {\Large {${\mbox{\boldmath$A$}}$}};
\draw (16.2,9) node  {\Large {${\mbox{\boldmath$B$}}$}};
\draw (13,2.8) node  {\Large {${\mbox{\boldmath$W$}}$}};
\draw (12,0) node  {\Large {${\mbox{\boldmath$D^2=29$}}$}};
\draw (12,-3) node  {\Large {${\mbox{\boldmath$\rm (b)$}}$}};

\path [draw=black, line width=0.4mm] (22,10)--(22,3)--(28,7)
--(22,10);
\path [draw=black, line width=0.4mm] (22,10)--(28,6)--(22,3);
\path [draw=black, line width=0.6mm] (22,10)--(22,3);

\draw (20.5,10) node  {\Large {${\mbox{\boldmath$O$}}$}};
\draw (28.2,8.5) node  {\Large {${\mbox{\boldmath$A$}}$}};
\draw (28.8,5) node  {\Large {${\mbox{\boldmath$B$}}$}};
\draw (22,1.8) node  {\Large {${\mbox{\boldmath$W$}}$}};
\draw (24,-1) node  {\Large {${\mbox{\boldmath$D^2=45$}}$}};
\draw (24,-4) node  {\Large {${\mbox{\boldmath$(c)$}}$}};

%\filldraw[white] (4, 12) circle (.4);
\end{tikzpicture}
%\end{center}
\caption
{\footnotesize{Sliding for $D^2=20$, $D^2=29$ and $D^2=45$. Note the
rectangular $OACW$ (a) and circular trapezes $OABW$ (b,c).}}
\label{Fig6}
\end{figure}

The list of sliding values $D$ in Table~2 is produced by {\tt MinimalTriangles.java}
and {\tt ZSliding.java}. 
% This program analyzes min-area 
% FTs for all values $D\leq 2000$ and reveals only $39$ cases of sliding for $D\leq 2000$, which 
% are shown in Table 2. The largest found value of $D$ with sliding is $195$.

Paper \cite{K} confirms our conjecture that the sliding occurs only for a finite number 
of values $D^2$ i.e., the Class S is finite. That is, there exists a constant $d^\star\in (0,\infty )$ such that
sliding does not occur for $D^2\geq d^\star$. An improvement of the result of \cite{K} is the following theorem. 

\bthms\label{TheoremS}
Table $2$ contains all instances of sliding, i.e., exhausts Class \rS.
\ethms

\begin{table}[H]%[ht]
\centering
%\begin{displaymath}
\boxed{
\begin{array}{llll}
D^2=m^2+n^2&W&[A,B,...]&S(D)\\
%\hline
4=2^2+0^2&(0,-2)&[(2,0),(2,-1),(2,-2)]&4\\
8=2^2+2^2&(2,-2)&[(2,2),(3,1),(4,0)]&8\\
9=3^2+0^2&(0,-3)&[(3,0),(3,-1),(3,-2),&\\
&&(3,-3)]&9\\
18=3^2+3^2&(3,-3)&[(3,3),(4,2),(5,1),(6,0)]&18\\
20=4^2+2^2&(2,-4)&[(4,2),(5,0),(6,-2)]&20\\
29=5^2+2^2&(4,-4)&[(5,2),(6,1)]&28\\
45=6^2+3^2&(0,-7)&[(6,-3)(6,-4)]&42\\
72=6^2+6^2&(6,-6)&[(6,5),(7,4),(8,3),(9,2),&\\
&&(10,1),(11,0)]&66\\
80=8^2+4^2&(0,-9)&[(4,8),(5,8)]&72\\
90=9^2+3^2&(7,-7)&[(9,3),(10,2)]&84\\
106=9^2+5^2&(0,-11)&[(9,-5),(9,-6)]&99\\
121=11^2+0^2&(0,-11)&[(10,-5),(10,-6)]&110\\
157=11^2+6^2&(0,-13)&[(11,-6),(11,-7)]&143\\
160=12^2+4^2&(9,-9)&[(12,4),(13,3)]&144\\
218=13^2+7^2&(0,-15)&[(13,-7),(13,-8)]&195\\
281=16^2+5^2&(12,-12)&[(16,5),(17,4)]&252\\
392=14^2+14^2&(14,-14)&[(19,6),(20,5)]&350\\
521=20^2+11^2&(0,-23)&[(20,12),(20,11)]&460\\
698=23^2+13^2&(0,-27)&[(23,-13),(23,-14)]&621\\
821=25^2+14^2&(0,-29)&[(25,-14),(25,-15)]&725\\
1042=31^2+9^2&(23,-23)&[(31,9)(32,8)]&920\\
1325=35^2+10^2=34^2+13^2&&&\\
\qquad\;=29^2+22^2&(26,-26)&[(35,10),(36,9)]&1170\\
1348=32^2+18^2&(0,-37)&[(32,-18),(32,-19)]&1184\\
1517=34^2+19^2=29^2+26^2&(0,-39)&[(34,-19),(34,-20)]&1326\\
1565=38^2+11^2=37^2+14^2&(28,-28)&[(38,11),(39,10)]&1372\\
2005=41^2+18^2=39^2+22^2&(0,-45)&[(39,-22),(39,-23)]&1755\\
2792=46^2+22^2&(0,-53)&[(46,-26),(46,-27)]&2438\\
3034=55^2+3^2=53^2+15^2&(39,-39)&[(53,15),(54,14)]&2652\\
3709=53^2+30^2&(0,-61)&[(53,-30),(53,-31)]&3233\\
4453=63^2+22^2=58^2+33^2&(0,-67)&[(58,-33),(58,-34)]&3886\\
4756=66^2+20^2=60^2+34^2&(0,-69)&[(60,-34),(60,-35)]&4140\\
6865=76^2+33^2=72^2+41^2&(0,-83)&[(72,-41),(72,-42)]&5976\\
11449=107^2+0^2&(0,-107)&[(93,-53),(93,-54)]&9951\\
12740=112^2+14^2=98^2+56^2&(0,-113)&[(98,-56),(98,-57)]&11074\\
13225=115^2+0^2=92^2+69^2&(0,-115)&[(100,-57),(100,-58)]&11500\\
15488=88^2+88^2&(88,-88)&[(120,33),(121,32)]&13464\\
22784=128^2+80^2&(0,-151)&[(131,-75),(131,-76)]&19781\\
29890=161^2+63^2=147^2+91^2&(0,-173)&[(150,-86),(150,-87)]&25950\\
37970=179^2+77^2=169^2+97^2&(0,-195)&[(169,-97),(169,-98)]&32955
\end{array}
}
%\end{displaymath}
\caption{\footnotesize{The values of $D^2$ with sliding.}}\label{SlidingDs} 
\end{table}
%\cl{\bf Table 2. The values of $D^2$ with sliding, for $D^2\leq 4040093$}

The proof of Theorem S is presented in Section \ref{Sec8}. A part of the proof is 
computer-assisted. The analytical upper bound upon constant $d^\star$ emerging from \cite{K}
was too high for performing an exhaustive enumeration of all cases of sliding for $D^2< d^\star$. 
We therefore needed (i) a better bound on $d^\star$ and (ii) an efficient algorithm of analysis
of values  $D^2< d^\star$. Cf. Section \ref{Sec8}.

%We {\bf conjecture} that for any value $D$ from Class S there is a unique Gibbs measure 
%for large $u$. 
%However, this seems a difficult question.
%that had enticed a number of colleagues for some time (with different opinions about
%possible answers).

%The list of sliding values $D$ in Table 2 is produced by using a program MinimalTriangles.java;
%cf. Section \ref{program} and the supplement to the paper.

\subsection{Main results}\label{SubSec2.3} 

In this section we state Lemma I and Theorems 1 - 3. The proofs of Theorems 1, 2 
are given in Section 3. The proof of Theorem 3 can be recovered from the material presented in Sections 5-7. In fact, Section 5.1 contains a finer version of this theorem (see Theorem~5-10).

A justification of problem (2.1) is provided in Lemma I (where
we use an observation from \cite{ChW}).

%{\bf Lemma I.}
\blI
{\rm{(i)}} For any attainable $D$, every  {\rm{PGS}} is obtained as a
tessellation by \ {\rm M}-triangles and their $\bbZ^2$-shifts.

{\rm{(ii)}} Furthermore, if \
$D$ \ is non-sliding then every \ {\rm{PGS}} \ is obtained from a  max-dense sub-lattice
by means of $\bbZ^2$-congruences. Consequently,
for any non-sliding $D$ the {\rm{PGS}} set \ $\cP(D)$ is finite.
\elI

Let us repeat: a value $D$ (or $D^2$) is of  Class  A if $D$ produces a unique triple
$(\ell_0,\ell_1,\ell_2)$, unique $\bbZ^2$-congruence class and no sliding. On the other hand,
a non-sliding $D$ (or $D^2$) is of Class B if value $D$ features several
equivalence classes in (2.1).

%{\bf Theorem 1.}
\bthma\label{Theorem1}\;{\rm{(Main result, I)}}
Assume an attainable value $D$ is of \  Class \rA. Then:
\begin{description}
  \item[(i)] The cardinality \ $\sharp\cP (D)$ equals $mS(D)$ where
$m=1,2$ or $4$. More precisely, {\rm{(a)}} $m=1$ for $D^2=1,2$,  {\rm{(b)}} $m=2$
if $D^2>2$ and the \rM-triangle is
isosceles, {\rm{(c)}} $m=4$ if the \ \rM-triangle  \ is non-isosceles.
The {\rm{PGS}}s are obtained from each other by $\bbZ^2$-congruences.
  \item[(ii)] There exists a value $u_*(D)\in (0,\infty )$ such that for $u\geq u_*(D)$ the following
assertions hold true. Every {\rm{EGM}} \ $\bmu\in\cE (D)$ is generated by a \ {\rm{PGS}}
$\varphi\in{\cP}(D)$ \ in the sense of {\rm{(2.3)}}: $\displaystyle\bmu
=\lim\limits_{\bbV\nearrow\bbZ^2}\mu_\bbV
(\,\cdot\,||\varphi )(=\bmu_\vphi)$.
The measures $\bmu_\vphi$ are mutually disjoint ($\bmu_{\vphi'}\perp\bmu_{\vphi''}$ for
$\vphi'\neq\vphi''$) and inherit the symmetry properties of their respective {\rm{PGS}}s.
Consequently, $\sharp{\mathcal E}(D)$ $=$ $\sharp{\cP}(D)$.
\end{description}
\ethma

Table 3 indicates the number of extreme Gibbs measures for values $D$ from Class A. See also Figures 3, 4. Values $D$ in Tables 3, 4 are produced by {\tt MinimalTriangles.java}.

%{\bf Theorem 2.}
\bthmb\label{Theorem2}\;{\rm{(Main result, II)}}
Suppose $D$ is attainable,
non-sliding and of Class {\rB}. Then the number of \ {\rm{PGS}}s \ in a given congruence class
is \ $2S(D)$ \ if the respective \rM-triangle is isosceles and \ $4S(D)$ \ if it is non-isosceles.
The total number of \ {\rm{PGS}}s equals the sum of the cardinalities of the equivalence
classes.

%However, in general not all PGSs generate extreme Gibbs measures.
Furthermore, there exist congruence classes (one or several) called {\bf  dominant} 
such that for $u\geq u_*(D)\in (0,\infty)$: \ {\rm{(i)}} \ every {\rm{PGS}} \ $\vphi$ from a dominant class
generates an {\rm{EGM}} $\bmu =\bmu_\vphi$ with  the same symmetries via the limit {\rm{(2.3)}},
{\rm{(ii)}} all \ {\rm{EGM}}s \ $\bmu\in\cE$ are obtained in such a way from the {\rm{PGS}}s
belonging to dominant classes. The {\rm{EGM}}s \ $\bmu_\vphi$ \ are disjoint for different
{\rm{PGS}}s \ $\vphi$.
\ethmb

\begin{table}
\centering
\begin{displaymath}
%\begin{array}{|rrrrl|}
%\boxed{
\begin{array}{|c|c|c|c||c|c|c|c||c|c|c|c||c|c|c|c|}
\hline
D^2&S&{\rm I}/{\rm N}&\sharp\,\cE&D^2&S&{\rm I}/{\rm N}&\sharp\,\cE&D^2&S&{\rm I}/{\rm N}
&\sharp\,\cE&D^2&S&{\rm I}/{\rm N}&\sharp\,\cE\\
\hline
1&1&{\rm I}&1           &64&56&{\rm I}     &112   &164&146&{\rm N}&584  &233&208&{\rm I}  &416\\
2&2&{\rm I}&2           &68&60&{\rm I}      &120  &170&150&{\rm I}&300   &241&209&{\rm I}  &418\\
5&5&{\rm I}&10         &74&68&{\rm N}    &272  &178&157&{\rm N}&627  &245&217&{\rm I}  &434\\
10&10&{\rm I}&20     &85&75&{\rm I}      &150  &180&162&{\rm N}&648  &256&224&{\rm I}  &448\\
13&12&{\rm I}&24     &89&80&{\rm I}      &160  &181&166&{\rm N}&664  &260&228&{\rm N} &912\\
17&15&{\rm I}&30     &97&86&{\rm N}     &344 &193&168&{\rm I}&336    &265&236&{\rm N} &944\\
25&23&{\rm N}&92    &100&90&{\rm I}    &180 &194&172&{\rm N}&688   &272&240&{\rm I}  &480\\
26&24&{\rm I}&48     &109&101&{\rm N} &404 &197&176&{\rm N}&704  &277&247&{\rm N} &988\\
34&30&{\rm I}&60     &113&102&{\rm N} &408 &200&180&{\rm I}&360    &290&254&{\rm N} &1016\\
37&34&{\rm N}&136 &117&105&{\rm N} &420 &202&181&{\rm N}&724  &293&258&{\rm N} &1032\\
40&38&{\rm N}&152 &128&112&{\rm I}   &224 &205&184&{\rm N}&736  &296&268&{\rm N} &1072\\
41&39&{\rm N}&156 &136&120&{\rm I}  &240 &208&188&{\rm N}&752   &305&269&{\rm N} &1076\\
50&45&{\rm I}&90     &137&124&{\rm N} &496 &212&194&{\rm N}&776 &306&270&{\rm I}  &540\\
52&48&{\rm I}&96     &145&127&{\rm N} &508 &221&195&{\rm I}&390   &320&280&{\rm I}  &560\\
53&52&{\rm N}&204 &148&134&{\rm N} &536 &225&198&{\rm N}&792  &328&292&{\rm N}&1168\\
58&53&{\rm N}&212 &153&135&{\rm I}   &270 &226&203&{\rm N} &812 &333&299&{\rm I} &598\\
\hline
\end{array}
\end{displaymath}
\vskip-10pt
\caption{\footnotesize{
An initial list of quadruples $(D^2,S,{\rm I}/{\rm N},\sharp\,\cE)$ for Class A with $D^2<337$.
Here $S$ stands for the area and I/N for Isosceles/Non-isosceles property of an M-triangle in (2.1).
Furthermore, $\sharp\,\cE$ indicates the number of extreme Gibbs measures. As above, we 
only list the maximal
value of $D^2$ consistent with a given M-triangle. E.g., the entry $(32,30,{\rm I},60)$
preceding $(34,30,{\rm I},60)$ is not listed.}}
\end{table}

 %{\rm I} {\rm N}

%\begin{table}[ht]
%\begin{displaymath}
%\begin{array}{|rrrrl|}
%    \hline
%    \text{factor}   &B    &C_\B    &n_\B    &\text{comment}\\
%    \hline
%    3/1             &1      &2      &2      & \bbZ_3\\
%    \hline
%    9/1             &6      &10     &6      & \bbZ_3^2\\
%    9/2             &7      &8      &2      & \bbZ_9\\
%    \hline
 %   27/1            &23     &34     &47     & \bbZ_3^3\\
%    27/2            &24     &28     &11     & \bbZ_3\times  \bbZ_9\\
%    27/3            &24     &27     &10     &14\text{ elements of order $3$}\\
%    27/4            &24     &28     &13     &20\text{ elements of order $3$}\\
%    27/5            &24     &27     &6      &2\text{ elements of order $3$}\\
%    27/6            &24     &27     &10     &8\text{ elements of order $3$}\\
%    27/7            &25     &26     &2      & \bbZ_{27}\\
%    \hline
%\end{array}
%\end{displaymath}

%\end{table}

%$$\boxed{\begin{array}{c|c|c|c||c|c|c|c||c|c|c|c}
\begin{table}[H]
\centering
\begin{displaymath}
%\begin{array}{|rrrrl|}
%\boxed{
\begin{array}{|c|c|c|c||c|c|c|c||c|c|c|c|}
\hline
D^2&S&\rB 0/1&\sharp\,\cP&D^2&S&\rB 0/1&\sharp\,\cP&D^2&S&\rB 0/1
&\sharp\,\cP\\
\hline
65&60&\rB 1&240             &1600&1400&\rB 0&5600    &3413&2971&\rB 1&23768\\
130&120&\rB 1&480         &1845&1605&\rB 0&12840  &3505&3061&\rB 1&24488\\
324&288&\rB 1&1728       &2098&1837&\rB 1&14696  &3690&3210&\rB 0&25680\\
425&375&\rB 0&1500       &2116&1840&\rB 1&14720  &3701&3225&\rB 1&19350\\
485&430&\rB 0&3440       &2213&1930&\rB 1&15440  &3770&3285&\rB 0&26280\\
562&493&\rB 1&2598       &2245&1960&\rB 0&15680  &3816&3318&\rB 1&26544\\
725&635&\rB 0&5080       &2468&2150&\rB 1&17200  &3865&3359&\rB 1&26872\\
832&728&\rB 1&5824       &2578&2247&\rB 1&17976  &4100&3572&\rB 1&28576\\
986&870&\rB 1&5220       &2609&2277&\rB 1&18216  &4210&3673&\rB 1&29384\\
1010&889&\rB 1&7112     &2650&2315&\rB 0&18520  &4232&3680&\rB 1&29440\\
1124&986&\rB 1&5916     &2725&2370&\rB 0&18960  &4250&3695&\rB 0&29560\\
1234&1075&\rB 1&6450   &2770&2419&\rB 1&19352  &4426&3860&\rB 1&30880\\
1297&1135&\rB 1&9080   &2993&2613&\rB 1&15678  &4441&3875&\rB 1&31000\\
1409&1236&\rB 1&9888   &3041&2654&\rB 1&21232  &4505&3919&\rB 1&31352\\
1489&1307&\rB 1&10456 &3060&2670&\rB 0&21360  &4624&4012&\rB 1&24072\\
1521&1329&\rB 1&10456 &3130&2717&\rB 1&16302  &4709&4102&\rB 1&32816\\
\hline
\end{array}
\end{displaymath}
\vskip-10pt
\caption{\footnotesize{Initial values of Class B with $D^2\leq 4709$. Here we feature a quadruple
$(D^2,S,\rB 0/1,\sharp\,\cP)$. (A higher degree of
degeneracy emerges for larger values of $D$.)
Again we pick the largest value of $D^2$ with a given $S$;
e.g. $S=635$ figures for $D^2=724$ and $D^2=725$; the former value $D^2=724$ is omitted.
In all entries, the number of non-congruent sub-lattices equals $2$.  Accordingly, the number of
the PGSs $\sharp\,\cP$ equals $4S$ if both sub-lattices are isosceles, $6S$ if
one sub-lattice is isosceles and the other not (this can occur in Case B1 but not B0),
and $8S$ if both sub-lattices are non-isosceles.
}}
\end{table}

As we said earlier, for a given $D$, the dominant classes could be determined by an
additional analysis of local \ excitations.  A general scheme, based on a number
of assumptions, was proposed in \cite{Z}. 
%A verification of the assumptions from
%\cite{BrS} for the HC model should be done case-by-case (which requires
%a combination of analytic and computer-assisted arguments). 
In this paper we comment
on the cases $D^2=65,130,324$ (an initial triple of values $D$ from Class B);
see Section 4.2, Figures 10-12. %and also Remark 8.3.
%{\color{red}\ref{Fig9Dom65}--\ref{Fig10dom324}}. 
%The formal proofs are postponed until later parts of this work and will be published elsewhere.
On a triangular lattice $\bbA_2$,
an analysis of dominance has been performed in \cite{MSS1} and on $\bbZ^3$ in \cite{MSSz3}.

%As earlier, the list of values $D$ in Table 4
%is produced by {\tt MinimalTriangles.java}.
%\hline
%\end{array}}$$
%\vskip .5cm
%\cl{\bf Table 4. The values of $D$ of Class B with $D^2\leq 4709$}

The following Theorem guarantees that all Classes A, B0, B1, and B2 are infinite. It also 
provides a coarse description of the algebraic structure of these classes. Here we use the 
cyclotomic ring $\bbZ[\sqrt[6]{-1}]$ and its group of units $\bbU[\sqrt[6]{-1}]$.

\bthmc\label{Theorem3}\;{\rm{(Main result, III)}}
Eeach of Classes {\rA}, {\rB}$0$, {\rB}$1$, {\rB}$2$ contains infinitely many values 
$D^2$. Each non-sliding value 
$D^2$ falls in one of these classes. All $\bbZ^2$-triangles can be partitioned into double infinite
sequences with the alternating signatures $a,b$ and $b,a$ for all triangles in the sequence. 
Moreover, inside each double sequence the number of non-M-triangles is finite. Algebraically, 
these sequences are identified with the cosets in ring $\bbZ[\sqrt[6]{-1}]$ by the unit group 
$\bbU[\sqrt[6]{-1}]$. 

Furthermore, the set of all double sequences, i.e. cosets, can be partitioned into maximal families 
such that for a given family there exists a $\bbZ$-indexation of every double sequence such that 
all M-triangles with the same index $j \in \bbZ$ have the same area. For a given index $j$ let $D^2_j$ 
be the shortest squared side-length among all M-triangles in the family indexed by $j$. Then all 
these $D^2_j$ (of which there are infinitely many) belong to the same Class ({\rA}, {\rB}$0$, 
{\rB}$1$, or {\rB}$2$). 
Every such family is finite. Among these families there are families corresponding (in the above 
sense) to each of Classes  {\rA}, {\rB}$0$, {\rB}$1$, and {\rB}$2$, and the size of the 
corresponding Class {\rB} families can be arbitrarily large.
\ethmc

Theorems~1-3 give a detailed description of the phase diagram of the hard-core model in $\bbZ^2$ for any attainable non-sliding $D^2$ and large enough $u$. %Although %there are some open mathematical problems, the physical picture is rather clear.

However, the phase diagram in the case of sliding remains an open question. Another question is an identification of dominant phases (Class B). Here we expect that in the hard-core model in $\bbZ^2$ there always exists a single dominant congruent class of PGSs. %In each specific case the identification procedure is algorithmically clear (following \cite{Z}) but requires an amount of computations growing exponentially with $D$. 

\section{The Peierls condition. Proof of the main results}\label{Sec3}

We begin this section by developing a technical background needed
for the proof of Lemma I; see Lemmas 3.1--3.6. Next, the same background
is used in the proof of Lemma II establishing the Peierls bound. After that,
Theorems 1 and 2 are deduced from the general PS theory.

Section 3.1 establishes some related definitions and known facts. In our own opinion,
it is quite elementary although seems rather tedious. It should be an easier
reading when one is familiar with the concept of the Delaunay triangulation.
Cf. \cite{DLRS}, Chapters 1--3.

\subsection{Voronoi cells, C-triangles and saturated configurations.}\label{SubSec3.1}
A given $D$-admissible configuration  $\phi\in{\mathcal A}(D)$ has a uniquely
defined collection of {\it Voronoi cells} ${\mathcal V}(x,\phi )$
constructed for the occupied sites $x\in\phi$. If $\phi$ has no unbounded
Voronoi cells then to each cell ${\mathcal V}(x,\phi )$ there is assigned a
finite set of circles centered at the vertices of ${\mathcal V}(x,\phi )$
and passing through $x$. We call them {\it V-circles} in $\phi$. Each
$x\in\phi$ lies in at least one of V-circles but no
$x\in\phi$ falls inside a circle. The sites $y\in\phi$ lying in a given
V-circle form the vertices of a {\it constituting polygon}. These polygons
form a tessellation of $\bbR^2$: they have disjoint interiors, and the union of
their closures gives the entire plane. If a constituting polygon has
$\geq 3$ vertices, it can be divided (non-uniquely) into
{\it constituting} triangles
(in short: C-triangles); this produces the Delaunay triangulation of
$\phi$ (and of $\bbR^2$).

A $D$-AC $\phi$ is called {\it saturated} if no occupied site can be added to
it without breaking
admissibility. A {\it saturation} of a given $D$-AC $\phi$ is a completion
of $\phi$ (in some uniquely defined way) with the maximal possible amount
of added occupied sites.

Clearly, every {\rm{PGS}} configuration is saturated. Saturated configurations are
convenient as a natural initial step in a procedure of
`identifying' PGSs within the set ${\mathcal A}(D)$ of admissible configurations.
The use of saturated configurations also makes more transparent the
derivation of the Peierls bound in Section 3.3.

The idea of a saturated configuration worked well in the study of dense-packed
circle configurations in $\bbR^2$; cf. \cite{ChW}. We attempt to emulate a similar
approach in $\bbZ^2$. It generates some technical complications but we manage to
get through, via Lemmas 3.1--3.5.

%{\bf Lemma 3.1.}
\bl A saturated configuration does not have V-circles of
radius $\geq D+\sqrt{2}/2$.
\el

\bp Suppose there exists a V-circle of radius $\geq D+\sqrt{2}/2$.
The center of the V-circle may not lie in $\bbZ^2$ but is at distance $\leq {\sqrt 2}/2$
from one of the $\bbZ^2$-sites.
Then an additional particle can be added at this site
without breaking admissibility. This contradicts the saturation
assumption.
\ep

We would like to note a difference between Lemma 3.1 and Lemma 2 from \cite{ChW}.
We have a lower bound $D+{\sqrt 2}/2$ whereas in \cite{ChW}, Lemma 2, one has $D$. This
creates a specific technical complication arising in $\bbZ^2$ compared with $\bbR^2$.

%{\bf Lemma 3.2.}
\bl Let $\triangle$ be a C-triangle in an $D$-AC $\phi$
and consider $3$ pair-wise disjoint disks of radius $D/2$
centered at the vertices of $\triangle$. Consider 3 sectors in these disks which
are intersections of the disks with the angles of $\triangle$ and let $\bbS (\triangle )$
denote the union of these sectors. Then the area of $\bbS (\triangle)$, i.e., the sum of
the areas of these 3 sectors, equals $\pi D^2/8$.
\el

\bp To avoid confusion we stress that $\bbS (\triangle )$ not necessarily fits completely inside the corresponding $\triangle$. Nevertheless, the collection of sets $\bbS (\triangle )$ where $\triangle$
runs over C-triangles of $\phi$ forms a partition of the union of the disks
$\operatornamewithlimits{\cup}\limits_{x\in\phi}\bbD(x,D/2)$
(modulo a set of measure $0$).
Here $\bbD(u,r)$ stands for the disk of radius $r>0$ centered at
$u\in\bbR^2$: $\bbD(u,r)=\{y\in\bbR^2:\rho (u,y)\leq r\}$.

For each angle of size $\a$ in $\triangle$ the intersection with
the corresponding disk is a full sector with the angular measure $\a$ and
area $\a D^2/8$. The sum of the triangle angles is $\pi$.
\ep

%{\bf Lemma 3.3.}
\bl The value $S(D)/2$ gives the minimal
area for the lattice triangles in problem {\rm{(2.1)}} with angles $\alpha_i\leq
2\pi/3$ instead of $\alpha_i\leq\pi/2$.
\el

\bp Consider a lattice triangle $\triangle$ with an angle $\a$ and the
opposite side $a$, where $\pi/2< \a \leq 2\pi/3$. The
union of $\triangle$ with its central-symmetric counterpart relative to the
middle of $a$ is a lattice
parallelogram whose longer diagonal is $a$. As $\a\leq 2\pi/3$, then the shorter
diagonal $b$
has length $\geq D$. Then $b$ divides
the parallelogram into two congruent non-obtuse admissible triangles $\triangle'$,
$\triangle''$. The area of the parallelogram is
twice the area of $\triangle'$ as well as twice the area of $\triangle$. By
construction, it is  $\geq S(D)$.
\ep

Note that if the circumradius about a C-triangle $\triangle$ is $\leq D$ (which is
the case for continuous dense-packing) then all angles
of this triangle are $\leq 2\pi /3$ and the maximal side-length is $\leq D\sqrt 3$.
Consequently,  Lemma 3.3 is applicable: the area of $\triangle$ is $\geq S(D)/2$.
However, due to Lemma 3.1, when working with saturated $D$-ACs, we need to
deal with C-triangles where the circumradius is
between $D$ and $D+{\sqrt 2}/2$ which may lead to the maximal side-length
$> D\sqrt 3$.  Such a C-triangle can have area $<S(D)/2$. However, it turns out that
in this case there will be an adjacent C-triangle (sharing a side) with a rather large area, so that the
area of their union is $\geq S(D)+1$. It may also happen that two or three C-triangles, of area
$<S(D)/2$ each, share a common  adjacent  triangle; in this case there will again be
a lower bound upon the area of their union. Such an observation allows us to circumspect
these issues; cf. Lemmas 3.4.1--3.4.5. The culmination is Lemma 3.4.5:
it asserts that triangles with obtuse angles $>2\pi/3$ could be circumvented with the help
of adjacent triangles of a large area.

%{\bf Lemma 3.4.1.}
\bla Suppose that a C-triangle $\triangle$ has the circumradius
$r=D+\delta$ where $0\leq\delta\leq{\sqrt 2}/2$. Then the minimal area of $\triangle$ is
$\diy >\frac{{\sqrt 3}D^2}{4}-\frac{D\delta}{2\sqrt 3}$. Also the longest side in
the min-area triangle has length $\diy <D{\sqrt 3}+\frac{\delta}{\sqrt 3}$.
\ela

\bp Suppose a C-triangle $\triangle$ with vertices $A,B,C$ satisfies the
assumptions of the lemma. Let the side-lengths be $AB=l_0$, $BC=l_1$, $CA=l_2$,
with $D\leq l_0\leq l_1\leq l_2\leq 2r$.
If two side-lengths are $>D$, say $l_1,l_2>D$, then the area of the $\triangle$ can be
made smaller by moving vertex $C$ along the circumcircle towards $B$, until
the length of side $BC$ becomes $D$. Indeed,
in the process of motion $l_0$ remains fixed but the height from $C$ to $AB$ shortens.
(The resulting triangle does not necessarily fit $\bbZ^2$.) Thus, the area of $\triangle$
is lower-bounded by the area of the triangle  with two side-lengths $D$ and the
third side-length $\diy 2D\sqrt{1-\frac{D^2}{4r^2}}$. A direct calculation shows that for $D\geq 1$
and $\delta\in [0,{\sqrt 2}/2)$ the bound $\diy 2D\sqrt{1-\frac{D^2}{4r^2}}
< D{\sqrt 3}+\frac{\delta}{\sqrt 3}$ holds true.
(The right-hand side is simply the Taylor expansion in $\delta$ up to order $1$.)
The area of such a triangle
equals $\diy\frac{D^3}{2r}\sqrt{1-\frac{D^2}{4r^2}}$. As above, we have the bound
$\diy\frac{D^3}{2r}\sqrt{1-\frac{D^2}{4r^2}} >\frac{{\sqrt 3}D^2}{4}-\frac{D\delta}{2\sqrt 3}$.
\ep

%{\bf Lemma 3.4.2.}
\blb Suppose that a C-triangle $\triangle$
with side-lengths $l_0,l_1,l_2$ has the circumradius $r=D+\delta$
where $0\leq\delta\leq{\sqrt 2}/2$. Consider the adjacent C-triangle $\triangle'$ that shares
with $\triangle$ the longest side (of length $l_2$). Then the area of the union
$\triangle\cup\triangle'$ is lower-bounded
by the area of a trapeze inscribed in a circle of radius $r$, with
three sides being of length $D$. Furthermore, the area of $\triangle\cup\triangle'$ is
$\diy\geq\frac{3{\sqrt 3}D^2}{4}-{\sqrt 3}\delta^2$.
\elb

\bp Again, we write
$D\leq l_0\leq l_1\leq l_2\leq 2r$. Two vertices of triangle $\triangle'$ are the
end-points of the side of length $l_2$ and lie in the V-circle of radius $r$ circumscribing
$\triangle$. The third vertex of $\triangle'$ cannot lie inside this V-circle
but can be placed on the circle. It also should lie outside the circles of radius $D$
centered at the end-points of the side of length $l_2$. Under these restrictions, the minimal area
of $\triangle'$ is no less than the area of a triangle inscribed in the V-circle which shares the side
 of length $l_2$ with $\triangle$ and has the other side of length $D$. (Cf. the proof of Lemma 3.4.1.) If we now minimize
the area of $\triangle$, we arrive at the pair $\triangle$, $\triangle'$ forming a trapeze,
as specified in the assertion of Lemma 3.4.2. (Again, the resulting triangle may not
fit $\bbZ^2$.)

The area of the trapeze in question equals
$$\frac{2D^3}{r}\left(\sqrt{1-\frac{D^2}{4r^2}}\right)^3=\frac{3{\sqrt 3}D^2}{4}-{\sqrt 3}\delta^2
+\frac{19\delta^3}{3{\sqrt 3}D}-\frac{113\delta^4}{9{\sqrt 3}D^2}+\ldots\;.$$
A tedious (but straightforward) calculation asserts that for $D\geq 1$ and
$0\leq\delta\leq{\sqrt 2}/2$ this expression is
$\diy\geq\frac{3{\sqrt 3}D^2}{4}-{\sqrt 3}\delta^2$, as claimed in the lemma.
\ep

%{\bf Lemma 3.4.3.}
\blc Suppose that a C-triangle $\triangle$ has the circumradius $r=D+\delta$
where $0\leq\delta\leq{\sqrt 2}/2$. Let $\triangle'$  be the adjacent C-triangle sharing the longest
side with $\triangle$ (cf. Lemma {\rm{3.4.2}}).
\begin{description}
  \item[(i)] Suppose that
$\triangle'$ is adjacent to another C-triangle, $\triangle_1$, with circumradius
$r_1=D+\delta_1$ where $0\leq\delta_1\leq{\sqrt 2}/2$. Then the area of $\triangle'$ is
$\geq D^2{\sqrt{11}}/4$.
\item[(ii)] Further, suppose $\triangle'$ is adjacent to other two C-triangles,
$\triangle_1$ and $\triangle_2$, with
circumradii $r_1=D+\delta_1$ and $r_2=D+\delta_2$ where $ 0\leq\delta_1,\delta_2\leq{\sqrt 2}/2$.
Then the area of $\triangle'$ is $\geq D^23{\sqrt 3}/4$.
\end{description}
\elc

\bp (i) Here the triangle $\triangle'$ has one side-length $\geq D$ and two others
$\geq D\sqrt 3$. (Because in each of its neighboring triangles, $\triangle$ and $\triangle_1$,
the angles opposite to the shared sides are $>2\pi/3$ by construction). The area of such
a triangle is, clearly,  $\geq D^2{\sqrt{11}}/4$.

(ii) In this case all side-lengths of $\triangle'$ are $\geq D\sqrt 3$. Hence, the area is
$\geq D^23{\sqrt 3}/4$.
\ep

%{\bf Lemma 3.4.4.}
\bld
The minimal area $S(D)$ in problem $(2.1)$
satisfies
$$\frac{{\sqrt 3}D^2}{2}<S(D)<\frac{{\sqrt 3}D^2}{2}+{\sqrt 2}D.\eqno (3.1)$$
Consequently, for $D^2\geq 19$, the following holds true.  In the situation of Lemma {\rm{3.4.2}}
we have:
$$\frac{3{\sqrt 3}D^2}{4}-{\sqrt 3}\delta^2\geq S(D)+1.$$
Next, in case {\rm{(i)}} of Lemma {\rm{3.4.3}},
$$\diy\frac{{\sqrt 3}D^2}{4}-\frac{D\delta}{2\sqrt 3}+D^2{\sqrt{11}}/4\geq
3(S(D)+1)/2.$$
Finally, in case {\rm{(ii)}} of Lemma {\rm{3.4.3}},
$$\diy\frac{{\sqrt 3}D^2}{4}-\frac{D\delta}{2\sqrt 3}+3{\sqrt 3}D^2/4\geq
4(S(D)+1)/2.$$
\eld

\bp First, let us prove the two-sided bound (3.1).
Consider a bisector of the line segment connecting two lattice sites at
distance $D$ from each other. On this bisector take a point at distance
$D{\sqrt 3}/2$ from the segment. Consider an angle of measure
$2\pi /3$ originating from this point satisfying the following
conditions: the angle is symmetric with respect to the bisector and the
angle does not contain original two sites at distance $D$ from each
other. Such an angle is uniquely defined and any lattice site inside it
is at distance larger than $D$ from both original sites.

The selected point at distance $\diy D\frac{\sqrt{3}}{2}$ from the segment
belongs to some unit square from the base lattice.  Our
angle of measure $2\pi /3$ contains at least one $\bbZ^2$-point
at distance at most $\diy D\frac{\sqrt{3}}{2} + \sqrt{2}$
from the line connecting the two original sites. Taking this vertex and the two
original sites we obtain the triangle with the double area smaller than
$\diy\frac{{\sqrt 3}D^2}{2}+{\sqrt 2}D$. The lower bound for $S(D)$
is obvious.

Now, in the situation of Lemma 3.4.2, we want to have
$$\frac{D^2}{4} 3\sqrt{3} - \sqrt{3}\; \delta^2 \ge S(D)+ 1.$$
This follows from
$$\frac{D^2 }{4} 3\sqrt{3} - \frac{\sqrt 3}{2} \geq\frac{D^2}{2}\sqrt{3} + D
\sqrt{2}+ 1,\;\hbox{ i.e., }\;
\frac{D^2}{4} \sqrt{3}  \geq D{\sqrt 2}+ 1 + \frac{\sqrt 3}{2},$$
which is true for $D \ge 4.3$, i.e. $D^2\geq 19$.

In the situation of Lemma 3.4.3(i) and Lemma 3.4.3(ii) the argument
follows the same line.
\ep

%{\bf Lemma 3.4.5.}
\ble
For all attainable values of $D$
the following holds true. In the situation of Lemma {\rm{3.4.2}}
we have:
$${\rm{area}}\,(\triangle\cup\triangle')\geq 2(S(D)+1)/2.$$
Next, in case {\rm{(i)}} of Lemma {\rm{3.4.3}},
$${\rm{area}}\,(\triangle\cup\triangle'\cup\triangle_1)\geq
3(S(D)+1)/2.$$
Finally, in case {\rm{(ii)}} of Lemma {\rm{3.4.3}},
$${\rm{area}}\,(\triangle\cup\triangle'\cup\triangle_1\cup\triangle_2)\geq
4(S(D)+1)/2.$$

Also, for each of the aforementioned triangle groups $\{\triangle,\triangle'\}$,
$\{\triangle,\triangle',\triangle_1\}$ and\\ $\{\triangle,\triangle',\triangle_1,\triangle_2\}$,
consider the union of the disks of radius $D/2$ centered at the vertices of the group.
Then the intersection of this union, respectively, with  $\triangle\cup\triangle'$,
$\triangle\cup\triangle'\cup\triangle_1$ and $\triangle\cup\triangle'\cup\triangle_1\cup\triangle_2$
has the area $k\pi D^2/8$ where $k$ is the number of triangles in the group.
\ele

\bp For $D^2\geq 19$ the assertion follows from Lemmas 3.4.2--3.4.4.
For $1\leq D^2\leq 18$ the proof is done by a direct enumeration of obtuse
lattice triangles with the longest side of length $<D{\sqrt 3}+{\sqrt{1/6}}$.
The latter value emerges from the bound $\diy 2D\sqrt{1-\frac{D^2}{4r^2}}$
$<D{\sqrt 3}+1/{\sqrt 6}$; cf. Lemma 3.4.1. The last assertion of Lemma 3.4.5
is straightforward; cf. the proof of Lemma 3.2.
\ep

Lemma 3.4.5 allows us to introduce a {\it re-distributed area} assigned to a C-triangle
which can be conveniently lower-bounded. Namely, for a C-triangle not mentioned
in Lemma 3.4.5, the re-distributed area coincides with its proper area. Next, if a triangle
falls into one of the groups mentioned in Lemma 3.4.5 (as $\triangle$, $\triangle'$,
$\triangle_1$ or $\triangle_2$) then its re-distributed area is set to be the
total area of the group ($\triangle\cup\triangle'$, $\triangle\cup\triangle'\cup\triangle_1$
or $\triangle\cup\triangle'\cup\triangle_1\cup\triangle_2$) divided by the number of
the members in the group. According to Lemma 3.4.5,  the redistributed area of a triangle
will be $>(S(D)+1)/2$.

%For a C-triangle with the angles
%$\leq 2\pi/3$ it coincides with its proper area. More precisely, the minimal C-triangles
%have the re-distributed area equal to $S(D)/2$, their proper area. Next,
%a triangle $\triangle'$ that has no adjacent triangles with an angle $>2\pi/3$ with which
%$\triangle'$ shares its longest side gets a re-distributed area $>(S(D)+1)/2$.
%Further, a triangle $\triangle'$ that has no adjacent triangle

\subsection{Maximum-density configurations. Proof of Lemma I.}\label{SubSec3.2} A $D$-admissible configuration containing
only C-triangles of area $S(D)/2$ (i.e., only M-triangles)
is called {\it perfect}. Clearly, a perfect configuration is saturated. An
example of a perfect configuration is a max-dense sub-lattice.

In general, we define the {\it density} of a $D$-AC $\phi$ as
$$\limsup_{L \to \infty} \frac{\sharp(L,\phi)}{L^2}.\eqno (3.2)$$
Here $\sharp(L,\phi)$ denotes the amount of occupied sites  $x\in\phi$
lying in $\bbQ_L(0,0)$. In turn,  $\bbQ_L(m,n)\subset \bbR^2$, $m,n,L\in\bbZ$,
$L\geq 1$, stands for
the square of an integer odd side-length ${\wt L}=L+1-(L\,{\rm{mod}}\;2)$, centered at
site $(m,n)\in\bbZ^2$: 
$$\bbQ_L(m,n) =\begin{cases} (m-\frac{L}{2},\; m+\frac{L}{2})\times (n-\frac{L}{2},\; 
n+\frac{L}{2}),&\hbox{$L$ odd,}\\
(m-\frac{L+1}{2},\; m+\frac{L+1}{2})\times (n-\frac{L+1}{2},\; 
n+\frac{L+1}{2}),&\hbox{$L$ even.}\end{cases}\eqno (3.3)$$
The intersection of $\bbQ_L(m,n)\cap\bbZ^2$ (with $L^2$  lattice sites) is also
denoted by $\bbQ_L(m,n)$; the specific context determines the meaning of this
notation in the sequel.

%{\bf Lemma 3.5.}
\blf The maximal possible density of a $D$-AC is $1/S(D)$.
This density is attained
on any perfect configuration, in particular, on any max-dense sub-lattice.
\elf

\bp (See \cite{ChW}.) The density of any non-saturated  admissible
configuration is not larger than the density of its
saturation. Therefore, it is enough to analyze the densities of saturated
configurations.

To calculate $\sharp(L,\phi)$ for saturated $\phi\in \cA$ we need to count all
occupied sites in $\phi$ which belong to $\bbQ_L(0,0)$.
%Some of the disks of diameter $D$ centered at these points are located
%partially outside of the $L\times L$ square but their total outside area is not
%larger than $4LD$, which for $L \to \infty$ is a quantity $E_1$ of order $O(L)$.
If $A_\rD(L,\phi)$ denotes the total area of disks of diameter $D$ centered at
these occupied sites then $A_\rD(L,\phi) = {\sharp(L, \phi)\pi{D^2/4}}$.

Consider the total area $A_\rT(L,\phi)$ of all C-triangles such that the triangle
itself or a member of its group determined in Lemma 3.4.5 has a vertex at an
occupied point inside $\bbQ_L(0,0)$. Then $A_\rT(L,\phi)=L^2 + E_1(L, \phi)$,
where $E_1(L, \phi )= O(L)$ (recall that the configuration $\phi$ is saturated
and apply Lemma 3.1). The total area $E_2(L, \phi)$ of the intersection of all
these triangles with the disks of diameter $D$ centered at triangle vertices
located outside of $\bbQ_L(0,0)$ is also a quantity of order $O(L)$.

Thus,
$$\frac{\sharp(L,\phi )\pi{D^2/4}}{L^2} = \frac{A_\rD(L, \phi)}{A_\rT(L, \phi )
-E_1(L, \phi)} = \frac{A_\rD(L, \phi)+ E_2(L, \phi)}{ A_\rT(L, \phi)}
+ O\left(1/L\right).$$
According to Lemma 3.2, 3.3 and 3.4.5, the ratio $\diy\frac{A_\rD(L, \phi)
+ E_2(L, \phi)}{A_\rT(L, \phi)}$ is not larger than
$\diy\frac{\pi D^2/4}{S(D)}$ and equals $\diy\frac{\pi D^2/4}{S(D)}$ for any
perfect configuration.
\ep

Observe that {\rm{PGS}} $\vphi$ is uniquely mapped into a configuration
in a torus $\bbT_{kS(D)}$ of size $kS(D)\times kS(D)$ with an integer
$k=k(\vphi )$. Here an $L\times L$  torus $\bbT_L$ is understood as
square $\bbQ_L(0,0)$ (cf. (3.3))
with the identified opposite sides and with the toroidal Euclidean metric $\rho_L^{\;\rT}$:
$$\rho_L^{\;\rT}(x,y)=\min\Big\{\rho (x,y),\rho (x,y\pm v_i),\;i=1,2 \Big\},
\;\;x,y\in\bbQ_{L}(0,0).\eqno (3.4)$$
Here vectors $v_1=({\wt L},0)$ and $v_2=(0,{\wt L})$ where, as before,
${\wt L}=L+1-(L\,{\rm{mod}}\;2)$.

The image configuration in torus
$\bbT_{kS(D)}$ has a maximal possible amount of occupied sites. Note
that the opposite implication is not
true. A configuration which for some integer $N$ has the maximal possible
amount of occupied sites
inside $\bbT_N$ does not necessarily generates a {\rm{PGS}} unless $N^2$ is
divisible by $S(D)$.

%{\bf Lemma 3.6.}
\blg
For any attainable $D$, all PGSs are perfect configurations.
\elg

\bp For a {\rm{PGS}} $\vphi$ the density is a $\lim$ rather than a $\limsup$. This
value must equal $1/S(D)$ as in the opposite case one can gain more occupied points
by replacing the part of $\vphi$ in a sufficiently large square with the part of  a
max-dense sub-lattice generated by M-triangles in a square of the same size.

Clearly, a {\rm{PGS}} $\vphi$  maps into the corresponding $D$-AC inside torus
$\bbT_{k(\vphi )S(D)}$. If inside this torus there is at least one
non-$D$-optimal triangle then the total number of
occupied sites inside the torus is at most
$\diy\frac{k(\vphi )^2S(D)^2}{S(D)} -1$; consequently, the density is $<1/S(D)$.
\ep

\bplemI %{\bf Proof of Lemma I.}
Lemma 3.6 implies assertion (i) in Lemma I. To deduce assertion (ii),
observe that, in absence of sliding, any two M-triangles sharing a common side are centrally
symmetric relative to the mid-point of the shared side. Thus, for any non-sliding $D$ and
any M-triangle $\triangle$ with a vertex at the origin, the corresponding min-area
sub-lattice is the only possible perfect configuration containing $\triangle$. The
$\bbZ^2$-congruences produce, from any max-dense sub-lattice, 1 or 3 additional
sub-lattices, depending on whether
$\triangle$ is isosceles or not.
\eplemI

We would like to note that here and below the absence of sliding is used in a particular
manner (specified in the proof of Lemma 3.6): it implies that if in an $D$-AC $\phi$ there
are two adjacent M-triangles  sharing a common side then they are centrally symmetric.

\subsection{The Peierls bound.} As was said, an application of the PS theory needs a
Peierls bound. We again begin with some auxiliary notions and statements. From
now on we assume that the value of $D$ is non-sliding.

%{\bf Lemma 3.7.}
\blh
Consider a max-dense sub-lattice $\vphi$. Then for any
$k\in\bbN$
and any lattice site $(m,n)\in\bbZ^2$, the square $\bbQ_{kS(D)}(m,n)$ contains
$k^2 S(D)$ occupied sites from $\vphi$.
\elh

\bp The restriction of $\vphi$ to $\bbQ_{kS(D)}(m,n)$ actually is a
$D$-AC lying inside the corresponding torus $\bbT_{kS(D)}$.
This torus is partitioned into $k^2S(D)$ fundamental parallelograms of
the sub-lattice.
\ep

Given an attainable $D$ and a site $(m,n)\in\bbZ^2$, the $S(D)\times S(D)$ square\\ $\bbQ_{S(D)}(mS(D), nS(D))$ (see (3.3))
is called a {\it template}. Given a $D$-AC $\phi\in{\mathcal A}(D)$ and PGS
$\vphi\in{\cP}(D)$, a template $\bbQ_{S(D)}(mS(D), nS(D))$ is said to be {\it $\vphi$-correct} in
$\phi$ if $\vphi$ and $\phi$ coincide inside $\bbQ_{3S(D)}(mS(D), nS(D))$. A template is called
{\it incorrect} (in $\phi$) if it is not $\vphi$-correct for some $\vphi\in{\cP}(D)$.

A {\it contour} in a $D$-admissible configuration $\phi$ is a connected component $\Gam$ of
incorrect templates in $\phi$, where the connectedness is understood in the $\bbR^2$ sense.
Cf. \cite{PiS}. In particular, a contour has an exterior, Ext $\Gam$, and an interior, Int $\Gam$,
which can be further divided into components Int$_\vphi (\Gam)$.

Contour-related pictures are featured in Figures \ref{Fig7}--\ref{Fig9}.

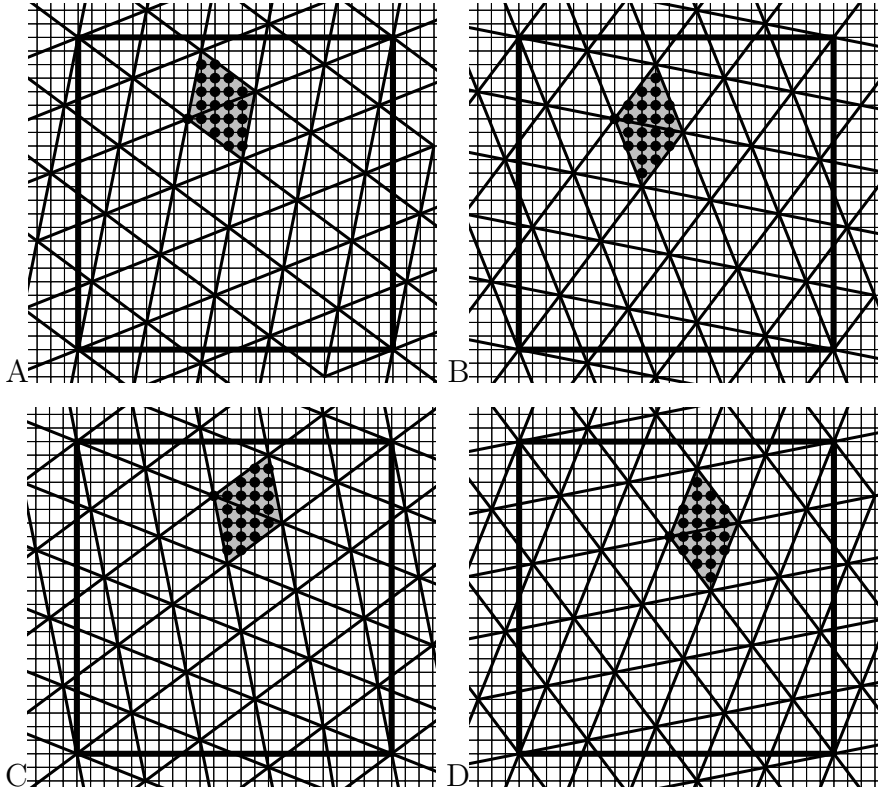
\begin{figure}
\centering
A\begin{tikzpicture}[scale=0.18]\label{1FParaD^2=25}
\clip (-5.6, -15.4) rectangle (24.2, 12.5);

%O: (-2,-13)
\path [draw=black, line width=0.8mm] (-2,-13)--(-2,10)--(21,10)--(21,-13)--(-2,-13);

%\filldraw[lightgray] (10,1)--(15,3)--(19,0)--(18,-5)--(13,-7)--(9,-4)--(10,1);
\filldraw[lightgray] (6,4)--(7,9)--(11,6)--(10,1)--(6,4);

\path [draw=black, line width=0.4mm] (-9,-2)--(26,12); %(5,2)
\path [draw=black, line width=0.4mm] (-5,-5)--(25,7);
\path [draw=black, line width=0.4mm] (-6,-10)--(24,2);
\path [draw=black, line width=0.4mm] (-7,-15)--(28,-1);
\path [draw=black, line width=0.4mm] (2,-16)--(27,-6);  %(5,2)
\path [draw=black, line width=0.4mm] (-8,3)--(22,15);   %(5,2)
\path [draw=black, line width=0.4mm] (-7,8)--(13,16);   %(5,2)
\path [draw=black, line width=0.4mm] (16,-15)--(21,-13)--(26,-11);   %(5,2)

\path [draw=black, line width=0.4mm] (-6,-10)--(-5,-5)--(-1,15); %(1,5)
\path [draw=black, line width=0.4mm] (-3,-18)--(4,17);
\path [draw=black, line width=0.4mm] (2,-16)--(9,19);
\path [draw=black, line width=0.4mm] (6,-19)--(7,-14)--(8,-9)--(9,-4)
--(10,1)--(13,17); %(1,5)
\path [draw=black, line width=0.4mm] (11,-17)--(12,-12)--(13,-7)
--(14,-2)--(15,3)--(17,13)--(18,18); %(1,5)
\path [draw=black, line width=0.4mm] (16,-15)--(17,-10)--(18,-5)--(19,0)--(20,5)
--(21,10)--(22,15); %(1,5)
\path [draw=black, line width=0.4mm] (20,-18)--(21,-13)--(22,-8)--(23,-3)
--(24,2); %(1,5)

\path [draw=black, line width=0.4mm] (-6,-10)--(-2,-13)--(2,-16); %(4,-3)
\path [draw=black, line width=0.4mm] (-9,-2)--(-5,-5)--(-1,-8)--(3,-11)--(7,-14)--(11,-17);
\path [draw=black, line width=0.4mm] (-8,3)--(-4,0)--(0,-3)--(8,-9)--(12,-12)
--(16,-15); %(4,-3)
\path [draw=black, line width=0.4mm] (-7,8)--(-3,5)--(1,2)--(5,-1)--(9,-4)--(13,-7)
--(17,-10)--(21,-13)--(25,-16); %(4,-3)
\path [draw=black, line width=0.4mm] (-6,13)--(-2,10)--(2,7)--(6,4)--(10,1)
--(14,-2)--(18,-5)--(22,-8)--(26,-11); %(4,-3)
\path [draw=black, line width=0.4mm] (-1,15)--(3,12)--(7,9)--(11,6)--(15,3)--(23,-3)--(27,-6);
\path [draw=black, line width=0.4mm] (4,17)--(8,14)--(12,11)--(16,8)--(20,5)--(24,2);
\path [draw=black, line width=0.4mm] (13,16)--(17,13)--(21,10)--(25,7); %(4,-3)

\filldraw[black] (7,8) circle (.35); \filldraw[black] (8,8) circle (.35);  %2
\filldraw[black] (7,7) circle (.35);   \filldraw[black] (8,7) circle (.35);  \filldraw[black] (9,7) circle (.35);  %3
\filldraw[black] (7,6) circle (.35);   \filldraw[black] (8,6) circle (.35);  \filldraw[black] (9,6) circle (.35);
\filldraw[black] (10,6) circle (.35);  %4
\filldraw[black] (7,5) circle (.35);   \filldraw[black] (8,5) circle (.35);  \filldraw[black] (9,5) circle (.35);
\filldraw[black] (10,5) circle (.35);  %4
\filldraw[black] (6,4) circle (.35);   \filldraw[black] (7,4) circle (.35);  \filldraw[black] (8,4) circle (.35);
\filldraw[black] (9,4) circle (.35);   \filldraw[black] (10,4) circle (.35);  %5
\filldraw[black] (8,3) circle (.35);   \filldraw[black] (9,3) circle (.35);  \filldraw[black] (10,3) circle (.35); %3
\filldraw[black] (9,2) circle (.35);   \filldraw[black] (10,2) circle (.35); %2 total: 23

\draw[yscale=1, xslant=0] (-90,-60) grid (90, 60);
\draw[yscale=1, yslant=0] (-90,-60) grid (90, 60);

%\filldraw[white] (-0.5,-11.5) circle (1);
%\draw (-0.5,-11.5) node   {{${\mbox{\boldmath$O$}}$}};
\end{tikzpicture} B\begin{tikzpicture}[scale=0.18]\label{1FParaD^2=8}
\clip (-5.6, -15.4) rectangle (24.2, 12.5);

%O: (-2,-13)
%\filldraw[lightgray] (7,-1)--(12,-2)--(14,-7)--(11,-11)--(6,-10)--(4,-5)--(7,-1);

\filldraw[lightgray] (5,4)--(8,8)--(10,3)--(7,-1)--(5,4);

\path [draw=black, line width=0.8mm] (-2,-13)--(-2,10)--(21,10)--(21,-13)--(-2,-13);

\path [draw=black, line width=0.4mm] (-6,-3)--(-2,-13)--(0,-18); %(2,-5)
\path [draw=black, line width=0.4mm] (-7,11)--(-5,6)--(-3,1)--(-1,-4)--(1,-9)--(3,-14)--(5,-19);
\path [draw=black, line width=0.4mm] (-4,15)--(-2,10)--(0,5)--(2,0)--(4,-5)--(6,-10)--(8,-15)--(10,-20);
\path [draw=black, line width=0.4mm] (1,14)--(3,9)--(5,4)--(7,-1)--(9,-6)--(11,-11)--(13,-16);
\path [draw=black, line width=0.4mm] (6,13)--(8,8)--(10,3)--(12,-2)--(14,-7)--(16,-12)--(18,-17); %(2,-5)
\path [draw=black, line width=0.4mm] (9,17)--(11,12)--(13,7)--(15,2)--(17,-3)--(19,-8)--(21,-13)--(23,-18);
\path [draw=black, line width=0.4mm] (14,16)--(16,11)--(18,6)--(20,1)--(22,-4)--(24,-9)--(26,-14); %(2,-5)

%\path [draw=black, line width=0.4mm] (2,-16)--(27,-6);  %(5,2)
%\path [draw=black, line width=0.4mm] (-8,3)--(22,15);   %(5,2)

\path [draw=black, line width=0.4mm] (-8,2)--(-5,6)--(-2,10)--(1,14);   %(3,4)
\path [draw=black, line width=0.4mm] (-6,-3)--(-3,1)--(0,5)--(3,9)--(6,13); %(3,4)
\path [draw=black, line width=0.4mm] (-7,-12)--(-4,-8)--(-1,-4)--(2,0)--(5,4)--(8,8)--(11,12)--(14,16); %(3,4)
\path [draw=black, line width=0.4mm] (-5,-17)--(-2,-13)--(10,3)--(13,7)--(16,11)--(19,15); %(3,4)
\path [draw=black, line width=0.4mm] (0,-18)--(3,-14)--(6,-10)--(9,-6)--(12,-2)--(15,2)--(18,6)
--(21,10)--(24,14); %(3,4)
\path [draw=black, line width=0.4mm] (5,-19)--(8,-15)--(11,-11)--(14,-7)--(17,-3)--(20,1)--(23,5)--(26,9); %(3,4)
\path [draw=black, line width=0.4mm] (13,-16)--(16,-12)--(19,-8)--(22,-4)--(25,0); %(3,4)
\path [draw=black, line width=0.4mm] (18,-17)--(21,-13)--(24,-9);

\path [draw=black, line width=0.4mm] (-7,-12)--(-2,-13)--(3,-14)--(8,-15)--(13,-16); %(5,-1)
\path [draw=black, line width=0.4mm] (-9,-7)--(-4,-8)--(1,-9)--(6,-10)--(11,-11)--(16,-12)--(21,-13)
--(26,-14); %(5,-1)
\path [draw=black, line width=0.4mm]  (-6,-3)--(-1,-4)--(4,-5)--(9,-6)--(14,-7)--(19,-8)--(24,-9);
\path [draw=black, line width=0.4mm] (-8,2)--(-3,1)--(2,0)--(7,-1)--(12,-2)--(17,-3)--(22,-4)--(27,-5); %(5,-1)
\path [draw=black, line width=0.4mm] (-10,7)--(-5,6)--(0,5)--(5,4)--(10,3)--(15,2)--(20,1)--(25,0); %(5,-1)
\path [draw=black, line width=0.4mm] (-7,11)--(-2,10)--(3,9)--(8,8)--(13,7)--(18,6)--(23,5)--(28,4); %(5,-1)
\path [draw=black, line width=0.4mm] (-9,16)--(-4,15)--(1,14)--(6,13)--(11,12)--(16,11)--(21,10)
--(26,9); %(5,-1)

\filldraw[black] (8,7) circle (.35);  %1
\filldraw[black] (7,6) circle (.35);   \filldraw[black] (8,6) circle (.35);  %2
\filldraw[black] (6,5) circle (.35);   \filldraw[black] (7,5) circle (.35);  \filldraw[black] (8,5) circle (.35);
\filldraw[black] (9,5) circle (.35); %4
\filldraw[black] (5,4) circle (.35); \filldraw[black] (6,4) circle (.35); \filldraw[black] (7,4) circle (.35);
\filldraw[black] (8,4) circle (.35); \filldraw[black] (9,4) circle (.35);   %5
\filldraw[black] (6,3) circle (.35); \filldraw[black] (7,3) circle (.35); \filldraw[black] (8,3) circle (.35);
\filldraw[black] (9,3) circle (.35); %4
\filldraw[black] (6,2) circle (.35); \filldraw[black] (7,2) circle (.35); \filldraw[black] (8,2) circle (.35);
\filldraw[black] (9,2) circle (.35); % 4
\filldraw[black] (7,1) circle (.35); \filldraw[black] (8,1) circle (.35); %2
\filldraw[black] (7,0) circle (.35); %1
%total: 23

\draw[yscale=1, xslant=0] (-90,-60) grid (90, 60);
\draw[yscale=1, yslant=0] (-90,-60) grid (90, 60);

%\filldraw[white] (-0.5,-11.5) circle (1);
%\draw (-0.5,-11.5) node   {{${\mbox{\boldmath$O$}}$}};
\end{tikzpicture}

\vskip 3 truemm

C\begin{tikzpicture}[scale=0.18]\label{1FParaD^2=8}
\clip (-5.6, -15.4) rectangle (24.2, 12.5);

%O: (-2,-13)
\path [draw=black, line width=0.8mm] (-2,-13)--(-2,10)--(21,10)--(21,-13)--(-2,-13);
%\filldraw[lightgray] (9,1)--(14,-1)--(15,-6)--(11,-9)--(6,-7)--(5,-2)--(9,1);
\filldraw[lightgray] (8,6)--(12,9)--(13,4)--(9,1)--(8,6);

\path [draw=black, line width=0.4mm] (-7,-11)--(-2,-13)--(3,-15)--(8,-17); %(5,-2)
\path [draw=black, line width=0.4mm] (-8,-6)--(-3,-8)--(2,-10)--(7,-12)--(12,-14)--(17,-16); %(5,-2)
\path [draw=black, line width=0.4mm] (-9,-1)--(-4,-3)--(1,-5)--(6,-7)--(11,-9)--(16,-11)
--(21,-13)--(26,-15); %(5,-2)
\path [draw=black, line width=0.4mm] (-10,4)--(-5,2)--(0,0)--(5,-2)--(10,-4)--(15,-6)--(20,-8)
--(25,-10); %(5,-2)
\path [draw=black, line width=0.4mm] (-6,7)--(-1,5)--(4,3)--(9,1)--(14,-1)--(19,-3)--(24,-5)
--(29,-7); %(5,-2)
\path [draw=black, line width=0.4mm] (-7,12)--(-2,10)--(3,8)--(8,6)--(13,4)--(18,2)--(23,0)
--(28,-2); %(5,-2)
\path [draw=black, line width=0.4mm] (2,13)--(7,11)--(12,9)--(17,7)--(22,5)--(27,3); %(5,-2)
\path [draw=black, line width=0.4mm] (11,14)--(16,12)--(21,10)--(26,8); %(5,-2)

\path [draw=black, line width=0.4mm] (-6,7)--(-5,2)--(-4,-3)--(-3,-8)--(-2,-13)--(-1,-18); %(1,-5)
\path [draw=black, line width=0.4mm] (-3,15)--(-2,10)--(-1,5)--(0,0)--(1,-5)--(2,-10)--(3,-15); %(1,-5)
\path [draw=black, line width=0.4mm] (2,13)--(3,8)--(4,3)--(5,-2)--(6,-7)--(7,-12)--(8,-17); %(1,-5)
\path [draw=black, line width=0.4mm] (6,16)--(7,11)--(8,6)--(9,1)--(10,-4)--(11,-9)--(12,-14)
--(13,-19); %(1,-5)
\path [draw=black, line width=0.4mm] (11,14)--(12,9)--(13,4)--(14,-1)--(15,-6)--(16,-11)
--(17,-17); %(1,-5)
\path [draw=black, line width=0.4mm] (15,17)--(16,12)--(17,7)--(18,2)--(19,-3)--(20,-8)--(21,-13)
--(22,-18); %(1,-5)
\path [draw=black, line width=0.4mm] (20,15)--(21,10)--(22,5)--(23,0)--(24,-5)--(25,-10); %(1,-5)

\path [draw=black, line width=0.4mm] (-6,7)--(-2,10)--(2,13); %(4,3)
\path [draw=black, line width=0.4mm] (-9,-1)--(-5,2)--(-1,5)--(3,8)--(7,11)--(11,14);
\path [draw=black, line width=0.4mm] (-8,-6)--(-4,-3)--(0,0)--(4,3)--(8,6)--(12,9)--(16,12)--(20,15);
\path [draw=black, line width=0.4mm] (-7,-11)--(-3,-8)--(1,-5)--(5,-2)--(9,1)--(13,4)--(17,7)--(21,10)
--(25,13); %(4,3)
\path [draw=black, line width=0.4mm] (-6,-16)--(-2,-13)--(2,-10)--(6,-7)--(10,-4)--(14,-1)
--(18,2)--(22,5)--(26,8); %(4,3)
\path [draw=black, line width=0.4mm] (-1,-18)--(3,-15)--(7,-12)--(11,-9)--(15,-6)--(19,-3)
--(23,0)--(27,3); %(4,3)
\path [draw=black, line width=0.4mm] (8,-17)--(12,-14)--(16,-11)--(20,-8)--(24,-5)
--(27,-2); %(4,3)
\path [draw=black, line width=0.4mm] (17,-16)--(21,-13)--(25,-10); %(4,3)

\filldraw[black] (11,8) circle (.35);   \filldraw[black] (12,8) circle (.35);  %2
\filldraw[black] (10,7) circle (.35);   \filldraw[black] (11,7) circle (.35);  \filldraw[black] (12,7) circle (.35);  %3
\filldraw[black] (8,6) circle (.35);  \filldraw[black] (9,6) circle (.35);
\filldraw[black] (10,6) circle (.35);  \filldraw[black] (11,6) circle (.35); \filldraw[black] (12,6) circle (.35);  %5
\filldraw[black] (9,5) circle (.35);
\filldraw[black] (10,5) circle (.35);  \filldraw[black] (11,5) circle (.35);  \filldraw[black] (12,5) circle (.35);  %4
\filldraw[black] (9,4) circle (.35);
\filldraw[black] (10,4) circle (.35);  \filldraw[black] (11,4) circle (.35);  \filldraw[black] (12,4) circle (.35);  %4
\filldraw[black] (9,3) circle (.35); \filldraw[black] (10,3) circle (.35); \filldraw[black] (11,3) circle (.35);  %3
\filldraw[black] (9,2) circle (.35); \filldraw[black] (10,2) circle (.35);  %2

\draw[yscale=1, xslant=0] (-90,-60) grid (90, 60);
\draw[yscale=1, yslant=0] (-90,-60) grid (90, 60);

%\filldraw[white] (-0.5,-11.5) circle (1);
%\draw (-0.5,-11.5) node   {{${\mbox{\boldmath$O$}}$}};
\end{tikzpicture} D\begin{tikzpicture}[scale=0.18]\label{1FParaD^2=8}
\clip (-5.6, -15.4) rectangle (24.2, 12.5);

%O: (-2,-13)
%\filldraw[lightgray] (7,-2)--(12,-1)--(15,-5)--(13,-10)--(8,-11)--(5,-7)--(7,-2);
%\filldraw[lightgray] (-2,-13)--(-2,9.3)--(20.3,9.3)--(20.3,-13)--(-2,-13);
\filldraw[lightgray] (9,3)--(11,8)--(14,4)--(12,-1)--(9,3);

\path [draw=black, line width=0.8mm] (-2,-13)--(-2,10)--(21,10)--(21,-13)--(-2,-13);

\path [draw=black, line width=0.4mm] (-6,0)--(-4,5)--(-2,10)--(0,15); %(2,5)
\path [draw=black, line width=0.4mm] (-7,-14)--(-5,-9)--(-3,-4)--(-1,1)--(1,6)--(3,11)--(5,16); %(2,5)

\path [draw=black, line width=0.4mm] (-4,-18)--(-2,-13)--(0,-8)--(2,-3)--(4,2)--(6,7)--(8,12)
--(10,17); %(2,5)
\path [draw=black, line width=0.4mm] (1,-17)--(3,-12)--(5,-7)--(7,-2)--(9,3)--(11,8)--(13,13);
\path [draw=black, line width=0.4mm] (6,-16)--(8,-11)--(10,-6)--(12,-1)--(14,4)--(16,9)--(18,14);
\path [draw=black, line width=0.4mm] (9,-20)--(11,-15)--(13,-10)--(15,-5)--(17,0)--(19,5)--(21,10)
--(23,15); %(2,5)
\path [draw=black, line width=0.4mm] (14,-19)--(16,-14)--(18,-9)--(20,-4)--(22,1)
--(24,6); %(2,5)
\path [draw=black, line width=0.4mm] (19,-18)--(21,-13)--(23,-8)--(25,-3); %(2,5)

\path [draw=black, line width=0.4mm] (-8,-5)--(-5,-9)--(-2,-13)--(1,-17); %(3,-4)
\path [draw=black, line width=0.4mm] (-6,0)--(-3,-4)--(0,-8)--(3,-12)--(6,-16); %(3,-4)
\path [draw=black, line width=0.4mm] (-7,9)--(-4,5)--(-1,1)--(2,-3)--(5,-7)--(8,-11)--(11,-15)--(14,-19); %(3,-4)

\path [draw=black, line width=0.4mm] (-8,18)--(-5,14)--(-2,10)--(1,6)--(4,2)--(7,-2)--(10,-6)--(13,-10)
--(16,-14)--(19,-18);   %(3,-4)

\path [draw=black, line width=0.4mm] (0,15)--(3,11)--(6,7)--(9,3)--(12,-1)--(15,-5)--(18,-9)--(21,-13)
--(24,-17); %(3,-4)
\path [draw=black, line width=0.4mm] (5,16)--(8,12)--(11,8)--(14,4)--(17,0)--(20,-4)--(23,-8)--(26,-12); %(3,-4)
\path [draw=black, line width=0.4mm] (13,13)--(16,9)--(19,5)--(22,1)--(25,-3); %(3,-4)
\path [draw=black, line width=0.4mm] (18,14)--(21,10)--(24,6); %(3,-4)

\path [draw=black, line width=0.4mm] (-9,-19)--(-4,-18)--(1,-17)--(6,-16)--(11,-15)--(16,-14)--(21,-13)
--(26,-12); %(5,1)
\path [draw=black, line width=0.4mm] (-7,-14)--(-2,-13)--(3,-12)--(8,-11)--(13,-10)--(18,-9)
--(23,-8)--(28,-7); %(5,1)
\path [draw=black, line width=0.4mm]  (-5,-9)--(0,-8)--(5,-7)--(10,-6)--(15,-5)--(20,-4)--(25,-3);
\path [draw=black, line width=0.4mm] (-8,-5)--(-3,-4)--(2,-3)--(7,-2)--(12,-1)--(17,0)--(22,1)--(27,2); %(5,1)
\path [draw=black, line width=0.4mm] (-6,0)--(-1,1)--(4,2)--(9,3)--(14,4)--(19,5)--(24,6)--(29,7); %(5,1)

\path [draw=black, line width=0.4mm] (-9,4)--(-4,5)--(1,6)--(6,7)--(11,8)--(16,9)--(21,10)
--(26,11); %(5,1)
\path [draw=black, line width=0.4mm] (-7,9)--(-2,10)--(3,11)--(8,12)--(13,13)--(18,14)--(23,15)--(28,16); %(5,1)

\filldraw[black] (11,7) circle (.35);
\filldraw[black] (11,6) circle (.35);   \filldraw[black] (12,6) circle (.35);  %2
\filldraw[black] (10,5) circle (.35);   \filldraw[black] (11,5) circle (.35);  \filldraw[black] (12,5) circle (.35);
\filldraw[black] (13,5) circle (.35); %4
\filldraw[black] (10,4) circle (.35); \filldraw[black] (11,4) circle (.35); \filldraw[black] (12,4) circle (.35);
\filldraw[black] (13,4) circle (.35); %4
\filldraw[black] (9,3) circle (.35); \filldraw[black] (10,3) circle (.35);  \filldraw[black] (11,3) circle (.35);
\filldraw[black] (12,3) circle (.35); \filldraw[black] (13,3) circle (.35); %5
\filldraw[black] (10,2) circle (.35); \filldraw[black] (11,2) circle (.35); \filldraw[black] (12,2) circle (.35);
\filldraw[black] (13,2) circle (.35); % 4
\filldraw[black] (11,1) circle (.35); \filldraw[black] (12,1) circle (.35); %2
\filldraw[black] (12,0) circle (.35);  %1

%total: 23

\draw[yscale=1, xslant=0] (-90,-60) grid (90, 60);
\draw[yscale=1, yslant=0] (-90,-60) grid (90, 60);

% \filldraw[white] (-0.5,-11.5) circle (1);
% \draw (-0.5,-11.5) node   {{${\mbox{\boldmath$O$}}$}};
\end{tikzpicture}
\caption{\footnotesize{Templates (encircled with thick lines) and FPs (gray parallelograms) for $D^2=25$, with
$S=23$. Here each template
(treated as a torus) has $23\times 23=529$ lattice sites and $23$ FPs,
and each FP  covers $23$ lattice sites (represented by thick black dots). Since the FTs
are not isosceles, there are $4$ types of min-area sub-lattices and $92$ PGSs.}}
\label{Fig7}
\end{figure}

Let $\phi^*$ be a saturation of a given $D$-AC $\phi$. If an added occupied site
$x\in\phi^*\setminus\phi$ lies in a template then, clearly, this template is incorrect (more precisely, non-$\vphi$-correct in $\vphi$ for each $\vphi\in\cP$). We say that such a template is an {\it s-defect} (in $\phi$). Another possibility for a defect is where, in the saturation $\phi^*$, a template has a non-empty intersection with one of C-triangles
that is not an M-triangle. We call it a {\it t-defect} (again in $\phi$).
Finally, an incorrect (but actually perfect) template can be simply a neighbor of an s- or a t-defect. We call it an {\it n-defect} (still in $\phi$).

We would like to note that C-triangles with obtuse angles $>2\pi/3$ (considered in Lemmas 3.4.1--3.4.5) lead to t-defects by definition.

\begin{figure}
\centering
\begin{tikzpicture}[scale=0.06] %\label{%FParaD^2=25}
\clip (-50.6, -40.4) rectangle (69.2, 38.5);
%O: (-2,-13)
\begin{scope}
\definecolor{gray7}{gray}{0.7}
\definecolor{gray5}{gray}{0.5}

\filldraw[gray5] (-2,10)--(21,10)--(21,-13)--(-2,-13)--(-2,10); %central square

%begin the $\vphi$-mesh
\path [draw=black, line width=0.2mm] (-3,5)--(12,11); %(22,15);   %(5,2)
\path [draw=black, line width=0.2mm] (-4,0)--(21,10); %(5,2)
\path [draw=black, line width=0.2mm] (-5,-5)--(25,7);
\path [draw=black, line width=0.2mm] (-6,-10)--(24,2); %(5,2)
\path [draw=black, line width=0.2mm] (-2,-13)--(23,-3); %(28,-1);
\path [draw=black, line width=0.2mm] (7,-14)--(22,-8);  %(5,2)
%\path [draw=black, line width=0.4mm] (16,-15)--(21,-13)--(26,-11);   %(5,2)

%\path [draw=black, line width=0.4mm] (-6,-10)--(-5,-5)--(-1,15); %(1,5)
\path [draw=black, line width=0.2mm] (-2,-13)--(3,12); %(4,17);
\path [draw=black, line width=0.2mm] (2,-16)--(8,14); %(9,19);
\path [draw=black, line width=0.2mm] (7,-14)--(8,-9)--(9,-4)
--(10,1)--(12,11); %(13,16); %(1,5)
\path [draw=black, line width=0.2mm] (11,-17)--(12,-12)--(13,-7)
--(14,-2)--(15,3)--(17,13); %(1,5)
\path [draw=black, line width=0.2mm] (16,-15)--(17,-10)--(18,-5)--(19,0)--(20,5)
--(21,10); %(1,5)
%\path [draw=black, line width=0.4mm] (20,-18)--(21,-13)--(22,-8)--(23,-3)
%--(24,2)--(25,7)--(26,12)--(27,17); %(1,5)

%\path [draw=black, line width=0.4mm] (-2,-13)--(2,-16); %(4,-3)
\path [draw=black, line width=0.2mm] (-5,-5)--(-1,-8)--(3,-11)--(7,-14);
\path [draw=black, line width=0.2mm] (-4,0)--(0,-3)--(8,-9)--(12,-12)
--(16,-15); %(4,-3)
\path [draw=black, line width=0.2mm] (-3,5)--(1,2)--(5,-1)--(9,-4)--(13,-7)
--(17,-10)--(21,-13); %(4,-3)
\path [draw=black, line width=0.2mm] (-2,10)--(2,7)--(6,4)--(10,1)
--(14,-2)--(18,-5)--(22,-8); %(4,-3)
\path [draw=black, line width=0.2mm] (3,12)--(7,9)--(11,6)--(15,3)--(23,-3);
\path [draw=black, line width=0.2mm] (12,11)--(16,8)--(20,5)--(24,2);
%\path [draw=black, line width=0.4mm] (13,16)--(17,13)--(21,10); %(4,-3)
%end the $\vphi$-mesh

\filldraw[white] (-2,-13)--(-2,-36)--(21,-36)--(21,-13)--(-2,-13);
\filldraw[white] (-2,10)--(21,10)--(21,33)--(-2,33)--(-2,10);
\filldraw[white] (-25,-13)--(-25,10)--(-2,10)--(-2,-13)--(-25,-13);
\filldraw[white] (21,-13)--(44,-13)--(44,10)--(21,10)--(21,-13);

\filldraw[gray7] (-2,10)--(21,10)--(21,33)--(-2,33)--(-2,10);          %left column of lightgray squares
\filldraw[gray7] (-2,-13)--(-2,-36)--(21,-36)--(21,-13)--(-2,-13);
\filldraw[gray7] (-2,-13)--(-2,-36)--(-25,-36)--(-25,-13)--(-2,-13); %left column of lightgray squares

\filldraw[gray7] (-25,-13)--(-25,10)--(-2,10)--(-2,-13)--(-25,-13);  %right column of lightgray squares
\filldraw[gray7] (21,-13)--(44,-13)--(44,10)--(21,10)--(21,-13);
\filldraw[gray7] (21,-13)--(44,-13)--(44,-36)--(21,-36)--(21,-13); %right column of lightgray squares

\filldraw[gray7] (-2,10)--(-2,33)--(-25,33)--(-25,10)--(-2,10);  %middle column of lightgray squares
\filldraw[gray7] (21,10)--(44,10)--(44,33)--(21,33)--(21,10);  %middle column of lightgray squares

\path [draw=black, line width=0.8mm] (-25,-36)--(-25,-13)--(-48,-13)--(-48,-36)--(-25,-36);
\path [draw=black, line width=0.8mm] (-2,-36)--(-2,-13)--(-25,-13)--(-25,-36)--(-2,-36);
\path [draw=black, line width=0.8mm] (-2,-13)--(-2,-36)--(21,-36)--(21,-13)--(-2,-13);
\path [draw=black, line width=0.8mm] (21,-13)--(21,-36)--(44,-36)--(44,-13)--(21,-13);
\path [draw=black, line width=0.8mm] (44,-13)--(44,-36)--(67,-36)--(67,-13)--(44,-13);

\path [draw=black, line width=0.8mm] (-48,-13)--(-48,10)--(-25,10)--(-25,-13)--(-48,-13);
\path [draw=black, line width=0.8mm] (-25,-13)--(-25,10)--(-2,10)--(-2,-13)--(-25,-13);
\path [draw=black, line width=0.8mm] (-2,-13)--(-2,10)--(21,10)--(21,-13)--(-2,-13);
\path [draw=black, line width=0.8mm] (21,-13)--(44,-13)--(44,10)--(21,10)--(21,-13);
\path [draw=black, line width=0.8mm] (44,-13)--(44,10)--(67,10)--(67,-13)--(44,-13);

\path [draw=black, line width=0.8mm] (-25,33)--(-48,33)--(-48,10);
\path [draw=black, line width=0.8mm] (-2,10)--(-2,33)--(-25,33)--(-25,10)--(-2,10);
\path [draw=black, line width=0.8mm] (-2,33)--(21,33);
\path [draw=black, line width=0.8mm] (21,10)--(21,33)--(44,33)--(44,10)--(21,10);
\path [draw=black, line width=0.8mm] (44,33)--(67,33)--(67,10);

\path [draw=black, line width=0.8mm] (-25,33)--(-25,39);
\path [draw=black, line width=0.8mm] (-2,33)--(-2,39);
\path [draw=black, line width=0.8mm] (21,33)--(21,39);
\path [draw=black, line width=0.8mm] (44,33)--(44,39);

\path [draw=black, line width=0.8mm] (-25,-36)--(-25,-40);
\path [draw=black, line width=0.8mm] (-2,-36)--(-2,-40);
\path [draw=black, line width=0.8mm] (21,-36)--(21,-40);
\path [draw=black, line width=0.8mm] (44,-36)--(44,-40);

\filldraw[black] (0,-3) circle (.8); \filldraw[black] (1,2) circle (.8);  \filldraw[black] (2,7) circle (.8);

\filldraw[black] (5,-1) circle (.8);

\filldraw[black] (10,1) circle (.8);  \filldraw[black] (15,3) circle (.8);
\filldraw[black] (14,-2) circle (.8); \filldraw[black] (16,8) circle (.8);
\filldraw[black] (9,-4) circle (.8); \filldraw[black] (13,-7) circle (.8);
\filldraw[black] (8,-9) circle (.8);
\filldraw[black] (17,-10) circle (.8); \filldraw[black] (18,-5) circle (.8);
\filldraw[black] (19,0) circle (.8);

\filldraw[black] (0,-11) circle (.8);  \filldraw[black] (2,-8) circle (.8);
\filldraw[black] (8,6) circle (.8);

%\draw[yscale=1, xslant=0] (-90,-60) grid (90, 60);
%\draw[yscale=1, yslant=0] (-90,-60) grid (90, 60);

\filldraw[white] (-2,-13) circle (.8); \filldraw[white] (3,-11) circle (.8);
\filldraw[white] (-1,-8) circle (.8); \filldraw[white] (4,-6) circle (.8);

\filldraw[white] (6,4) circle (.8);  \filldraw[white] (7,9) circle (.8);
\filldraw[white] (11,6) circle (.8);

\filldraw[white] (20,5) circle (.8);  \filldraw[white] (12,-12) circle (.8);

%\draw (-0.5,-11.5) node   {\Large {${\mbox{\boldmath$O$}}$}};
\end{scope}
\end{tikzpicture}

\caption{\footnotesize{A snapshot of templates for $D^2=25$ with $S(D)=23$. All squares are $23\times 23$
and represent templates for this value of $D$. The white squares indicate
$\vphi$-correct templates in a given $D$-AC $\phi$ while gray and light-gray squares
indicate non-$\vphi$-correct ones.
The original lattice $\bbZ^2$ is not shown due to its small scale.
The {\rm{PGS}} $\vphi$ is represented by the mesh over the central (gray) square.
The black dots indicate the occupied $\bbZ^2$-sites over the central template square.
The white dots indicate vacant sites in $\vphi$  in the central square. All light-gray and
white common squares indicated in the diagram
are assumed to preserve the structure of the {\rm{PGS}} $\vphi$ unperturbed. However, the $8$ light-gray
templates around the central square are declared $\vphi$-non-correct by adjacency.}}
\label{Fig8}
\end{figure}
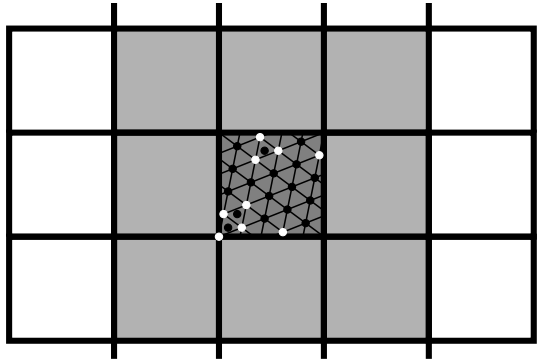

\begin{figure}
\centering
\begin{tikzpicture}[scale=0.25]
\clip (-8, -4.5) rectangle (30, 27);

\begin{scope}
\definecolor{gray3}{gray}{0.3}
\definecolor{gray5}{gray}{0.5}
\path [fill=gray5, draw=black, very thick] (2, 6) -- (5, 6) -- (5, 9) -- (4, 9) -- (4, 11) -- (1, 11) -- (1, 8) -- (2, 8) -- cycle;
\path [fill=gray3, draw=black, thick] (3,7) -- (4, 7) -- (4, 8) -- (3, 8) -- cycle;
\path [fill=gray3, draw=black, thick] (2,9) -- (3, 9) -- (3, 10) -- (2, 10) -- cycle;
\path [fill=gray5, draw=black, very thick] (5,3) -- (8,3) -- (8,4) -- (11,4) -- (11,6) -- (14,6) -- (14,9) -- (15,9) --
(15,12) -- (16,12) -- (16,13) -- (18,13) -- (18,12) -- (19,12) -- (19,6) -- (14,6) -- (14,3) -- (15,3) -- (15,2) -- (20,2) --
(20,3) -- (23,3) -- (23,6) -- (22, 6) -- (22,14) -- (21,14) -- (21,15) -- (20,15) -- (20,16) -- (16,16) -- (16,21) -- (15,21) --
(15,22) -- (14,22) -- (14,23) -- (13,23) -- (13,25) -- (12,25) -- (12,26) -- (2,26) -- (2, 16) -- (6,16) -- (6,14) -- (4,14) --
(4,11) -- (7,11) -- (7,10) -- (8,10) -- (8,9) -- (9,9) -- (9,8) -- (10,8) -- (10,7) -- (8,7) -- (8,6) -- (5,6) -- cycle;
\path [fill=gray3, draw=black, thick] (6,4) -- (7,4) -- (7,5) -- (6,5) -- cycle;
\path [fill=gray3, draw=black, thick] (9,5) -- (10,5) -- (10,6) -- (9,6) -- cycle;
\path [fill=gray3, draw=black, thick] (5,12) -- (6,12) -- (6,13) -- (5,13) -- cycle;
\path [fill=gray3, draw=black, thick] (15,4) -- (16,4) -- (16,3) -- (19,3) -- (19,4) -- (22,4) -- (22,5) -- (21,5) -- (21,13)--
(20,13) -- (20,14) -- (19,14) -- (19,15) -- (15,15) -- (15,20) -- (14,20) -- (14,21) -- (13,21) -- (13,22) -- (12,22) -- (12,24)--
(11,24) -- (11,25) -- (3,25) -- (3,17) -- (7,17) -- (7,12) -- (8,12) -- (8,11) -- (9,11) -- (9,10) -- (10,10) -- (10,9) --
(11,9) -- (11,8) -- (12,8) --(12,7) -- (13,7) -- (13,10) -- (14,10) -- (14,11) -- (12,11) -- (12,10) -- (10,10) -- (10,11) --
(9,11) -- (9,12) -- (8,12) -- (8,13) -- (8,18) -- (4,18) -- (4,24) -- (9,24) -- (9,21) -- (13,21) -- (13,20) -- (14,20) --
(14,16) -- (13,16) -- (13,13) -- (15,13) -- (15,14) -- (19,14) -- (19,13) -- (20,13) -- (20,5) -- (17,5) -- (17,4) -- (16,4)--
(16,5) -- (15,5) -- cycle;
\path [fill=white, draw=black, very thick] (10,12) -- (12,12) -- (12,17) -- (13,17) -- (13,19) -- (12,19) -- (12,20) -- (8,20) --
(8,23) -- (5,23) -- (5,19) -- (9,19) -- (9,13) -- (10,13) -- cycle;

\draw (-10,0) grid (34,28);

%\path [fill=white, draw=white] (6,21) -- (9.5,21.5) -- (9.5,21.5) -- (10,21) -- (9.5,20.5) -- (9.5,20.5) -- (6,21) -- cycle;

\draw (-2.5, 21) node   {\Large {{\bf{Ext}}({\mbox{\boldmath$\Gamma$}})}};
\end{scope}
\end{tikzpicture}
\caption{\footnotesize{A contour support (the union of gray and dark gray common squares). Here
the internal area ${\rm{Int}}(\Gamma)$ includes three components ${\rm{Int}}_{\varphi_i}(\Gamma)$,
$i=1,2,3$. The boundary layers are shown as the union of gray squares.
For a $\vphi$-contour, the configuration over every white common square forming the external
area ${\rm{Ext}}\,(\Gamma)$ is  supposed to be the restriction of $\vphi$.
}}
\label{Fig9}
\end{figure}
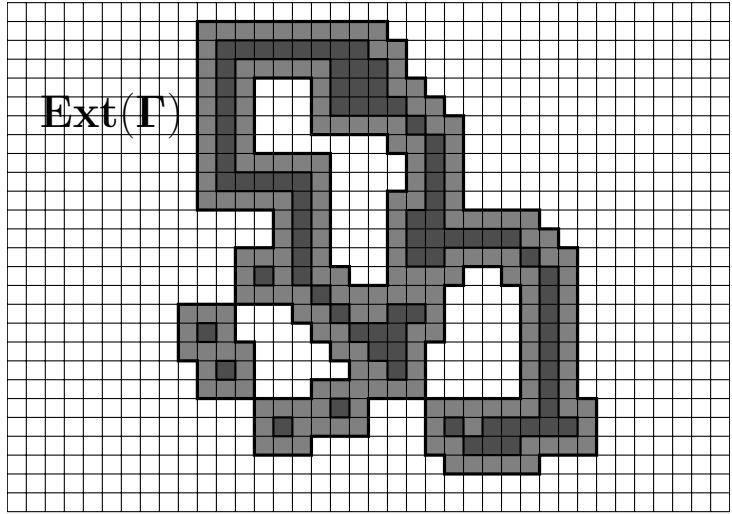

%{\bf Lemma 3.8.}
\bli {\rm{(A Peierls bound for defects.)}}
Consider a $D$-AC $\phi\in\cA$ and
assume that for configuration $\phi$ there exists a contour $\Gam$ containing $m$
incorrect templates enclosed by $n$ adjacent $\vphi$-correct templates (for the same
or different values of $\vphi$). Additionally, assume that $m = i + j + k$ where $i,j,k$
give the amount of s-, t- and n-defects, respectively. Then the amount of occupied sites inside
$\Gam$ is at most
$$(m+n)S(D)-i-\max\left(1, \frac{j}{8S(D)}\right).\eqno (3.5)$$
\eli

\bp
As the contribution $i$ given by s-defects is straightforward, we
consider the saturation $\phi^*$ and its t-defects only. We confine the entire
connected component of templates to a sufficiently large square
$\bbQ_{kS(D)}(0,0)$ by filling all additional $k^2 - m - n$ templates with the
appropriate PGS. Here the appropriate means that for each added $\vphi$-correct
template all neighboring templates (which by construction are correct) have
the same value of $\vphi$. Afterwards we wrap $\bbQ_{kS(D)}(0,0)$ into
the $kS(D)\times kS(D)$ torus $\bbT_{kS(D)}$ and still obtain a
$D$-admissible configuration in $\bbT_{kS(D)}$.

Due to Lemma 3.2, for any saturated configuration in a torus the total amount
of corresponding
C-triangles inside this torus is twice the amount of the occupied sites in the
configuration. In particular, the total amount of C-triangles is always even. Also the
existence of at least one non-$D$-optimal
C-triangle inside the torus reduces the maximal amount of occupied sites by at
least one.

Observe that, according to Lemma 3.4.5,  the re-distributed area of any C-triangle
that is not an M-triangle is at least $(S(D)+\del (D))/2$, where
$\del (D) \ge 1$. Further, a C-triangle that is not an M-triangle can be shared by at most
4 templates. Therefore, $j$ templates with
t-defects contain at least $j/4$ C-triangles that are not M-triangles. Thus, one half of
the maximal possible
amount of C-triangles inside $T_{kS(D)}$ is $\diy k^2S(D) - \max \left(1, \frac{j}{8
S(D)}\right)$.

The term $(m+n)S(D)$ corresponds to the absence of defects and is calculated
according to Lemma 3.6.
\ep

Informally, Lemma 3.8 states that the increment of `energy' (i.e., decrease in the number
of particles) caused by a deviation from a PGS
is lower-bounded proportionally to the `size' of the deviation. This is the gist of Peierls
bounds used in the Pirogov--Sinai theory and its applications.

%{\bf Lemma 3.9.}
\blj
Let $\vphi',\vphi''\in{\cP}(D)$ be two distinct
PGSs. Let \ $\Gam$ \ be a connected component of $\vphi'$-correct templates
enclosed by a connected component of  $\vphi''$-correct templates. Then any extension
of this restricted configuration to a
$D$-AC $\phi\in{\mathcal A}(D)$ in $\bbZ^2$ contains a closed chain of adjacent
non-minimal C-triangles enclosing \ $\Gam$.
\elj

\bp Due to the absence of sliding, M-triangles from different classes
cannot share a side in a $D$-AC.
\ep

%{\bf Lemma II.}
\blII {\rm{(A Peierls estimate for contours)}}
Consider a finite contour containing
$m$ incorrect templates. Then the amount of occupied sites inside this contour is at
most
$$m S(D) -\max \left(1,\; \frac{ m}{72 S(D)} \right).\eqno (3.6)$$
\elII

\bplemII %{\bf Proof of Lemma II.}
The lemma is a direct consequence of Lemmas 3.8~and 3.9 with an
additional factor
$1/9$ accounting for the possibility for each s- or t-defect to be surrounded by 8
n-defects.
\eplemII

Lemma II completes the verification of assumptions of the PS theory for HC models
in $\bbZ^2$.

%{\bf Remark 3.1.}
\bra {\rm The Peierls bound for the HC model on $\bbA_2$ is considerably
simpler and is stated in terms of Voronoi cells.}  \quad $\blacktriangle$
\era

\subsection{The proof of Theorems 1, 2.} With Lemmas I, II at hand, the proof of
Theorems 1, 2 is obtained by a direct application of the PS-theory. Cf. \cite{PiS,
Z, DoS} and \cite{Si}.

We refer the
reader to \cite{MSS1} for the scheme of such an application.

\section{Additional results. The issue of dominance}\label{Sec4}

\subsection{More on extreme Gibbs measures.}\label{SubSec4.1} As an annex to Theorem 1, based
on straightforward applications of the
approach from \cite{Z}, we offer three more statements in Theorem 3. Consider
a sequence of $nS(D)\times nS(D)$ tori $\bbT^{(n)}=\bbT_{nS(D)}$. Let $\cA(D,\bbT^{(n)})$ denote the set of $D$-admissible configurations in $\bbT^{(n)}$:
$$\beal\cA (D,\bbT^{(n)})=\Big\{\psi :\bbQ_{2nS(D)}(0,0)\to\{0,1\},\qquad{}\\
\qquad{}\psi (x)\psi(y)=0\;\;\forall\;
x,y\in\bbQ_{nS(D)}(0,0)\;\hbox{ with }\;\rho^{\;\rT}_{nS(D)} (x,y)< D\Big\}.\ena\eqno (4.1)$$

Let $\mu_n^{\rm{per}}=\mu_{\bbT^{(n)}}$ denote the probability measure on $\cA (D,\bbT^{(n)})$
with
$$\mu_{\bbT^{(n)}}(\psi )=\frac{u^{\sharp\,\psi}}{\BZ (\bbT^{(n)})},\;\;\psi\in \cA (D,\bbT^{(n)}),\eqno (4.2)$$
where the partition function $\BZ(\bbT^{(n)} )=\sum\limits_{\psi\in\cA(\bbT^{(n)})}u^{\sharp\,\psi}$.
Then $\mu_{\bbT^{(n)}}$ can be treated as a measure on $\cA(D,\bbZ^2)$ concentrated on
$\cA (D,\bbT^{(n)})$.

%{\bf Theorem 3.}
\bthmd\label{Theorem4}
Under assumptions of Theorem $1$ or Theorem $2$,
%$\exists$ $u_*\in (0,\infty)$ such that for $u\geq u_*$
the following
properties {\rm{(i)}}, {\rm{(ii)}} are satisfied:
\begin{description}
  \item[(i)] For every extreme measure $\bmu_\vphi\in{\mathcal E}(D)$ $\exists$ a set
$\sI=\sI(\vphi )\subset{\mathcal A}(D)$ with $\vphi\in\sI$ and
$\bmu_\vphi (\sI)=1$ such that the set of $\vphi$-incorrect templates in any $\phi\in\sI$ has
no connected components of infinite diameters in $\bbR^2$. Consequently,
 {\rm{(a)}} $\forall$ $\phi\in\sI$ the set of $\vphi$-
correct templates in $\phi$ has an infinite connected component, and {\rm{(b)}}
$\vphi$ is the only $D$-AC with this property.
  \item[(ii)] Measures $\bmu_\vphi$ have an exponential decay of spatial correlations:
$\forall$ \ $l>1$, finite sets $\bbV_{1,2}\subset\bbZ^2$ with $\rho (\bbV_1,\bbV_2)\geq l$
and \ $\forall$ \ $D$-ACs
$\psi^{(\bbV_i)}\in{\mathcal A}(D,\bbV_i)$, $i=1,2$,
$$\left|\bmu_\vphi (\psi^{(\bbV_1)}\vee\psi^{(\bbV_2)})-\bmu_\vphi (\psi^{(\bbV_2)})
\bmu_\vphi (\psi^{(\bbV_2)})\right|\leq \exp (-Cl)\bmu_\vphi (\psi^{(\bbV_2)})
\bmu_\vphi (\psi^{(\bbV_2)})\eqno (4.3)$$
Here $C=C(D)>0$ is a constant and $\bmu (\psi^{(\bbV_1)}\vee\psi^{(\bbV_2)})$ and
$\bmu (\psi^{(\bbV_i)})$ stand for the $\bmu_\vphi$-probabilities of the cylinder
events that $\phi\upharp_{\bbV_1\cup\bbV_2}=\psi^{(\bbV_1)}\vee\psi^{(\bbV_2)}$
and $\phi\upharp_{\bbV_i}=\psi^{(\bbV_i)}$, respectively.

Furthermore, assume that all  extreme measures $\bmu_\vphi$ are
generated by PGSs $\vphi$ belonging to a single equivalence class. (In particular, this covers
the case of Theorem $1$.) Then
%$\exists$ $u_*\in (0,\infty)$ such that for $u\geq u_*$
the following property is satisfied.
  \item[(iii)] The weak limit $\bmu^{\rm{per}}
=\lim\limits_{n\to\infty}\mu_{\bbT^{(n)}}$ exists and determines a probability measure
on $\sX$ which is the averaged sum: $\bmu^{\rm{per}}=\diy\frac{1}{\sharp\cE (D)}
\sum\limits_{\bmu_\vphi\in\cE(D)}\bmu_\vphi$. In particular, $\bmu^{\rm{per}}$ is shift-invariant.
\end{description}
\ethmd

\bp %{\bf Proof of Theorem 3.}
Assertion (i) is a standard corollary of the Pirogov--Sinai
theory \cite{PiS,Si,Z,BrS}: it implies that (1) the $\bmu_\vphi$-probability that a given site
$x\in\bbZ^2$ is encircled by infinitely many contours equals $0$ and (2) the $\bmu_\vphi$-probability
that there is an infinite contour equals $0$ as well.

(ii) This statement follows from standard polymer expansions for probabilities\\
$\bmu (\psi^{(\bbV_1)}\vee\psi^{(\bbV_2)})$ and $\bmu (\psi^{(\bbV_i)})$; cf. \cite{Z}. The polymers
are determined as collections of overlapping contours, and the statistical weight of a polymer
decays exponentially with the polymer size.

(iii) Observe that, under the stated assumption, $\forall$ PGSs $\vphi_1,\vphi_2\in{\cP}(D)$
such that $\mu_{\vphi_1},\mu_{\vphi_2}\in\cE(D)$,
there exists a 1-1 map $\Phi_{\vphi_1,\vphi_2}:\bbZ^2\to\bbZ^2$ (a $\bbZ^2$-congruence)
taking $\vphi_1\to\vphi_2$ such that $\bmu_{\vphi_1}$ is taken to $\bmu_{\vphi_2}$
by the induced map ${\mathcal A}(D)\to{\mathcal A}(D)$. On the other hand, $\forall$ PGSs
$\vphi_1,\vphi_2\in{\cP}(D)$ as above, the measure $\mu_{\bbT^{(n)}}$
is invariant under the map induced by $\Phi_{\vphi_1,\vphi_2}$ restricted to $\bbT^{(n)}$.
Hence, any limit point for the sequence $\mu_{\bbT^{(n)}}$ is a uniform mixture of measures
$\bmu_\vphi$. Hence, the limiting measure $\bmu^{\rm{per}}$ is a uniform mixture,
as claimed.
\ep

%{\bf Remark 4.1.}
\br {\rm A more detailed analysis of dependence of Gibbs measures $\bmu$ upon
boundary conditions requires a deeper technical involvement and will be given elsewhere.
The same can be said about a similar question on the HC model in $\bbA_2$ and $\mathbb{H}_2$. Cf. \cite{MSS1}.}
 \quad $\blacktriangle$
\er

\subsection{A discussion of dominance.}\label{SubSec4.2} This section contains some discussions and comments but no formally proven results. First, as was said earlier,
we conjecture that, under the assumptions of Theorem 2, there is only one
dominant class. The reason is that co-existence of several dominant classes requires
countably many `balance identities' of a number-theoretical nature,
in each order of an emerging perturbation series. In other words, it would require
some (rather strong) form of equivalence between contributions from local excitations
of a given weight coming from two or more distinct classes.
It seems unlikely since the definition of an equivalence {\rm{PGS}} class probably captures
all involved forms of symmetry in the hard-core model on $\bbZ^2$. The identification of the dominant PGS class is an algorithmically finite procedure. According to \cite{Z}, one needs to iteratively calculate truncated statistical weights of contours and the corresponding approximations to the free energies of the emerging truncated contour models. The procedure stops as soon as it detects a single class having the lowest free energy. The amount of calculations required by the direct approach is not accessible to modern computers. In cases where the above program was carried out to the end, one needed a number of tricks (both analytical and computational), as can be seen in \cite{MSS1}.

As an illustration, we briefly discuss the cases $D^2=65$, $D^2=130$ and $D^2=324$:
these are the lowest values of $D$ from Class B. The algebraic nature of these values
is commented upon in Section 6--8; here we focus on the structure of local excitations.

The numerical values presented below have been calculated by {\tt CountTuples.java}.%; cf. Section \ref{program} and the supplement to this paper.

For $D^2=65=8^2+1^2=7^2+4^2$, we have $S(D)=60$. There are two
minimizing triangle types. The corresponding
side-length triples from (2.2) have $\ell_0^2=\ell_1^2=65$, $\ell_2^2=80$ and
$\ell_0^2=\ell_1^2=68$, $\ell_2^2=72$. It is convenient to say that one type
gives a $[65|65|80]$-triangle and the other a $[68|68|72]$-triangle. (This notation and
terminology will be further elaborated in Sections 6, 7.) A $[68|68|72]$-triangle
is also an M-triangle for a Class A value $D^2=68$ (cf. Table 3).

Both $[65|65|80]$- and $[68|68|72]$-types admit a single congruence class or
a single implementation. (In Section 7.2 we speak of Class B1 in more general
terms and provide details of a number-theoretical background.)
For type $[65|65|80]$, a representative is an M-triangle with vertices $\{(0,0),(8,4),(-1,8)\}$,
while for type $[68|68|72]$ the vertices are $\{(0,0), (8,2), (2,8)\}$. Recall, an M-triangle gives
rise to a class of PGSs obtained as $\bbZ^2$-shifts and reflections of the corresponding
min-area sub-lattice.

As was said earlier, the dominance of a particular {\rm{PGS}} class is established by
comparing the number of local excitations for both minimizing types. (Actually, what we 
compare is the naturally emerging densities of excitations.) It is convenient to categorize 
the local excitations into `strains' bearing a distinctive
geometric character. The initial excitation strain is where we simply remove one particle
from a PGS; this yields an excitation of (relative) weight $u^{-1}$. However, all PGSs are
`equal in rights' with respect to this strain of excitations as the particle density in all PGSs
is the same. In other words, there is no discrimination between the PGSs in the
perturbation order $u^{-1}$. A similar conclusion can be reached when one considers
removing a pair of particles from a {\rm{PGS}} which produces a strain of excitations of weight
$u^{-2}$.

Next, we can remove three particles at vertices of a fundamental parallelogram $\Pi$ in
a max-dense sub-lattice and add one
particle inside $\Pi$ maintaining admissibility. We call this strain a single or 1-site insertion.
This produces another strain of local excitations of weight $u^{-2}$ (3 occupied
sites removed, 1 added). Another strain of excitations of weight $u^{-2}$ is a double or
2-site insertion
where we add 2 particles in an FP $\Pi$  and remove 4 occupied sites at the vertices of
the $\Pi$. Next,  we can attempt to add 3 particles and remove 5, then add 4 and remove 6,
and so on. The PGSs must be compared by counting all strains of excitations per an FP;
the ones with a larger total count will be dominant.

To complete the argument, we must check that there is no other possibilities to produce
a local excitation of weight $u^{-2}$ (which is the most difficult part) and prove that other
strains of excitations (of higher weights) can be made `insignificant' if $u$ is large enough
(which is usually not a difficult part).  Such a program has been
carried through for a triangular lattice $\bbA_2$ in \cite{MSS1} with the help of an
analytical construction and a computer-assisted argument. For the lattice $\bbZ^2$
it is a work in progress which will be published separately: here we only attempt to explain the
place of the dominance argument in the whole study and comment on key points of
the forthcoming argument.

The count of excitations of weight $u^{-2}$ yields these results: there are $40$ single
insertions per an FP for both types. Next, there are  $109$ double insertions per an FP
for type  $[65|65|80]$ and $126$ for type $[68|68|72]$. Further, there are
$104$ triple insertions and no quadruple ones for type $[65|65|80]$ whereas type $[68|68|72]$
offers $140$ triple insertions  and $10$. Finally, no insertions of $5$ or more particles of
weight $u^{-2}$ exist for either type.  The overall score is $253$
for $[65|65|80]$ and 316 for $[68|68|72]$. The tentative conclusion is that type
$[68|68|72]$, with an M-triangle with vertices
$\{(0,0),(8,2),(2,8)\}$, is dominant. Consequently, for $u$ large enough, we expect that
the number of the extreme Gibbs measures for $D^2=65$ is $120$, and they all are
generated by the PGSs coming from $[68|68|72]$-triangles. %, in accordance with Eqns (2.4.1), (2.4.2).
Cf. Figure \ref{Fig10Dom65}.

Now, consider the value $D^2=130=11^2+3^2=9^2+7^2$, with $S(D)=120$. Here we
again have two minimizing types. They are represented by M-triangles with
squared side-lengths (I) $\ell_0^2=\ell_1^2=130$, $\ell_2^2=160$ and (II)
$\ell_0^2=\ell_1^2=136$, $\ell_2^2=144$, respectively, both isosceles. We again refer
to them as $[130|130|160]$- and  $[136|136|144]$-types or triangles.  It turns out that a
$[136|136|144]$-triangle also serves as the M-triangle for the value $D^2=136$, which is
of Class A.

As above, both types admit a single implementation.  For type $[130|130|160]$,
a representative is an M-triangle with vertices $\{(0,0),(12,4), (3,11)\}$, while for type
$[136|136|144]$ the vertices are $\{(0,0),(12,0),(6,10)\}$. In fact, it is easy to see that
the picture of M-triangles (or min-area sub-lattices) for $D^2=130$ is obtained from
that for $D^2=65$ by applying the linear transformation ${\tt T}:\;\bbZ^2\to\bbZ^2$
effectuated by the matrix $\diy{\mathbf T}=\begin{pmatrix}1&-1\\1&1\end{pmatrix}$ of
determinant $2$. %(The transposed matrix ${\mathbf T}^{\rm T}$ can also be used for that purpose.) %Cf. Remark 8.1 in Section {\color{red} 7}. 
In other words, the case $D^2=130$ is
obtained from $D^2=65$ by scaling
the side-length by a factor ${\sqrt 2}$ and rotating by $\pm\pi/4$. However, it does not
guarantee that all excitations are mapped to each other; in fact we will see that it is not the
case, although type $[136|136|144]$ again will look dominant.

Indeed, type $[130|130|160]$ generates $72$ single insertions, $303$ double insertions,
$284$ triple insertions and  no quadruple insertions per an FP, 659 in total. Further, for type\\
$[136|136|144]$, the corresponding numbers are $72$,  $356$, $532$ and
$50$, all-in-all $1010$. As above, none of the types produces admissible insertions PGS
of $5$ or more particles of weight $u^{-2}$. Then, for $u$ large enough, the number of the EGMs is $240$, and they are generated by the
PGSs of type $[136|136|144]$. Cf. Figure \ref{Fig11dom130}. 

The next case is $D^2=324$, with $S=288$. In Section 6 we show that this case
gives rise to an infinite sequence of values of $D^2$ of Class B1. Here we have two
minimizing types, $[324|337|337]$ (isosceles), and
$[325|333|340]$ (non-isosceles). For congruence class representatives, one can take M-triangles
with vertices $\{(0,0),(18,0),(9,16)\}$ for type $[324|337|337]$ and $\{(0,0),(18,3),(6,17)\}$ for
$[325|333|340]$. Here $[325|333|340]$-triangles give rise to $1152$ PGSs, twice as much as
$[324|337|337]$-ones.  Cf. Figure \ref{Fig12dom324}. The total count of excitations of weight
$u^{-2}$ yields $3797$ per an FP for type $[324|337|337]$ and $2684$ for
$[325|333|340]$. As a result, $[324|337|337]$-triangles are expected to win, and the number
of EGMs for $u$ large will be $576$, all of them being generated by the respective PGSs.

\begin{figure}[H]
\centering
A

\caption{\footnotesize{Excitation counts for $D^2=65$, with $S=60$. The crucial excitations
are those of weight $u^{-2}$. Frame A shows the $65\times 65$ squares, for $65=8^2+1^2$
and $65=7^2+4^2$. In frame B we show the M-triangles of type $[65|65|80]$ (left)
and $[68|68|72]$ (right) and their FPs. The gray regions indicate
the single-insertion positions repelling $3$ sites in a PGS (black circles); the white marks and black segments outline examples of such
excitations. The total number of single insertions per an FP equals 40 for both types.
Frame C shows double insertions (pairs of white marks connected with white lines)
repelling $4$ sites (the vertices of an FP). The total number of double insertions per
an FP equals $109$ and $126$ for types $[65|65|80]$ and  $[68|68|72]$, respectively.
Frame D  shows triple insertions (white marks and triangles) repelling $5$ sites (the
vertices of 3 adjacent M-triangles forming a trapezoid). The total number of triple insertions
per an FP equals $104$ for type $[65|65|80]$ and $140$ for $[68|68|72]$. In frames E,F we
attempt at quadruple insertions repelling $6$ sites (they should be the vertices of $6$
adjacent M-triangles forming a triangle of a double size). Their number of such insertions
is $0$ and $10$ per an FP, respectively.  In both cases there is no $u^{-2}$-excitations
with $5$ or more insertions. Total total of excitations of weight $u^{-2}$ amounts to $253$ for
$[65|65|80]$ and $316$ for $[68|68|72]$; the latter is expected to be the dominant one.
So, for $D^2=65$ the number of extreme Gibbs measures for $u\geq u_*$ equals 120, and 
they are generated by min-area sub-lattices of type $[68|68|72]$.}}
\label{Fig10Dom65}
\end{figure}

\begin{figure}[H]
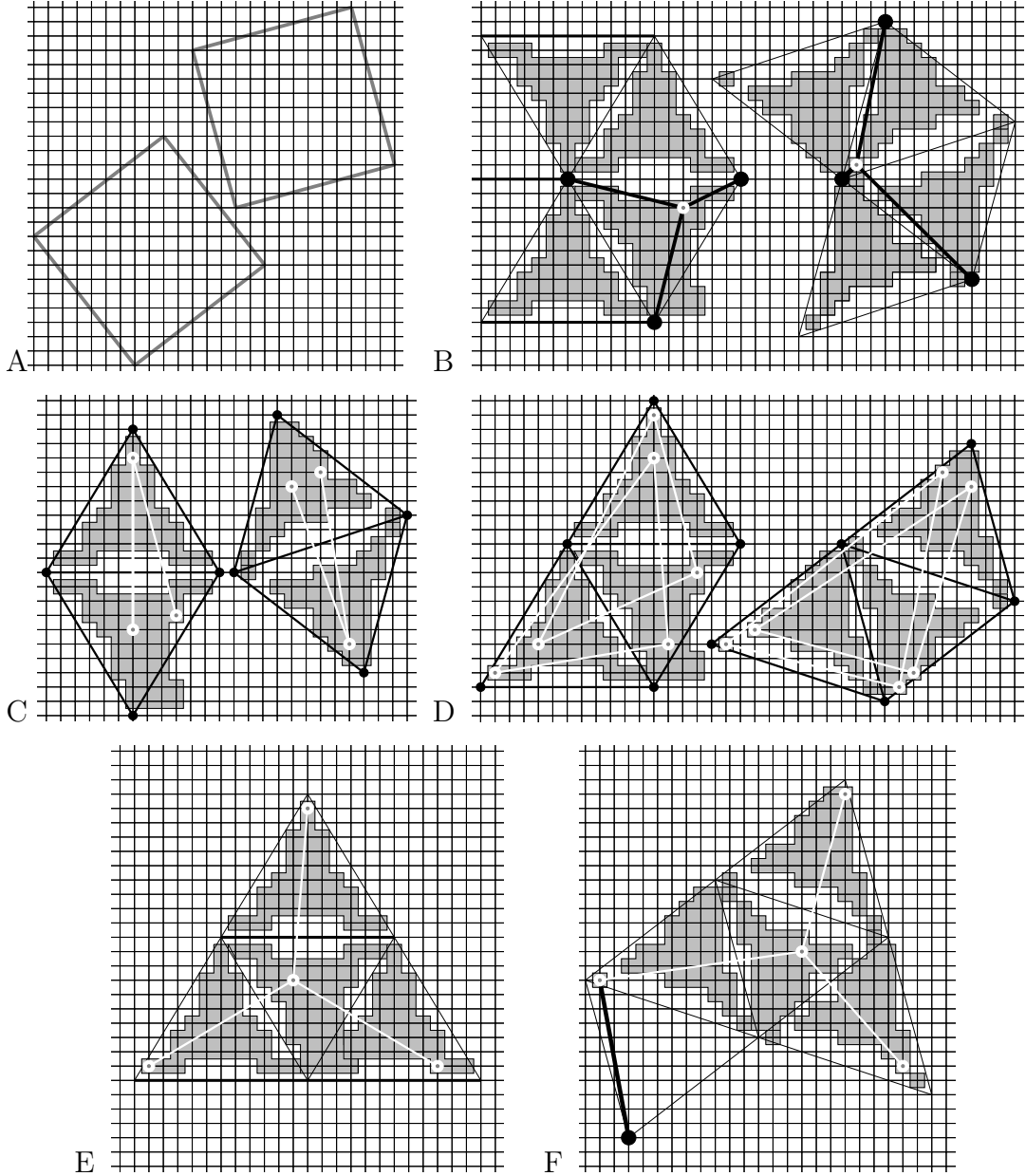

\centering
A

\caption{\footnotesize{Excitation counts for $D^2=130$, with $S=120$. As before, we look
for excitations of weight $u^{-2}$. Frame A shows the $130\times 130$
squares, for $130=11^2+3^2$ and $130=9^2+7^2$. In
frame B we show the M-triangles $[136|136|144]$ (left) and
$[130|130|160]$(right) and their fundamental parallelograms;
we also display 1-site insertions. Here again, the
gray regions indicate the single-site insertions excluding $3$ vertices in a PGS;
the white marks and black lines show examples of such excitations. The total
number of single insertions in an FP equals $72$ for both type. Frame
C shows double insertions (pairs of white marks connected with lines)
repelling $4$ sites (as above, -- vertices of an FP). The total number of double
insertions within an FP equals $356$ and $303$ for
$[136|136|144]$- and $[130|130|160]$-types, respectively.
Frame D  shows triple insertions (white marks forming triangles) which repel $5$
sites placed on a trapezium. The total number of such excitations per an FP equals
$532$ and $284$, respectively. On frames E,F we comment on quadruple insertions
that repel $6$ vertices; their numbers are $50$ per an FP and $0$, respectively.
Indeed, an attempt to construct a 4-site excitation in a $[130|130|160]$-type PGS
in frame F leads to a repulsion of 7 vertices, hence to the weight $u^{-3}$. As with
$D^2=65$, neither type for $D^2=130$ yields  an admissible $u^{-2}$-excitation
with $\geq 5$ insertions. The overall count of $u^{-2}$-excitations per an FP yields
$1010$ for type $[136|136|144]$ and $659$ for $[130|130|160]$. Consequently, type
I is dominant, and, for  $u\geq u_*$, the number of extreme Gibbs measures for
$D^2=130$ equals $240$, and they are
generated by the generated by the min-area sub-lattice of type $[136|136|144]$.}}
\label{Fig11dom130}
\end{figure}

\begin{figure}[H]
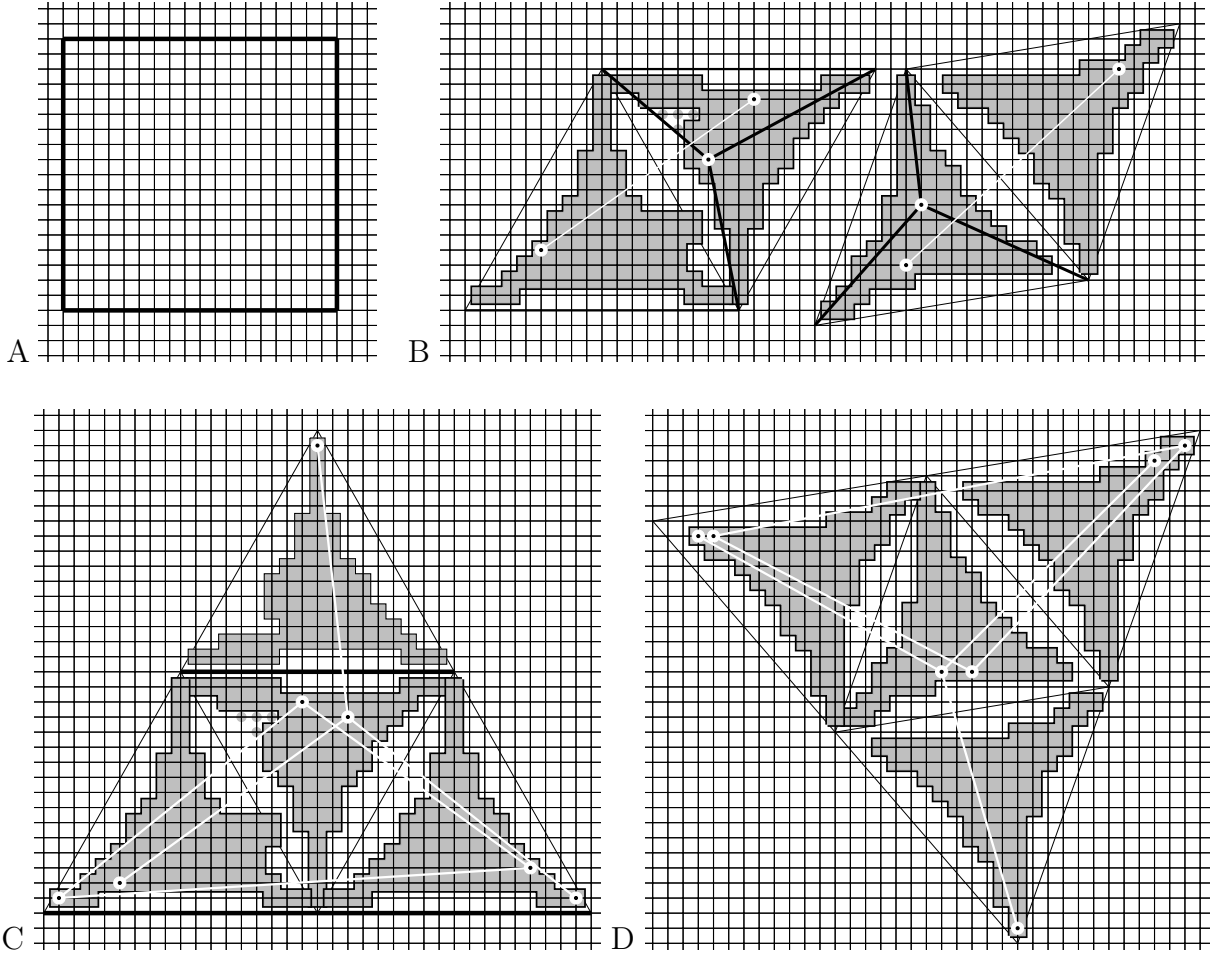

\centering
A\;

\caption{\footnotesize{Excitation count, of weight $u^{-2}$, for $D^2=324=18^2$, with
$S=288$. Frame A shows a $324\times 324$ square. Frame B indicates M-triangles
$[324|337|337]$ (left, isosceles) and
$[325|333|340]$ (right, non-isosceles). The gray areas again indicate the sites in
an M-triangle which repel  $3$ vertices: examples are highlighted with thick black lines. The
total number of 1-site insertions within an FP equals $144$ for both types. Next, the white
lines indicate double-insertions repelling four vertices of an FP; their numbers are $393$ for
$[324|337|337]$ and $408$ for $[325|333|340]$. In frame C we show 3- and 4-site insertions
for the type $[324|337|337]$ sub-lattice. There are $2510$ triple and $750$ quadruple
insertions per an FP. Frame D addresses the similar issues for type $[325|333|340]$. Here
we have $1982$ triple and $150$ quadruple insertions per an FP. As in the
previous cases,  there are no excitations of weight $u^{-2}$ with at least five insertions.
As a result, the total number of $u^{-2}$-excitations in a  {\rm{PGS}} equals 3797 for
$[324|337|337]$ and 2684 for $[325|333|340]$. Consequently, the isosceles
type is dominant.  For $u\geq u_*$, the number of the
extreme Gibbs measures for $D^2=324$ equals $576$, and they are generated by a type
$[324|337|337]$ min-area sub-lattice.}}
\label{Fig12dom324}
\end{figure}

%\begin{section}{\bf Analysis of M-triangles, 1. Class A}\label{Sec5}

\section{Minimal triangles and norm equations in $\bbZ\left[ {\sqrt[6]{-1}} \right]$}
\label{Sec5}

In Sections \ref{Sec5} and \ref{Sec6} we focus upon number-theoretical aspects of the present work.
To start with, in Section \ref{SubSec5.1} we identify $\bbZ^2$-triangles with the elements of ring 
$\bbZ\left[{\sqrt[6]{-1}}\right]$ and group them into double-infinite sequences corresponding to cosets by
the unit group in the ring $\bbZ\left[{\sqrt[6]{-1}}\right]$; this is done via the solution cosets to the  norm
equations in $\bbZ\left[{\sqrt[6]{-1}}\right]$. At the end of Section \ref{SubSec5.1} we state the relevant 
results, Theorems 5--10, which establish connections between solutions of the norm equation
and classes of values $D$.  We then prove an eventual minimality
property for $\bbZ^2$-triangles emerging from a coset; cf. Lemma 5.1 in Section 5.2.
Next, in Section \ref{Sec6} we give examples of sequences of values $D$ from various classes.
Section \ref{Sec7} contains the proof of Theorems 5--10.

%the prime factorization in ring $\bbZ\left[{\sqrt[6]{-1}}\right]$ and relate it with an analysis
%of uniqueness and non-uniqueness in problem (2.1). {\color{red} The last two sentences are outdated 
%and must be replaced:} In Sections 7 and 8 we provide
%examples of how the prime decomposition in {\color{red}$\bbZ\left[{\sqrt[6]{-1}}\right]$} can be used
%for $\bbZ^2$-implementations of isosceles triangles. This
%information is used for studying values $D$ of Classes A, B0 and B1.

We believe that the number-theoretical material regarding ring $\bbZ\left[{\sqrt[6]{-1}}\right]$ collected 
and used in Sections \ref{Sec5}--\ref{Sec7} is well-known to adepts of the algebraic number theory. In fact, a considerable part of this material can be derived
directly from the basics contained in standard textbooks \cite{BoS, Sie}, let alone such
more advanced sources as \cite{Gr, Ne}, to name a few. As a brief introduction, one can
refer the reader to papers \cite{De,IP,Le}. Nevertheless, we failed to find a single reference
where all facts that we need for our analysis are discussed in full detail.

\subsection{Triangles in $\bbZ^2$ as elements of ring $\bbZ [{\sqrt[6]{-1}}]$.}\label{SubSec5.1}

Consider a triangle with side-lengths $\ell^2_i$, $i=1,2,3$ (cf. (2.2)). It is convenient to introduce the 
notation $[\ell^2_1|\ell^2_2|\ell^2_3]$ where $\ell^2_1 \le \ell^2_2 \le\ell^2_3$ for the set of $\bbZ^2$ 
triangles with given squared lengths of sides. Clearly, all triangles in the set $[\ell^2_1|\ell^2_2|\ell^2_3]$ 
are congruent but not necessarily $\bbZ^2$-congruent.

Let us repeat the definition from Section \ref{SubSec2.1}: a pair of integers $a= \ell^2_2 - \ell^2_1$ 
and $b = \ell^2_3 - \ell^2_3$ is called a {\it signature} (of a triangle). 
All triangles with a given signature form a collection which we denote by $[\circ|\circ+a|\circ+a+b]$ 
assuming that $a$ and $b$ are fixed while the shortest squared side-length in the triangle varies.

Given non-negative $a$ and $b$ with $0 < a + b$, the collection $[\circ|\circ+a|\circ+a+b]$ may be empty. 
Later in this section we demonstrate that such collections (if they are not empty) can be partitioned into infinite sequences. Accordingly, we 
use a sequence notation, e.g. $d_j$, for running values of $\circ$. Without loss of generality we assume 
that $d_j$ is indefinitely increasing. The existence of an infinite sequence $d_j$ is the key fact behind 
the construction of infinite sequences of
values $D_j$ of Classes A and B. The second crucial fact is that all triangles in the the triangle sequence 
$[d_j|d_j+a|d_j+a+b]$ become M-triangles starting from some finite index $j^*$. The proof is given in 
the next section. The emerging theory clarifies the meaning of parameters $a$, $b$: they characterize 
closeness of the triangle to the equilateral one. 

% For a general $[d+a|d+b|d+c]$-triangle, Heron's formula gives the following 
% representation for its doubled area $s$:
% $$4s^2=3d^2 + 2d(a+b+c) +(2ab+2bc+2ca-a^2-b^2-c^2).\eqno (5.1)$$

Triangles from collection $[\circ|\circ+a|\circ+a+b]$ are constructed in terms 
of coordinates of their vertices. First, consider a single $ \bbZ^2$-triangle with 
vertices
$$(0,0),\;(m,n),\;(k,l)\eqno (5.1)$$
referred to as triangle $\{m,n,k,l\}$. (This includes degenerate triangles 
which have all 3 vertices lying in a single line.) Its squared side-lengths are
$$d= (m-k)^2 + (n-l)^2, \;\; e= k^2+l^2, \;\; f= m^2+n^2.\eqno (5.2.1)$$
The squared doubled area of triangle $\{m,n,k,l\}$ is
$$s^2 = (nk-ml)^2\eqno (5.2.2)$$
and the squared half-sum of its squared side-lengths
$$t^2 = (m^2+n^2+k^2+l^2-mk-nl)^2.\eqno (5.2.3)$$
Assuming for definiteness that $d \le e \le f$, we observe that the signature is given by
$$a^2 = (e-d)^2 = (2 m k + 2 n l- m^2 - n^2)^2, \;\;
b^2 = (f-e)^2 = (m^2 + n^2 - k^2 - l^2)^2.\eqno (5.2.4)$$
The quantity
$$r = a^2+b^2+ab\eqno (5.2.5)$$
is called the {\it norm} of triangle $\{m,n,k,l\}$ for the reason explained later. It plays a central role 
in our algebraic considerations.

It is evident that triangle $\{m,-n, k, -l\}$ is the 
image of $\{m,n, k, l\}$ under reflection about the horizontal axis in $ \bbZ^2$. Triangle 
$\{k,l,m,n\}$ corresponds to swapping two triangle vertices (which geometrically but not 
algebraically leaves the triangle intact), and triangle $\{k,-l,m,-n\}$ is the image of the swapped 
triangle under reflection about the horizontal axis in $ \bbZ^2$.

The central idea of forthcoming considerations is that triangles $\{m,n, k, l\}$ can be mapped into a well-known algebraic structure. Let $ \zeta =\diy
\frac{\sqrt{-1}+\sqrt{3}}{2}$ be the primary 12-th root of unity. From now on we use a standard 
notation $ \zeta$ instead of $\sqrt[6]{-1}$. It is convenient for us to fix the basis $\{\zeta^{10},\zeta^1,\zeta^2,\zeta^5\}$ in the corresponding cyclotomic field $\bbQ [ \zeta ]$:
$$\zeta^{10} =\frac{1-\sqrt{-3}}{2},\;\;
 \zeta^1 =\frac{\sqrt{-1}+\sqrt{3}}{2},\;\; \zeta^2 =\frac{1+\sqrt{-3}}{2},\;\; 
 \zeta^5 =\frac{\sqrt{-1}-\sqrt{3}}{2}\eqno (5.3)$$
The reason is that in a complex plane 
the angles between $\zeta^{10}$ and $\zeta^{1}$ and angles between $\zeta^2$ and $\zeta^5$ are 
$\pi /2$ while the 
angles between $\zeta^2$ and angles between $\zeta^{10}$ and $\zeta^1$ and $\zeta^5$ are $2 \pi /3$, i.e.
pairs of square lattices and pairs of triangular lattices are naturally fitting this basis. 
The triangle $\{m,n,k,l\}$ can be 
identified with an algebraic integer
$$\alpha = m  \zeta^{10}+n \zeta^1 + k  \zeta^2 + l  \zeta^5\eqno(5.4)$$
in $\bbZ [ \zeta ]$. The choice of the basis (5.3)
allows us to easily translate between manipulations 
with algebraic integers and those with triangles.

As the first example we investigate how the torsion group generated by $ \zeta^1=\{0,1,0,0\}$ acts on 
$\alpha =\{m,n,k,l\}$. Observe that
$$ \zeta^0= \zeta^2+ \zeta^{10},\;  \zeta^3= \zeta^1 +  \zeta^5,\;  \zeta^4=- \zeta^{10},\; 
 \zeta^6= - \zeta^2- \zeta^{10},\eqno (5.5.1)$$
$$ \zeta^7 = - \zeta^1,\;  \zeta^8=- \zeta^2,\;  \zeta^9 = - \zeta^1 -  \zeta^5,\;  \zeta^{11} = - \zeta^5.
\eqno (5.5.2)$$
Accordingly,
$$\{m,n,k,l\} \zeta^1 = \{-l,k, n-l,k-m\},
\eqno (5.6.1)$$
i.e. multiplication by $\zeta^1$ rotates the triangle 
by $-\pi /2$ and then shifts the image of 
the vertex $(k,l)$ to the origin. Next,
$$\{m,n,k,l\} \zeta^2 = \{m-k,n-l,m,k\},\eqno (5.6.2)$$
i.e. multiplication by $\zeta^2$ rotates the triangle 
by $\pi$ and then shifts the image of the vertex 
$(m,n)$ to the origin. Further,
$$\{m,n,k,l\} \zeta^3 = \{-n,m,-l,k\},\eqno (5.6.3)$$
i.e. multiplication by $\zeta^3$ rotates the triangle by $\pi/2$. Consequently, the multiplication by 
$ \zeta^6$ and $ \zeta^9$ corresponds to rotations by $\pi$ and $3\pi /2$. Similarly, the 
multiplication by $ \zeta^4$, $ \zeta^7$ or $ \zeta^{10}$ rotates the triangle obtained after 
multiplication by $\zeta$. In the same way, multiplication by $ \zeta^5$, $ \zeta^8$ or $ \zeta^{11}$ 
rotates the triangle obtained after multiplication by $ \zeta^2$.

The triangle symmetries 
$$\beac\{m,n,k,l\}\mapsto \{m,-n, k, -l\},\\
\{m,n,k,l\}\mapsto \{k,l,m,n\},\\
\{m,n,k,l\}\mapsto \{k,-l,m,-n\}\ena\eqno (5.7)$$
correspond to algebraic conjugates of $\a$. Combining (5.7) with multiplication by 
elements of the torsion group we obtain all 24 geometrical and all 48 algebraic congruent 
versions of a given triangle.
\medskip

The norm $N(\alpha )$ of $\a$ is the product of all 4 conjugates of an algebraic integer 
$\alpha$. It can be calculated by grouping conjugates in pairs:
$$\begin{array}{cl}
N(\alpha ) &=\big(\{m,n,k,l\}\{m,-n, k, -l\}\big)\big(\{k,l,m,n\}\{k,-l,m,-n\}\big)\\
\;&= \big( (m^2  \zeta^{-4}  + k^2  \zeta^4 +2  m k) -(n^2  \zeta^2 + l^2  \zeta^{-2} - 2 n l) \big)\\
\;&\quad\times\big( (k^2  \zeta^{-4}  + m^2  \zeta^4 + 2 m k) -(l^2  \zeta^2 + n^2  \zeta^{-2} - 2 n l) \big) \\
\;&=\big((2 m k + 2 n l - m^2 - n^2) - (m^2 + n^2 - k^2 - l^2 ) \zeta^4 \big)\\
\;&\quad\times\big((2 m k + 2 n l - m^2 - n^2) - (m^2 + n^2 - k^2 - l^2 ) \zeta^{-4} \big)\\
\;&= (m^2 + n^2 - k^2 - l^2 )^2  + (2 m k + 2 n l - m^2 - n^2)^2 \\
\;&\quad + (m^2 + n^2 - k^2 - l^2 ) (2 m k + 2 n l - m^2 - n^2). \end{array}\eqno (5.8.1)$$
%&= k^4 + 2 k^2 l^2 + l^4 - 2 k^3 m - 2 k l^2 m + 3 k^2 m^2 - l^2 m^2 -2 k m^3 + m^4 - 2 k^2 l n  \\ 
%&\hskip34pt - 2 l^3 n + 8 k l m n - 2 l m^2 n -  k^2 n^2 + 3 l^2 n^2 - 2 k m n^2 + 2 m^2 n^2 - 2 l n^3 + n^4 \\ $$
Different groupings of conjugates produce different representations for $N(\a)$:
$$\begin{array}{cl}
N(\a) &=\big(\{m,n,k,l\}\{k,l,m,n\}\big)\big(\{m,-n, k, -l\}\{k,-l,m,-n\}\big)\\
\;&= \big((m^2 - n^2 + k^2 - l^2 - m k + n l) + (2 m n + 2 k l - m l - n k) \zeta^3 \big) \\
\;&\quad\times\big((m^2 - n^2 + k^2 - l^2 - m k + n l) - (2 m n + 2 k l - m l - n k) \zeta^3 \big)\\
\;&=(m^2 - n^2 + k^2 - l^2 - m k + n l)^2 + (2 m n + 2 k l - m l - n k)^2 \end{array}
\eqno (5.8.2)$$
and
$$\begin{array}{cl}
N(\a) &=\big(\{m,n,k,l\}\{k,-l,m,-n\}\big)\big(\{k,l,m,n\}\{m,-n, k, -l\}\big)\\
\;&=\big((m^2 + n^2 + k^2 + l^2 - m k - n l) - (m l - n k)(\zeta - \zeta^5) \big) \\
\;&\quad\times\big((m^2 + n^2 + k^2 + l^2 - m k - n l) + (m l - n k)(\zeta - \zeta^5) \big) \\
\;& = (m^2 + n^2 + k^2 + l^2 - m k - n l)^2 - 3 (m l - n k)^2. \end{array}\eqno (5.8.3)$$
%&= \big((4 m k + 4 n l - m^2 - n^2 - k^2 - l^2)^2 +  3 (m^2 - k^2 + n^2 - l^2)^2\big) {\textstyle \frac{1}{ 4}}, \\ 
%&= (m^2 + n^2 - k^2 - l^2 )^2  + (2 m k + 2 n l - k^2 \hskip3pt - l^2\hskip4pt)^2  - (m^2 + n^2 - k^2 - l^2 ) (2 m k + 2 n l - k^2 \hskip3pt - l^2\hskip4pt)^2. \\ }$$

Comparing with the earlier definitions (5.8.1-3), we can see that
$$N(\a) = t^2-3s^2 = a^2 + b^2 +  a b = r.\eqno (5.9)$$
Similarly to $a$ and $b$, the value $r$ (which is always positive) can be considered as a measure 
of closeness of a $ \bbZ^2$-triangle to an equilateral one or, equivalently, a measure of  closeness 
between a triangular sublattice of $ \bbZ^2$ and a perfect triangular lattice. 

A dual problem is closeness between a sub-lattice of a triangular lattice and a square 
lattice. In this case one needs to study triangles in a triangular lattice  which are close to isosceles 
right ones. This duality is covered by expression (5.8.2) for the norm above. We interpret $(m,k)$
and $(n,l)$ (note the coordinate switch) as two sites of a unit triangular lattice $\bbA_2$ (with the 
angle $2\pi /3$ between coordinate axes) which we call a {\it dual lattice}. Accordingly, we treat 
the triple
$$(0,0),\;(m,k),\;(n,l)\eqno (5.10)$$
as a triangle in the dual lattice. For this triangle the value
$$(m^2 + n^2 + k^2 + l^2 - m k - n l)^2\eqno (5.11.1)$$
is equal to the square of the sum of squared lengths of two  of the triangle sides, while
$$3 (m l - n k)^2\eqno (5.11.2)$$
equals the squared area of the triangle multiplied by 16. Moreover,
$$ g ^2=(m^2 - n^2 + k^2 - l^2 - m k + n l)^2\eqno (5.11.3)$$
is equal to the square of the difference between squared lengths of two triangle sides while
$$ h ^2=(2 m n + 2 k l - m l - n k)^2\eqno (5.11.4)$$
is equal to the square of the difference between the squared side length of the third side and 
the sum of squared side lengths of the first two sides. The quantities $g^2$ and $h^2$ in (5.11.3,4)
again can be interpreted as a triangle signature. The sum $g^2+h^2$ is the measure of 
closeness of a triangle to an equilateral right one. 
%Indeed, it is the squared deviation from the Pythagorean theorem (the sum of squared ``catheti'' lengths minus squared ``hypotenuse'' length) plus the squared deviation from being equilateral (the difference between squared ``catheti'' lengths). With this right triangle terminology it becomes clear that for the triangular lattice the geometrical meaning of the quantity
%$$(m^2 + n^2 + k^2 + l^2 - m k - n l)^2$$
%is similar to its meanig for the square lattice. It is still something close to the half sum of all 3 squared side lengths, just remember that in the right triangle the squared hypothenuse length is a sum of the squared catheti lengths.
%Observe that the pair $( a, b)$ can be interpreted as a site of a triangular lattice such that $r$ becomes the distance of this site from the origin (assuming that the coordinate axes in the triangular lattice are at angle $\pi \over 3$). We call this lattice the {\it lattice of types}.
%We use the term a {\it dual lattice} for the triangular lattice discussed above.

The full group of units $\bbU [\zeta ]$ of the cyclotomic field $\bbQ[\zeta ]$ consists of elements
$$ \zeta^i (1+\zeta )^j,\quad i,j \in  \bbZ,\eqno (5.12.1)$$
or equivalently of elements
$$  \zeta^i (2+\sqrt{3})^j\quad {\rm and}\quad  \zeta^i (2+\sqrt{3})^j(\zeta + \zeta^2),\quad i,j \in  \bbZ.$$
Denoting by $\bbU ^+[\zeta ]$ the sub-group generated by $\zeta$ and  $(2+\sqrt{3})$ we obtain that 
$$\bbU [\zeta ] = \bbU ^+[\zeta ] \cup (\zeta + \zeta^2)\bbU ^+[\zeta ] \in  \bbZ[\zeta ].
\eqno (5.12.2)$$
Note that for brevity we slightly deviate from traditional notations which reserve superscript $+$ for the sub-group of real units. Our notation is used for a larger (not real subgroup) which contains all real units together with their multiplications by the elements of the torsion sub-group. The relationship between  $\bbU ^+[\zeta ]$ and $(\zeta + \zeta^2)\bbU ^+[\zeta ]$ can be understood from the identity
$$\begin{array}{cl}
\wt\a : &= \a (\zeta + \zeta^2) \\
&= (m  \zeta^{10}+n \zeta^1 + k  \zeta^2 + l  \zeta^5)(\zeta + \zeta^2) \\
%&= m  \zeta^{10}+n \zeta^1 + k  \zeta^2 + l  \zeta^5 + m  \zeta^{11}+n \zeta^2 + k  \zeta^3 + l  \zeta^6 \\ 
%&= m  \zeta^{10}+n \zeta^1 + k  \zeta^2 + l  \zeta^5 + m (- \zeta^5)+n \zeta^2 + k ( \zeta^1+ \zeta^5) + l (- \zeta^2- \zeta^{10}) \\ 
&= (m-k-l) \zeta^{10}+(n+k-l) \zeta^1+(m+n-l) \zeta^2+(-m+n+k) \zeta^5,\end{array}
\eqno (5.13.1)$$
i.e. multiplication by $(\zeta + \zeta^2)$ corresponds to the map
$$\{m,n,k,l\} \mapsto \{m-k-l,\; n+k-l,\; m+n-l,\; -m+n+k \}.\eqno (5.13.2)$$

After the multiplication by $(\zeta + \zeta^2)$ the norm remains the same (as $(\zeta + \zeta^2)$ is a unit) 
but the signature changes. It is a direct calculation to verify that if the original triangle signature is $a, b$
 then the resulting signature is $b,a$.

\medskip
Let 
$$l_j = \frac{1}{2}((2+\sqrt{3})^j + (2-\sqrt{3})^j)\quad{\rm and}\quad k_j 
= \frac{1}{2\sqrt{3}}((2+\sqrt{3})^j - (2-\sqrt{3})^j).\eqno (5.14)$$
The multiplication of $\a = \{m,n,k,l\}$ by 
$$(2+\sqrt{3})^j \;=\; l_j+ k_j \sqrt{3} \;=\; l_j ( \zeta^2- \zeta^4) + k_j ( \zeta^1- \zeta^5)
\eqno (5.15.1)$$ 
corresponds to the map
$$\beal\{m,n,k,l\} \\
\;\;\mapsto \{l_j m + k_j(n-2l), l_j n + k_j(2k-m),\;  l_j k + k_j(2n-l), l_j l + k_j(k-2m)\}.\ena
\eqno (5.15.2)$$
The image in (5.15.2) defines a triangle 
%$$\big(0,0\big) -\hskip-4pt-\; \big(l_j m + k_j(n-2l),\;l_j n + k_j(2k-m) \big) -\hskip-4pt-\; \big(l_j k + k_j(2n-l),\; l_j l + k_j(k-2m) \big)$$
with the same signature as the original triangle.
%$$(0,0) -\hskip-4pt-\; (m,n) -\hskip-4pt-\; (k,l).$$
Indeed, the squared side-lengths of this triangle are
$$\begin{array}{l}(m^2  + n^2) l_j^2 + (m^2+  n^2 + 4 k^2  + 4 l^2  - 4 m k - 4 n l  ) k_j^2 
- 4 (m l  - n k) k_j l_j,\\
(k^2  + l^2) l_j^2 + (4 m^2 +  4 n^2 + k^2  + l^2  - 4 m k  - 4 n l  ) k_j^2 - 4 (m l  - n k) k_j l_j, \\
(m^2+ n^2 + k^2  + l^2  -  2 m k    - 2 n l  ) l_j^2\\
\qquad\qquad\qquad + (m^2+  n^2 + k^2  + l^2  + 2 m k  + 2 n l  ) 
k_j^2  - 4 (m l  - n k) k_j l_j.\end{array}\eqno (5.16)$$
Using the relation $l_j^2 = 3 k_j^2 +1$ (a Pell equation), it is not hard to verify that the 
difference between the first 
and second squared lengths in (5.16) is $ b$ whereas the difference between the second and 
third squared lengths is $ a$.

Thus, any triangle $\{m,n,k,l\}$
%$$(0,0) -\hskip-4pt-\; (m,n) -\hskip-4pt-\; (k,l)$$
with a given signature $a,b$ is always a member of a coset  of triangles 
$\{m,n,k,l\}\bbU ^+[\zeta ]$ (explicitly listed above), with the same signature and norm.
%The full coordinate-wise (in $ \bbZ^2$) description of this coset is given above. 
Considering a coset (of $\bbU [\zeta ]$ or $\bbU ^+[\zeta ]$) as a sequence  labeled by 
$j \in  \bbZ$ (up to multiplication by an element of the torsion group) we can see that the larger 
$|j|$ the closer the corresponding triangle is to the equilateral one. In particular, each coset 
contains only a finite number of obtuse triangles.

If we treat quadruples $\{m,n,k,l\}$ as a row vector in $\bbR^4$ then it is a straightforward calculation to verify that a multiplication by an algebraic number $\{m,n,k,l\}$ is equivalent to  multiplication by the matrix
$$ T=\begin{pmatrix}k & l & k-m & l-n \\-l & k & n-l & k-m \\
m-k & n-l & m & n \\ l-n & m-k & -n & m \end{pmatrix}.\eqno (5.17)$$
In particular, relationships (5.13.2) and (5.15.2) are the results of the multiplication by the corresponding matrices.

A natural approach to identifying all triangles with a given signature $a$, $b$ consists in finding 
all solutions to the corresponding norm equation
$$N(\a)=r,\eqno (5.18)$$
where $r=a^2+b^2+ a b$. If solutions exist they can be further classified into Classes A, B0, B1 
and B2 (see below) closely related to eponymous classes defined earlier for values $D_n$.

Being cyclotomic, the field $\bbQ[\zeta ]$ possesses a unique factorization property for the
corresponding algebraic integers in $ \bbZ[\zeta ]$. This implies the following algorithm of solving 
the norm equation. (I) First, analyze how rational primes are decomposed into algebraic primes. 
(II) After that decompose $r$ into the product of rational primes and then decompose each rational 
prime into the product of algebraic primes. (III) Finally, partition the entire product of algebraic 
primes into 4 sub-products such that these sub-products give 4 conjugated   algebraic integers 
from $ \bbZ[\zeta ]$. Each such partition corresponds to some solution coset $\{m,n,k,l\}\bbU [\zeta ]$ 
to the norm equation, and all solutions to the norm equation can be found in this way. Given $r$, 
the above algorithm  allows us to identify all the corresponding solution cosets and therefore calculate
 their total amount (which is always finite). It also clarifies that a larger amount of different prime 
 factors in the prime decomposition of $r$ leads to a larger number of solution cosets to the norm
  equation (5.18). 

It is convenient to pick a unique representative of a solution coset $\{m,n,k,l\}\bbU [\zeta ]$  
which we call a {\it leader} or a {\it leading} triangle.  We choose it to be a triangle with the minimal 
area within the coset. Moreover, among all 48 possible representations of such a triangle (12 
members of the torsion group times 4 conjugations) we select the triangle, denoted by  $\{m|n|k|l\}$, 
with $m,n \ge 0$ and $m^2+n^2 \ge k^2+l^2 \ge (m-k)^2 + (n-l)^2$. It is not hard to verify that the 
leading triangle is always unique. 
%In other words, $\{m,n,k,l\}$ denotes any quadruple of integers while $\{m|n|k|l\}$ additionally claims that this quadruple represents the leading triangle in the corresponding coset by $\bbU [\zeta ]$. 
When we speak about a full coset, $\{m|n|k|l\}\bbU [\zeta ]$ we usually assume that the triangle 
$\{m|n|k|l\}$ corresponds to $j=0$. The notion of a leader allows us to use a ``canonical'' notation 
for different types of solution cosets: $\{m|n|k|l\}\bbU [\zeta ]$, $\{m|n|k|l\}\bbU ^+[\zeta ]$, $\{m|n|k|l\}\bbU ^+[\zeta ](\zeta + \zeta^2)$, etc.

For the future reference we summarize that the coset leader $\{m|n|k|l\}$ uniquely defines
$$\begin{array}{c}
d = (m-k)^2 + (n-l)^2, \;\; e= k^2+l^2 \;\; f= m^2+n^2, \\
 a^2 = (2 m k + 2 n l- m^2 - n^2)^2, \;\; b^2 = (m^2 + n^2 - k^2 - l^2)^2,\\
 g^2 =(m^2 - n^2 + k^2 - l^2 - m k + n l)^2, \;\; h^2 =(2 m n + 2 k l - m l - n k)^2, \\
s^2 = (nk-ml)^2,\;\;t^2 = (m^2+n^2+k^2+l^2-mk-nl)^2 \end{array}\eqno (5.19)$$
with
$$r=a^2+b^2+ab=g^2+h^2=t^2-3s^2.\eqno (5.20)$$
A full coset $\{m|n|k|l\}\bbU [\zeta ]$ is the union 
$$\{m|n|k|l\}\bbU ^+[\zeta ]\cup \{m|n|k|l\}\bbU ^+[\zeta ](\zeta + \zeta^2).\eqno (5.21)$$
Here 
$$\begin{array}{l} \{m|n|k|l\}\bbU ^+[\zeta ] \\
= \big\{l_j m + k_j(n-2l),l_j n + k_j(2k-m),l_j k + k_j(2n-l),l_j l + k_j(k-2m)\big\},\\
=\{l_j, k_j, l_j, -k_j \}T,\quad j\in\bbZ,\end{array}\eqno (5.22)$$
is a coset with signature $a,b$ while $\{m|n|k|l\}\bbU ^+[\zeta ](\zeta + \zeta^2)$ 
is a coset with signature $b,a$. The latter is 
obtained from the former via the map 
$$\begin{array}{l}\{m,n,k,l\} \\
\;\;\mapsto \{m-k-l,\; n+k-l,\; m+n-l,\; -m+n+k \}=\{0,1,1,0\} T.
\end{array}\eqno (5.23)$$
It will be convenient to refer to $\{m|n|k|l\}\bbU ^+[\zeta ]$ and $\{m|n|k|l\}\bbU ^+[\zeta ](\zeta + \zeta^2)$
in (5.21) as the $\{m|n|k|l\}\bbU ^+[\zeta ]$- and $\{m|n|k|l\}\bbU ^+[\zeta ](\zeta + \zeta^2)$-halves
of the full coset $\{m|n|k|l\}\bbU [\zeta ]$, respectively.
\medskip

As was mentioned earlier, an important property of a solution coset is the existence of a finite 
index $j^*=j^*(r)$ such that for any $j$ with $|j| \ge j^*=j^*(r)$ the area $s_j$ of the corresponding 
triangle is minimal among all acute $\bbZ^2$-triangles having their squared side-lengths at 
least $d_j$. This assigns a quantitative meaning to the claim that the triangles in the coset 
become close to an equilateral one as $|j|\to \infty$. The proof of this property 
is presented in the next section (see Lemma 5.1). 

We say that a full solution coset $\{m|n|k|l\}\bbU [\zeta ]$ is of {\it Class A} if there is no other 
coset $\{m',n',k',l'\}\bbU [\zeta ]$ with $s'_j = s_j$ for all $j \in  \bbZ$. The meaning is that for 
each $d_j$ with $|j| >j^*$, the minimal area of any acute $\bbZ^2$-triangle
with squared side-lengths 
$\ge d_j$ is equal to $s_j$, and it is achieved specifically on the corresponding triangle 
$\{m_j,n_j,k_j,l_j\}$ from the coset. (Here $k_j$ and $l_j$ should not be confused with the ones defined in (5.15.1).) Referring to the representation (5.21), we will say that 
$\{m|n|k|l\}\bbU ^+[\zeta ]$ and $\{m|n|k|l\}\bbU ^+[\zeta ](\zeta + \zeta^2)$ are
solution cosets of Class A with constant signatures $a,b$ and $b,a$, respectively. 

Next, a family ${\mathfrak B}_0$ of solution cosets $\{m^{(i)}|n^{(i)}|k^{(i)}|l^{(i)}\}\bbU [\zeta ]$ 
is said to be of {\it Class B0} if for each $j$ the triangles 
labeled by $j$ taken from any two cosets in the family are congruent to each 
other, and ${\mathfrak B}_0$ is the maximal one with this property. We will see 
later that it may happen only if the cosets in the family are obtained from each 
other via some $\bbR^2$-rotation by an angle not multiple of $\pi /2$. For brevity, we often 
say that ${\mathfrak B}_0$ is a B0 family; a similar agreement is in place for families
of Classes B1 and B2; see below.

We say that a family ${\mathfrak B}_1$ of solution cosets 
$\{m^{(i)}|n^{(i)}|k^{(i)}|l^{(i)}\}\bbU [\zeta ]$ 
is a {\it Class B1} family if for each $j$ the triangles labeled by $j$ taken from any two 
cosets in the family are not congruent but have the same area, 
and the family is the maximal one with this property. This is 
related to $\bbR^2$-rotations of the corresponding triangles in the dual lattice $\bbA_2$.

Finally, a family ${\mathfrak B}_2$ of solution cosets $\{m^{(i)}|n^{(i)}|k^{(i)}|l^{(i)}\}\bbU [\zeta ]$ 
is said to be of  {\it Class B2}  if, for any $j$, all triangles labeled by $j$ in these cosets have the 
same area but some of them are congruent while some other are not,  and the family is the 
maximal one with this property. This is related to superposition of rotations in both main and 
dual lattices.

The results of this section are given in Theorems 5 -- 10 below. Here we suppose that the 
values $\vrho (\cdot )$ are non-negative integers.

\bthme\label{Theorem5} A solution to the norm equation $N(\a) = r$
exists only for positive integers $r$ with the prime decomposition
$$r = 2^{2\vrho (2)} 3^{2\vrho (3)} \prod_i p_i^{\vrho (p_i)} \prod_j q_j^{2\vrho (q_j)} \prod_k 
w_k^{2\vrho (w_k)} \prod_l z_l^{2\vrho (z_l)}.\eqno (5.24)$$
Here all $p_i \bmod 12 = 1$, $q_j \bmod 12 = 11$,  $w_k \bmod 12 = 5$ and $z_l \bmod 12 = 7$.
\ethme 

The further relation between values of $r$ of the form (5.24) and solution cosets 
to the norm equation is classified in Theorems 6--10.

\bthmf\label{Theorem6} {\bf (Class A).} Assume  that 
$$r = 2^{2\vrho (2)} 3^{2\vrho (3)} p^{\vrho (p)} w^{2\vrho (w)} z^{2\vrho (z)} 
\prod_j q_j^{2\vrho (q_j)}, \eqno (5.25)$$
where $\vrho (p)+\vrho (w)+\vrho (z) \in \{0,1\}$. 
Then the corresponding norm equation $N(\a)=r$ has at least one pair of solution cosets
with constant signatures $a,b$ and $b,a$ such that $r=a^2+b^2+ab$:
$$\{m|n|k|l\}\bbU ^+[\zeta ]\quad{\rm and}\quad \{m|n|k|l\}\bbU ^+[\zeta ](\zeta + \zeta^2), 
\eqno (5.26)$$
and the union $\{m|n|k|l\}\bbU ^+[\zeta ]\cup\{m|n|k|l\}\bbU ^+[\zeta ](\zeta + \zeta^2)$ gives
a full solution coset of Class {\rA}. 
Conversely, any Class {\rA} solution coset is obtained as a union of cosets {\rm{(5.26)}} for
$a,b$ and $b,a$ with $a^2+b^2+ab=r$ where $r$ has the form {\rm{(5.25)}}.

Let $\sharp(r)$ denote the total number of pairs of solution cosets {\rm{(5.26)}} of
constant signatures $a,b$ and $b,a$ (that is, of full solution cosets). Then 
$$\sharp (r) = 1\;\hbox{ iff }\;\vrho (p)+\vrho (w)+\vrho (z)+\sum_j \vrho (q_j) \le 1.$$ 
If $\vrho (q_j)=1$ for all $j$ then
$$\sharp (r) = 2^{-1+\vrho (p)+\vrho (w)+\vrho (z)+\sum_j \vrho (q_j)}.$$
If all $\vrho (\cdot)=0$ except for one $q$ then 
$$\sharp(r) = [\vrho (q)/2]+1.$$ 
If 
$\vrho (p)+\vrho (w)+\vrho (z)=1$ and the remaining $\vrho (q_j)=0$ except for  
$q_1$ then 
$$\sharp (r) = \vrho (q_1)+1.$$ 
In a general case $\sharp (r)$ can be 
expressed via the collection of divisors of $r$ in a fashion similar to Theorem {\rm 3}, 
Chapter {\rm 1},~\S{\rm 4} in {\rm\cite{G}}.
\ethmf

\bthmg\label{Theorem7} {\bf (Class B0).} Assume that
$$r = 2^{2\vrho (2)} 3^{2\vrho (3)} p^{\vrho (p)} \prod_j q_j^{2\vrho (q_j)}   
\prod_k w_k^{2\vrho (w_k)} 
z^{2\vrho (z)}. \eqno (5.27)$$
Here $\vrho (p), \vrho (z)\in \{0,1\}$,  $\vrho (p) + \vrho (z) \le 1$, and 
$\vrho (p)+\vrho (z) + \sum_k \vrho (w_k) \ge 2$. Then the corresponding solution 
cosets are partitioned into one or several Class {\rB}$0$ families. Conversely, every 
{\rB}$0$ family corresponds to an $r$ as in {\rm{(5.27)}}. Given an $r$ of the 
form {\rm{(5.27)}}, the number $\sharp (r)$ of cosets in each emerging {\rB}$0$ family  
is the same. If additionally $\vrho (w_k)\le 1$ for all $k$ then  
$$\sharp (r)=2^{-1+\vrho (p) + \vrho (z)+\sum_k \vrho (w_k)}.$$ 
Increasing the value 
$\vrho (w_k)>1$ increases $\sharp (r)$. In particular, if $\vrho (w_k)=0$ for $k\not=1$ then 
$$\sharp (r)=\begin{cases} \vrho (w_1)+1,&\hbox{if $\vrho (p) + \vrho (z)\not=0$,}\\
 [\vrho (w_1)/2]+1, &\hbox{otherwise.}\end{cases}$$
\ethmg

\bthmh\label{Theorem8} {\bf (Class B1).} Assume that
$$r = 2^{2\vrho (2)} 3^{2\vrho (3)} \prod_i p_i^{\vrho (p_i)} \prod_j q_j^{2
\vrho (q_j)} w^{2\vrho (w)} \prod_l z_l^{2\vrho (z_l)} ,\eqno (5.28)$$
where $\vrho (w), \vrho (p_i)\in \{0,1\}$, and $\sum_i \vrho (p_i) + \sum_l \vrho (z_l) \ge 2$. 
Then the corresponding cosets are partitioned into one or 
several Class {\rB}$1$ families. The number $\sharp (r)$ of cosets in each 
emerging {\rB}$1$ family is the same.  If additionally 
$\vrho (z_l)\le 1$ then 
$$\sharp(r)=2^{-1+\vrho (w)+\sum_i \vrho (p_i)+\sum_l \vrho (z_l)}.$$ 
Increasing the value $\vrho (z_l) > 1$ increases $\sharp(r)$. In particular, if all 
$\vrho (p_i)=0$ and  $\vrho (z_l)=0$ for $l\not=1$ then 
$$\sharp(r)=\begin{cases} \vrho (z_1)+1,&\hbox{if $\vrho (w)\not=0$,}\\
[\vrho (z_1)/2]+1, &\hbox{otherwise.}\end{cases}$$
\ethmh

\bthmi\label{Theorem9} {\bf (Class B2).} Assume that 
$$r = 2^{2\vrho (2)} 3^{2\vrho (3)} \prod_i p_i^{\vrho (p_i)} \prod_j q_j^{2\vrho (q_j)} 
\prod_k w_k^{2\vrho (w_k)}  \prod_l z_l^{2\vrho (z_l)},$$
where $\vrho (p_i)\in \{0,1\}$, and $ \sum_k \vrho (w_k) \ge 1$, $\sum_l \vrho (z_l) \ge 1$, 
$\sum_i \vrho (p_i) + \sum_k \vrho (w_k) \ge 2$, $\sum_i \vrho (p_i) +  \sum_l \vrho (z_l) \ge 2$. 
Then the corresponding cosets are partitioned into one or several {\rB}$2$ families. The 
number $\flat(r)$ of cosets in each emerging {\rB}$2$  family is the same. If 
$\vrho (w_k)=1$ for all $k$ and $\vrho (z_l)=1$ for all $l$ then 
$$\flat(r)=2^{-2+\sum_i \vrho (p_i)+\sum_k \vrho (w_k)+\sum_l \vrho (z_l)}.$$
Increasing $\vrho (w_k) > 1$ and/or $\vrho (z_l) > 1$ increases $\flat (r)$.
\ethmi

\bthmj\label{Theorem10} {\bf (Classes B1,2).} Assume that 
$$r = p^{\vrho (p)},$$
where $\vrho (p) > 2$. Then the whole set of solution cosets has cardinality
$$\sharp (r)=\big((\vrho (p) + 2)(2 \vrho (p)^2+ 8 \vrho (p) + 15 + 9 (-1)^{\vrho (p)} \big)/ 48$$
and is partitioned into $\flat(r)=[\vrho (p)/2]+1$ families, where one family is of Class {\rB}$1$ 
and the remaining families are of Class {\rB}$2$.
\ethmj

The proof of Theorems 5 -- 10 is given in Section 7: it is
essentially an application of well-known facts 
from the algebraic number theory (see \cite{BoS}, \cite{W}, \cite{L}) to the triangle 
representation developed above. A feature of 
the proof is that the statements of Theorems 5 -- 10 are established via 
a common argument that goes in parallel for several or all of them. 

\subsection{Coset intersections and minimal triangles.}\label{SubSec5.2} 
In Section \ref{SubSec5.2} we prove eventual minimality of triangles in a solution coset. This finishes the proof of Theorem~3.
 
%{\bf Lemma 5.1.}
\bl Consider a coset $\{m|n|k|l\}\bbU[ \zeta]$ with norm $r$, where
$$r =(m^2 - n^2 + k^2 - l^2 - m k + n l)^2 + (2 m n + 2 k l - m l - n k)^2.$$
Then there exists an integer $j_*$, with $0 \leq j_* \leq \log_{2+\sqrt{3}}(3200 r^{3/2})$,  such 
that for $|j| > j_*$ the members of the coset indexed by $j$ are M-triangles (for all 
corresponding $D^2$). 
\el

\bp Consider a non-obtuse triangle
$[d+a|d+b|d+c]$ where $0 \leq a\leq b\leq c$ and $d \geq c-a-b$. Let $s=s(a,b,c)$ be the 
doubled area of this triangle. Then 
$$4s^2=3d^2 + 2d(a+b+c) +(2ab+2bc+2ca-a^2-b^2-c^2).$$
In particular, 
$$\frac{\partial (4s^2)}{\partial c} = 2d +2(a+b)-2c = 2(d -(c-a-b))>0,$$
and the area of the triangle is
increasing in $c$. It is increasing in $a$ and $b$ for the same reason.

Consider another triangle $[d+a'|d+b'|d+c']$. Our aim is to compare the areas of these two 
triangles; without loss of generality we may assume that $a=0$. Then
$$s(d, d+b, d+c) \le s(d,d+c,d+c) < s(d, d, d+4c) \le s(d+a',d+b',d+4c')$$
whenever $b \ge 0$, $c \ge \max(1,b)$, $a' \ge 0$, $b' \ge b$, $c' \ge c$ and
$d>4c$. The only non-trivial inequality here is in the middle:
$$4s(d,d+c,d+c)^2 = 3d^2+4dc< 3d^2 +8dc-16c^2=4s(d,d,d+4c)^2$$
which follows from the acuteness requirement $d > 4c$ for a
$[d|d|d+4c]$-triangle.
Thus, given $b \ge 0$ and $c \ge \max(1,b)$, the area of a non-obtuse
$[d'+a'|d'+b'|d'+c']$-triangle
may be smaller than that of a non-obtuse $[d|d+b|d+c]$-triangle for at most
$(4c)^3$ triples $a',b',c' \ge 0$.
In fact, since no dependency of $d$ upon
$a'$, $b'$ and $c'$ is assumed, some of acute $[d'+a'|d'+b'|d'+c']$-triangles
may not exist, which makes the upper-bound $(4c)^3$ rather excessive.

For the rest of the argument we continue to assume that $a=0$. Given $b \ge 0$ and
$c \ge \max(1,b)$, we consider the corresponding norm equation $N(\{m,n,k,l\})=r$ with
$$r= a^2+b^2+c^2-ab-bc-ca.$$
As we are primarily interested in areas of triangles it is convenient to understand 
the norm equation in terms of variables $s^2 = (nk-ml)^2$, $t^2 = (m^2+n^2+k^2+l^2-mk-nl)^2$.
Cf. (5.2.1,2). In these variables the norm equation is $t^2-3s^2=r$, cf. (5.9).
For each coset $\{m|n|k|l\}\bbU[ \zeta]$ solving this equation we denote by $(x,y)$ 
the corresponding leading solution (here $x=x(0,b,c)$, $y=y(0,b,c)$ and $y^2-3x^2=r$).
In these variables the triangle $\{m|n|k|l\}$ can be written as $[d|d+b|d+c]$ where
$d=(2y-(b+c))/3$. The double area of this triangle is equal to $x$.

Denote by $\triangle$ an element of the triangle collection $[\circ|\circ+b|\circ+c]$. 
The double area of $\triangle$ is denoted by $s(\triangle )$
and equals to $x$. For a different collection $[\circ+a'|\circ+b'|\circ+c']$ with elements 
denoted by $\triangle'$ we would like to understand for which $a', b', c'$ one may have 
$s(\triangle' ) < s(\triangle)$. We already know that there exist only
finitely many triples $(a',b',c')$ for which the double area $s(\triangle') < s(\triangle )$
for at least one value $d$ (i.e. value of $\circ$). Next, for each such triangle $[d+a'|d+b'|d+c']$
there exists only a finite number of cosets solving the corresponding norm
equation  ${t'}^2-3{s'}^2=r'$ where $r'=r'(a', b', c')$. Indeed, given a leading  solution 
$(x,y)$ of ${t}^2-3{s}^2=r$ and a leading
solution $(x',y')$ of ${t'}^2-3{s'}^2=r'$, we have
$$0< x < \sqrt{r},\quad 0< y < 2\sqrt{r},\;\hbox{ and }\; 0< x' < \sqrt{r'},\quad 0< y' <
2\sqrt{r'}$$
as the leader is the triangle with the smallest area.
On the other hand,
$$0 \le r =a^2+b^2+c^2-ab-bc-ca \le c^2,$$
and we are interested only in $0 \le c' \le 4c$ \ and, consequently, in
$0 \le r' \le 16 c^2$.

If $a'+b'+c' = a+b+c$ then \ ${\rm Sign}(s(\triangle') - s(\triangle)) = {\rm Sign}(r-r')$ \ since
$$4s(\triangle')^2=\left(\sqrt{3}d + \frac{1}{\sqrt{3}}(a'+b'+c')\right)^2 - \frac{4}{3}r',$$
and
$$4s(\triangle)^2=\left(\sqrt{3}d + \frac{1}{\sqrt{3}}(a+b+c)\right)^2 - \frac{4}{3}r.$$
In the opposite case $a'+b'+c' \not = a+b+c$ we have equality $s(\triangle') = s(\triangle)$
at $d=d^*$ where
$$2d^*= -\frac{1}{3}\left( a'+b'+c'+a+b+c\right) + \frac{4}{3}\frac{r'-r}{a'+b'+c'-a-b-c}.$$
Therefore, for
$$d > 16c^2\ge\left|d^* \right|$$ % \left|\frac{r'-r}{a'+b'+c'-a-b-c}\right|$$
the function ${\rm Sign}(s(\triangle') - s(\triangle))$ is constant.
If $a'+b'+c' > a+b+c$ then for $d$  large enough
$$\left(\sqrt{3}d + \frac{1}{\sqrt{3}}(a'+b'+c')\right)^2 - \frac{4}{3}r' > \left(\sqrt{3}d + \frac{1}{\sqrt{3}}(a+b+c)\right)^2 - \frac{4}{3}r \eqno (5.29)$$
and consequently $s(\triangle') > s(\triangle)$. Due to the previous argument,
inequality (5.29) is true already for $d > 16c^2$. As we are interested only
in the case when $s(\triangle') \le s(\triangle)$, we assume from now on that
$$a'+b'+c' \le a+b+c.$$

Specifically, we are looking for the situation where  $3d'_i=2y'_i-(a'+b'+c')=2y_n -(a+b+c)=3d_n$
and $x'_i=x_n$. In particular, in this case we have
$$0 \le 2y_n - 2y'_i = (a+b+c) - (a'+b'+c') < 3c.$$
Thus, to finish the proof of lemma it is enough to verify that
$$2y_n - 2y'_i > 3c,\;\hbox{ whenever \ $i,n >
\log_{2+\sqrt{3}}(3200 c^3)$.} \eqno (5.30) $$

We start from the representation
$$ \begin{array}{l}
2y_n - 2y'_i = 2yl_n+6xk_n-2y'l_i-6x'k_i \\
\quad= \big((y-y'l_{i-n}-3x' k_{i-n}) + \sqrt{3}(x-x'l_{i-n}-y'
k_{i-n})\big)(2+\sqrt{3})^n \\
\quad\quad + \big((y-y'l_{i-n}+3x' k_{i-n}) + \sqrt{3}(x-x'l_{i-n}+y'
k_{i-n})\big)(2-\sqrt{3})^n,\end{array}$$
where $k_j$ and $l_j$ are given by (5.15.1).
%$$l_j = {1/2}((2+\sqrt{3})^j + (2-\sqrt{3})^j)\quad{\rm and}\quad k_j 
%=\frac{1}{2\sqrt{3}}((2+\sqrt{3})^j + (2-\sqrt{3})^j)$$
Next we observe that $l_{n+j} >\diy (2+\sqrt{3})^j l_n/2$ and $k_{n+j} >
(2+\sqrt{3})^j k_n$ for $n,j > 0$. Therefore, if $(2+\sqrt{3})^{i - n} > 4c $
then
$$ 2yl_n+6xk_n-2y'l_i-6x'k_i < 4 c l_n+6c k_n -2l_i-6k_i < 0.$$
Similarly, if $(2+\sqrt{3})^{n - i} > 19c $ then
$$  2yl_n+6xk_n-2y'l_i-6x'k_i > 2l_n+6k_n-16cl_i-24ck_i > 3c.$$
Thus, we may safely assume that $(2+\sqrt{3})^{|n - i|} < 19c$
and consequently $|l_{i-n}|, |k_{i-n}| < 10c$.

It is not hard to verify that for positive integer $p$ and $q$, we have
$|p\pm q\sqrt{3} | > 1/(4 q)$.
Then for $n > \log_{2+\sqrt{3}}(400 c^2)$
$$\begin{array}{l}\Big|(y-y'l_{i-n}+3x' k_{i-n})\\
\qquad + \sqrt{3}(x-x'l_{i-n}+y'
k_{i-n})\Big|(2-\sqrt{3})^n
\le 400 c^2 (2-\sqrt{3})^n < 1,\end{array}$$
while for $n > \log_{2+\sqrt{3}}(3200 c^3)$
$$\begin{array}{r}\Big|(y-y'l_{i-n}-3x' k_{i-n}) + \sqrt{3}(x-x'l_{i-n}-y'
k_{i-n})\Big|(2+\sqrt{3})^n\qquad{}\\
\diy>\frac{(2+\sqrt{3})^n}{4(x-x'l_{i-n}-y'k_{i-n})}
>\frac{(2+\sqrt{3})^n}{800 c^2} > 4c > 3c+1.\end{array}$$
The straightforward estimate $c < \sqrt{r}$ completes the proof of Lemma 5.1.
\ep

%Let us discuss
%basic facts about the structure of
%$\bbZ^2$-triangles, with squared side-lengths $\ell^2_i$, $i=1,2,3$ (cf. (2.2)); here
%each $\ell^2_i$ is a Gaussian number.

\section{Examples of Classes A, B0, B1 and B2}\label{examples}
\label{Sec6}

\bigskip

In this section we provide series of simplest examples related to Classes A and 
B0 - B2. The notations $k_j$ and $l_j$ refer to quantities defined in (5.15.1).
\medskip

{\bf 6.1.} The coset $\{1|0|1|0\}\bbU^+[\zeta](\zeta+\zeta^2)=\{0,1,1,0\}\bbU^+[\zeta]$ 
corresponding to $r=1$ is of Class A; it contains all triangles from collection $[\circ|\circ|\circ+1]$,
i.e, with signature $0,1$. These triangles have vertices  
$$(0,0),\;(2k_j-l_j, k_j),\;(k_j, 2k_j-l_j), \quad j\in\bbN;$$
they provide the best $\bbZ^2$-approximations to an equilateral triangle. These are M-triangles for 
the corresponding $D^2_j = l_j^2 + 5k_j^2 -4l_j k_j$ for all $j$ and possibly for some smaller 
values ${\wt D}$ with $D^*_j\leq{\wt D}<D_j$. (The latter possibility is admitted in all examples 
but not stressed 
every time again. The range $j\in\bbN$ is also purported throughout the whole section.)
\medskip

In Examples {\bf 6.2}, {\bf 6.3} we construct families containing two cosets each. 
The families in Examples {\bf 6.4}, {\bf 6.5}  have larger cardinalities. 
\medskip

{\bf 6.2.} The cosets  $\{5|0|5|0\}\bbU^+[\zeta]$ and $\{4|3|4|3\}\bbU^+[\zeta]$ corresponding 
to $r=5^4$ form a B0 family containing triangles from collection $[\circ|\circ+25|\circ+25]$.
These triangles have vertices
$$(0,0),\;(5l_j, 5k_j),\;(5l_j, -5k_j)$$
and
$$(0,0),\;(4l_j-3k_j, 3l_j+4k_j),\;(4l_j+3k_j, 3l_j-4k_j),$$
respectively. These are M-triangles for $j\geq 3$ and $D^2_j =100k_j^2$.
\medskip

{\bf 6.3.} The cosets  $\{2|5|3|4\}\bbU^+[\zeta]$ and $\{5|3|4|1\}\bbU^+[\zeta]$ corresponding to 
$r=13\cdot 7^2$ form a B1 family containing triangles from collections $[\circ|\circ+23|\circ+4]$ 
and $[\circ|\circ+12|\circ+17]$. The vertices are
$$(0,0),\;(2l_j-3 k_j, 5l_j+ 4k_j),\;(3l_j+6k_j, 4l_j-k_j)$$
and
$$(0,0),\;(5l_j+k_j, 3l_j+3k_j),\;(4l_j+5k_j, l_j-6k_j),$$
respectively. Here we obtain M-triangles for $j \geq 2$ and $D^2_j =2l_j^2+106k_j^2 +28 k_j l_j$. 
\medskip

{\bf 6.4.} The cosets  $\{9|8|10|5\}\bbU^+[\zeta]$, $\{12|1|11|-2\}\bbU^+[\zeta]$, 
$\{7|11|7|6\}\bbU^+[\zeta]$ and \\
$\{13|1|9|-2\}\bbU^+[\zeta]$ corresponding to $r=13\cdot5^2\cdot 7^2$ form a B2 family
containing triangles from collections $[\circ|\circ+60|\circ+145]$, $[\circ|\circ+60|\circ+145]$, 
$[\circ|\circ+115|\circ+135]$ and $[\circ|\circ+115|\circ+135]$ respectively.
These triangles have vertices 
$$(0,0),\;(9l_j-2k_j, 8l_j+11k_j),\;(10l_j+11k_j, 5l_j-8k_j), $$
$$(0,0),\;(12l_j+5k_j, l_j+10k_j),\;(11l_j+4k_j, -2l_j-13k_j), $$
$$(0,0),\;(7l_j-k_j, 11l_j+7k_j),\;(7l_j+16k_j, 6l_j-7k_j), $$
and
$$(0,0),\;(13l_j+5k_j, l_j+ 5k_j),\;(9l_j+4k_j, -2l_j -17k_j) $$
respectively. They are M-triangles for $j\geq 3$ and $D^2_j =10l_j^2+530k_j^2 +140 k_j l_j$.
\medskip

{\bf 6.5.} The cosets  $\{6|5|4|3\}\bbU^+[\zeta]$ and $\{2|7|2|6\}\bbU^+[\zeta]$ corresponding 
to $r=13^3$ form a B1 family ${\mathfrak B}_1$ containing triangles 
from collections $[\circ|\circ+17|\circ+53]$ and $[\circ|\circ+39|\circ+52]$ respectively, 
while the cosets $\{8|1|5|-1\}\bbU^+[\zeta]$, $\{7|4|5|1\}\bbU^+[\zeta]$ and 
$\{3|7|4|5\}\bbU^+[\zeta]$ corresponding to the same $r=13^3$ form a B2 
family ${\mathfrak B}_2$ containing triangles from collections $[\circ|\circ+39|\circ+52]$, 
$[\circ|\circ+39|\circ+52]$ and $[\circ|\circ+17|\circ+53]$ respectively. The triangles in  
family ${\mathfrak B}_1$ have vertex triples  
$$(0,0),\;(6l_j-k_j, 5l_j+ 2k_j),\;(4l_j+7k_j, 3l_j -8k_j), $$
and
$$(0,0),\;(2l_j-5k_j, 7l_j+ 2k_j),\;(2l_j+8k_j, 6l_j -2k_j),$$
respectively. They are M-triangles for $j \geq 3$ and $D^2_j =l_j^2+185k_j^2 +8 k_j l_j$. 
The triangles in family ${\mathfrak B}_2$ have vertices
$$(0,0),\; (8l_j+3k_j, l_j+ 2k_j),\;(5l_j+3k_j, -l_j -11k_j), $$
$$(0,0) ,\;(7l_j+2k_j, 4l_j+ 3k_j),\;(5l_j+7k_j, l_j -9k_j) $$
and
$$(0,0),\;(3l_j-3k_j, 7l_j+ 5k_j),\;(4l_j+9k_j, 5l_j -2k_j), $$
respectively. They are M-triangles for $j \geq 3$ and $D^2_j =5l_j^2+193k_j^2 +52 k_j l_j$. 
\medskip

In all examples {\bf 6.1} - {\bf 6.5} the minimal value of $j$ has been calculated numerically by {\tt  MinimalTriangles.java}.

\section{Proof of Theorems 5-10}\label{Sec7}

We start with a review of decomposition of rational primes into the product of algebraic primes 
in $\bbZ [\zeta ]$ and the corresponding solutions to the norm equation. Then we 
proceed to the powers of rational primes and finally to the products of different rational primes, 
i.e. to the case of  general $r$. In this section the notations $k_j$ and $l_j$ refer to quantities defined in (5.15.1).
\bigskip

{\bf 7.1.} For $r=1$ the solutions to the corresponding norm equation $N(\a)=1$ (Pell's equation)
constitute the group of 
units $\bbU[\zeta ]=\{1|0|1|0\}\bbU[\zeta ]$. They represent the isosceles triangles 
in collections $[\circ|\circ|\circ+1]$ and $[\circ|\circ+1|\circ+1]$. Among them the triangles 
$\{1|0|1|0\}\bbU^+[\zeta ]$, i.e.,  $\{l_j,\; k_j,\; l_j,\; -k_j\}$, 
form the collection $[\circ|\circ+1|\circ+1]$ whereas triangles $\{0,1,1,0\}\bbU^+[\zeta ]$, i.e.,
$\{k_j, l_j+2k_j, l_j+2k_j, k_j\}$, form the collection $[\circ|\circ|\circ+1]$. 
%The coordinates of the vertices are
%$$\big(0,0\big) -\hskip-4pt-\; \big(l_j,\;k_j \big) -\hskip-4pt-\; \big(l_j ,\; - k_j \big) \quad {\rm and}\quad
%\big(0,0\big) -\hskip-4pt-\; \big(k_j,\;l_j+2k_j \big) -\hskip-4pt-\; \big(l_j+2k_j ,\; k_j \big)$$
%respectively. 
%Accordingly
%$$d_j=4 k_j^2\quad {\rm and}\quad d_j=k_j^2 + (l_j+2k_j)^2$$
%and
%$$s_j=2 k_j l_j\quad {\rm and}\quad s_j=(l_j+2k_j)^2 - k_j^2.$$
Obviously, the full solution coset $\{1|0|1|0\}\bbU[\zeta ]$ is of Class A, and the corresponding 
$j^*(r) = 0$.

If $r$ is a rational prime with $r \bmod 12 \in \{5,7,11 \}$ then the corresponding norm equation 
does not have solutions as the decomposition $r=a^2+b^2+ab$ is impossible for 
$r \bmod 12 = 5, 11$, and the decomposition $r=a^2+b^2$ is impossible for $r \bmod 12 = 7$. 
For $r=2,3$ the corresponding norm equation also does not have solutions for exactly the same 
reason.

The remaining rational primes are $r \bmod 12 = 1$. Each rational prime of this form 
has a unique decomposition
(up to the multiplication by the unit from $\bbU[\zeta ]$)  into the product of 4 
conjugated algebraic primes
$$r = \{m|n|k|l\} \{m,-n, k, -l\} \{k,l,m,n\} \{k,-l,m,-n\},\eqno (7.1)$$
which define a single full coset of $\bbU[\zeta ]$. In other words, the integers $m,n,k,l$ 
constituting the leader of the coset are uniquely defined by $r$. As before, 
the coset of $\bbU[\zeta ]$ is partitioned into two cosets of $\bbU^+[\zeta ]$ corresponding to two 
triangle collections. 

%Comparing the resulting coset with the original coset $\bbU[\zeta ]$ one can see that the prime $p$ acts simultaneously as the magnifying rotation 
%$$T_{gh}(p)=\bpma \og & \oog \\ -\oog & \og\epma$$
%in $\bbZ[\e^3]$ subring of $\bbZ[\zeta ]$, the magnifying rotation
%$$T_{ab}(p)=\bpma \orr & \oor \\ -\oor & \orr\epma$$
%in $\bbZ[\e^4]$ subring of $\bbZ[\zeta ]$, and the linear transformation
%$$T_{st}(p)=\bpma t & s \\ -3s & t\epma$$
%in $\bbZ[\sqrt{3}]$ subring of $\bbZ[\zeta ]$. 
Comparing the resulting coset with the original coset $\bbU[\zeta ]$ one can see that the 
corresponding collections change from $[\circ|\circ|\circ+1]$ and $[\circ|\circ+1|\circ+1]$ to 
$[\circ|\circ+ a|\circ+ a+ b]$ and $[\circ|\circ+ b|\circ+ a+ b]$, respectively, where $a,b$ are
as in Eqn (5.19), (5.20). 

The above argument verifies assertions of Theorems 5 and 6 for a rational prime $r$.
\bigskip

{\bf 7.2.} The next step is to consider $r$ which is a square of a rational prime. Every rational 
prime $w$ with $w \bmod 12 = 5$ admits a decomposition into the product of two self-conjugated 
algebraic primes: $w = \{m|n|m|n\}\{m,-n,m,-n\}$, 
and such a decomposition is unique. In addition, our definition of the leader implies that 
$0 <m < n$. Performing the multiplication we obtain $w=m^2+n^2$ (a Gaussian 
representation). Consequently, for $r=w^2$ there exists a unique decomposition 
$$w^2 = \{m|n|m|n\}^2 \{m,-n,m,-n\}^2\eqno (7.2)$$
into the product of 4 conjugated algebraic integers,
which generates a single full solution coset to the corresponding norm equation. Comparing 
this solution coset with the solution coset $\{1|0|1|0\}\bbU^+[\zeta ]$ for $r=1$, we observe 
that the isosceles triangles in the coset  $\{m|n|m|n\}\bbU^+[\zeta ]$ have the vertices
$\{l_j m - k_j n,l_j n + k_j m,l_j m + k_j n,l_j n - k_j m\}$. They 
are obtained from the triangles $\{l_j,k_j,l_j,-k_j\}$ from the coset $\{1|0|1|0\}\bbU^+[\zeta ]$  
 by the magnifying rotation
$$T_\sq(w)=\begin{pmatrix} m & n & 0 & 0 \\ -n & m & 0 & 0 \\
0 & 0 & m & n \\ 0 & 0 & -n & m \end{pmatrix}.\eqno (7.3)$$
which can be understood as the magnifying rotation 
$\diy \begin{pmatrix}m & n \\ -n & m\end{pmatrix}$  in $\bbZ^2$.
%This matrix performs the rotation by the angle $\b(p)=\arctan{a\over b}$ 
%and the magnification by the factor $\sqrt{p}$. 
Under $T_\sq(w)$, the triangle collection is changing from $[\circ|\circ|\circ+1]$ 
to $[\circ|\circ|\circ+w]$. The same observation is true for the coset $\{1|0|1|0\}\bbU^
+[\zeta ](\zeta +\zeta ^2) = \{0|1|1|0\}\bbU^+[\zeta ]$ where under $T_\sq(w)$ 
the triangle collection is changing from $[\circ|\circ+1|\circ+1]$ to $[\circ|\circ+w|\circ+w]$.

Further, every rational prime $z$ with $z \bmod 12 = 7$ admits a decomposition into 
the product of two self-conjugated algebraic primes:
$$z = \{m|0|k|0\}\{k,0,m,0\},$$
and such a decomposition is unique. In addition, our definition of the leader implies that 
$0 <k < m$. Performing the multiplication we obtain $z=m^2+k^2-mk=(m-k)^2+k^2+(m-k)k$ 
(a L\"oschian number).
Consequently, for $r=z^2$ there exists a unique decomposition 
$$z^2 = \{m|0|k|0\}^2\{k,0,m,0\}^2\eqno (7.4)$$
into the product of 4 conjugated algebraic integers,
which generates a single full solution coset to the corresponding norm equation. The 
$\bbU[\zeta ]^+$-half of this coset (i.e. $\{m|0|k|0\}\bbU^+[\zeta ]$) consist of the triangles
$$\{l_j m, \;\; k_j(2k-m),\;\;  l_j k,\;\; -k_j(2m-k)\}$$
from collection $[\circ|\circ+2mk-m^2|\circ+2mk-k^2]$. This corresponds to the magnifying 
rotation 
$$T_\tr(z)= \begin{pmatrix} k & 0 & k-m & 0 \\ 0 & k & 0 & k-m \\
m-k & 0 & m & 0 \\ 0 & m-k & 0 & m \end{pmatrix}\eqno (7.5)$$
which can be understood as the magnifying rotation 
$\diy T_\tr(z)= \begin{pmatrix}m & k -m\\ m-k & m\end{pmatrix}$
%by the angle $\g(p)=\arctan{a\over b}$ and the magnification by the factor $\sqrt{p}$ 
in the dual lattice ${\bbA}_2$. The $\bbU[\zeta ]^+(\zeta +\zeta ^2)$-half of the full coset 
(i.e. $\{m|0|k|0\}\bbU^+[\zeta ](\zeta +\zeta ^2) = \{m-k,k,m,k-m\}\bbU^+[\zeta ]$) consists 
of the triangles from collection $[\circ|\circ+m^2-k^2|\circ+2mk-k^2]$. The original cosets 
$\{1|0|1|0\}\bbU^+[\zeta ]$ and $\{0|1|1|0\}\bbU^+[\zeta ]$ have the corresponding values 
of $ a$ and $ b$ swapped, and the same property remains true for the resulting cosets 
(under the transformation $T_\tr (z)$).

Next, every rational prime $q$ with $q \bmod 12 = 11$ admits a decomposition into the 
product of two self-conjugated algebraic primes
$$q = \{m|n|m|-n\}\{m,-n,m,n\},$$
and this decomposition is unique. Performing the multiplication, we obtain $m^2-3n^2=q$.
Consequently, for $r=q^2$ there exists a unique decomposition 
$$q^2 = \{m|n|m|-n\}^2\{m,-n,m,n\}^2\eqno (7.6)$$
into the product of 4 conjugated algebraic integers. This generates a single full solution coset 
$\{m|n|m|-n\}\bbU[\zeta ]$ to the corresponding norm equation. The 
$\{m|n|m|-n\}\bbU^+[\zeta ]$-half of this coset consist of isosceles triangles
$$\{l_j m + k_j 3n, \; l_j n + k_j m,\;  l_j m + k_j 3n,\; -l_j n  - k_j m\}$$
from collection $[\circ|\circ|\circ+q]$. The $\{m|n|m|-n\}\bbU^+[\zeta ](\zeta +\zeta ^2)$-half 
of this coset consists of isosceles triangles from collection $[\circ|\circ+q|\circ+q]$. This corresponds to the linear transformation
$$T_{\lz}(q)=\begin{pmatrix} m & -n & 0 & -2n \\ n & m & 2n & 0 \\
0 & 2n & m & n \\ -2n & 0 & -n & m\end{pmatrix}\eqno (7.7)$$
which also can be understood as the linear transformation 
$\diy T(q)= \begin{pmatrix}m & n \\ 3n & m\end{pmatrix}$
in the $(t,s)$ lattice, i.e $\bbZ[\sqrt{3}]$ sub-ring of $\bbZ[\zeta ]$.

The rational prime $2$ is a product of two self-conjugated algebraic primes
$$2 = \{1|1|1|1\}\{1,-1,1,-1\}$$
which are the same up to the multiplication by a unit
$\{1,-1,1,-1\}=\{1|1|1|1\}\zeta ^9$.
Consequently, for $r=2^2$ there exists a unique decomposition 
$$2^2 = \{1|1|1|1\}^2\{1,-1,1,-1\}^2\eqno (7.8)$$
into the product of 4 conjugated algebraic integers, which generates a single full 
solution coset $\{1|1|1|1\}\bbU[\zeta ]$ to the norm equation $N(\alpha )=2^2$. This 
coset consist of isosceles triangles
which are triangles from $\bbU[\zeta ]$ transformed by the magnifying rotation
$$T_\sq(2) = \begin{pmatrix}1 & 1 & 0 & 0\\ -1 & 1 & 0 & 0 \\ 
0 & 0 & 1 & 1 \\ 0 & 0 & -1 & 1\end{pmatrix}.\eqno (7.9)$$
The resulting collections are $[\circ|\circ|\circ+2]$ and $[\circ|\circ+2|\circ+2]$ for each of 
the two halves of the full coset, respectively. E.g., the $\{1|1|1|1\}\bbU^+[\zeta ]$ half is
$\{l_j - k_j, l_j  + k_j, l_j + k_j,l_j - k_j\}$.

The rational prime $3$ is a product of two self-conjugated algebraic primes
$$3 = \{2|0|1|0\}\{1,0,2,0\}$$
which are the same up to the multiplication by a unit $\{2|0|1|0\}=\{1,0,2,0\}\zeta ^{10}$.
Consequently, for $r=3^2$ there exists a unique decomposition 
$$3^2 = \{2|0|1|0\}^2\{1,0,2,0\}^2\eqno (7.10)$$
into the product of 4 conjugated algebraic integers. This again generates a single full solution coset 
$\{2|0|1|0\}\bbU[\zeta ]$ to the norm equation $N(\alpha )=3^2$. The resulting collections 
of isosceles triangles are $[\circ|\circ|\circ+3]$ and $[\circ|\circ+3|\circ+3]$ for each two halves 
of the full coset, respectively. E.g., the $\{2|0|1|0\}\bbU^+[\zeta ]$ half is $\{l_j 2, \; 0,\;  l_j,\; k_j 3\}$
or, equivalently, $\{k_j 3, \; l_j,\;  k_j 3,\; -l_j\}$. 
This correspond to the magnifying rotation
$$T_\tr(3) =\begin{pmatrix}1 & 0 & -1 & 0 \\ 0 & 1 & 0 & -1 \\ 
1 & 0 & 2 & 0 \\ 0 & 1 & 0 & 2 \end{pmatrix}\eqno (7.11)$$
in the dual lattice $\bbA_2$. 

By this moment we have verified the assertions of Theorems 5--8 for a value $r$ that is a minimal 
power of a rational prime. The next portion of the argument deals with values $r$ that are 
squares of the above. 
\bigskip

{\bf 7.3.} For $r = w^4$ where $w \bmod 12 = 5$ there are 2 different possibilities to represent
$w^4$ as the product of 4 conjugated algebraic integers. The first possibility contains all 4 conjugates 
of $\{m|n|m|n\}\{m,-n,m,-n\}=w$, which corresponds to $w^2I$, where $I$ is 
the identity matrix.  The second possibility contains all 4 conjugates of $\{m|n|m|n\}^2
=\{m^2-n^2,2mn,m^2-n^2,2mn\}$, which corresponds to $T_\sq(w)^2$. In the first case the 
resulting solution is the full coset $\{w|0|w|0\}\bbU[\zeta ]$
%$$\{l_j p, \; k_j p,\;  l_j p,\; -k_j p\}$$
containing  $w$ times magnified isosceles triangles $\{1|0|1|0\}\bbU^+[\zeta ]$ and 
$\{0|1|1|0\}\bbU^+[\zeta ]$ obtained earlier for $r=1$. The corresponding triangle collections 
are $[\circ|\circ+w^2|\circ+w^2]$ and $[\circ|\circ|\circ+w^2]$ respectively. The second case 
results in the isosceles triangles obtained from triangles in $\{1|0|1|0\}\bbU[\zeta ]$ by the 
magnifying rotation $T_\sq(w)^2$.
All details of this case can be recalculated from the case of $T_\sq(w)$ via the map
$$m \mapsto m^2-n^2,\quad n \mapsto 2mn.\eqno (7.12)$$
%$$\bpma a & b \\ -b & a\epma^2 = \bpma a^2-b^2 & 2ab \\ -2ab & a^2-b^2\epma.$$
%This matrix performs the rotation by the angle $2\b(p)$ and the magnification by the factor $p$. 
%In particular, the corresponding triangle types are $[\circ|\circ|\circ+p^2]$ and $[\circ|\circ+p^2|\circ+p^2]$. The $\bbU[\zeta ]^+$ half of this coset is
%$$\{l_j (m^2-n^2) - k_j 2 mn, \quad l_j 2mn + k_j (m^2-n^2),\quad  l_j 2mn + k_j (m^2-n^2),\quad l_j (m^2-n^2) - k_j 2 mn\}$$
%and the second half of the full coset can be obtained from it via the map $M$.
It is clear that two full cosets (for two cases above) contain congruent (at given value of $j$) 
triangles rotated with respect to each other by $2\arctan\,(n /m)$. This is a general mechanism 
behind a B0 family of solution cosets as stated in Theorem 7.

For $r = z^4$ where $z \bmod 12 = 7$ there are 2 different possibilities to represent $z^4$ as the
product of 4 conjugated algebraic integers. The first possibility contains all 4 conjugates of 
$\{m|0|k|0\}\{k,0,m,0\}=z$, which corresponds to $z^2I$. The second possibility
contains all 4 conjugates of $\{m|0|k|0\}^2 =\{2mk-k^2,0,2mk-m^2,0 \}$, which corresponds to $T_\tr(z)^2$. 
In the first case the resulting solution is the coset $\{z|0|z|0\}\bbU[\zeta ]$
%$$\{l_j p, \; k_j p,\;  l_j p,\; -k_j p\}$$
containing  $z$ times magnified isosceles triangles $\{1|0|1|0\}\bbU[\zeta ]$ obtained earlier for 
$r=1$.  The corresponding triangle collections are $[\circ|\circ|\circ+z^2]$ and 
$[\circ|\circ+z^2|\circ+z^2]$. 
All details of the second case can be recalculated from the case of $T_\tr (z)$ via the map
$$m \mapsto 2mk -k^2,\quad k \mapsto 2mk -m^2.\eqno (7.13)$$
It is clear that two full cosets (for the two cases above) contain triangles of the same area 
(again for the same value of $j$). This is true because in both cases the area of the triangle is the 
solution to the same equation $t^2-3s^2 = 1$, and this solution is unique for a given $j$.  This is a 
general mechanics behind a B1 family of solution cosets as stated in Theorem 8.

For $r = q^4$ where $q \bmod 12 = 11$ there are 2 different possibilities to represent $q^4$ as the
product of 4 conjugated algebraic integers. The first possibility contains all 4 conjugates of 
$\{m|n|m|-n\}\{m,-n,m,n\}=p$, which corresponds to $q^2I$. The second possibility
contains all 4 conjugates of $\{m|n|m|-n\}^2 =\{2m+3n,m+2n,2m+3n,-m-2n\}$, which corresponds 
to $T_\lz (q)^2$.  In the first case the resulting solution is the full coset $\{q|0|q|0\}\bbU[\zeta ]$
%$$\{l_j p, \; k_j p,\;  l_j p,\; -k_j p\}$$
containing  $q$ times magnified isosceles triangles $\{1|0|1|0\}\bbU[\zeta ]$ obtained earlier for $r=1$. 
All details of the second case can be recalculated from the case of $T_\lz (q)$ via the map
$$m \mapsto 2m+3n,\quad n \mapsto m+2n.\eqno (7.14)$$
The triangle area (for a given value of $j$) is different for the two resulting full cosets of isosceles triangles 
though the triangle collection is the same for both full cosets. To be more specific, the $\bbU^+[\zeta ]$-  
and $\bbU^+[\zeta ](\zeta +\zeta ^2)$-halves form collections $[\circ|\circ+q^2|\circ+q^2]$ and 
$[\circ|\circ|\circ+q^2]$ respectively. Note that triangles in each collections are mapped into itself 
under some $\bbZ^2$-symmetry,

For $r = 2^4$ there is only one possibility to represent $2^4$ as the product of 4 conjugated algebraic 
integers. This possibility contains all 4 conjugates of $\{1|1|1|1\}^2$ (which coincide with $\{1|1|1|1\}^2$ 
up to the multiplication by a unit). This is the case of $T_\sq(2)^2$, and the resulting solution is the full 
coset $\{2|0|2|0\}\bbU[\zeta ]$ containing triangles from collections $[\circ|\circ|\circ+4]$ and 
$[\circ|\circ+4|\circ+4]$.

For $r = 3^4$ there is only one possibility to represent $3^4$ as the product of 4 conjugated 
algebraic integers. This possibility contains all 4 conjugates of $\{2|0|1|0\}$ (which coincide with 
$\{2|0|1|0\}^2$ up to the multiplication by a unit). This is the case of $T_\tr(3)^2$, and the resulting 
solution is the full coset $\{3|0|3|0\}\bbU[\zeta ]$ containing triangles from collections 
$[\circ|\circ|\circ+9]$ and $[\circ|\circ+9|\circ+9]$.

The case of $r=p^2$ where a rational prime $p$ has $p \bmod 12 = 1$ is more diverse. In this case there 
are 4 distinct possibilities to represent $r$ as a product of 4 conjugated algebraic integers. Starting from 
the decomposition
$$p = \{m|n|k|l\} \{m,-n, k, -l\} \{k,l,m,n\} \{k,-l,m,-n\}$$
we can form 4 products
$$\begin{array}{l}\{m|n|k|l\}\{m|n|k|l\},\;\; \{m,-n, k, -l\}\{m|n|k|l\},\\
\{k,l,m,n\}\{m|n|k|l\},\quad \{k,-l,m,-n\}\{m|n|k|l\}\end{array}\eqno (7.15)$$
coupling $\{m|n|k|l\}$ with any of its conjugates (including itself). The conclusion is that each couple 
together with its 3 conjugates constitutes the representation of $r$ as the product of 4 conjugated 
algebraic integers. This corresponds to linear transformations $T_\sq(p)T_\tr(p)T_{\lz}(p)/p$, 
$T_\tr(p)$, $T_\sq(p)$ and 
$T_{\lz}(p)$ respectively and defines 4 distinct cosets of $\bbU [\zeta]$. Here 
$$\begin{array}{c}
T_\sq (p)=\begin{pmatrix}g & h & 0 & 0 \\ -h & g & 0 & 0 \\ 0 & 0 & g & h \\ 0& 0 & -h & g \end{pmatrix},\quad
T_\tr(p) = \begin{pmatrix} a & 0 & -b & 0 \\ 0 & a & 0 & -b \\ b & 0 & a & 0 \\ 0 & b & 0 & a \end{pmatrix},\\
T_\lz (p) = \begin{pmatrix}t & s & 0 & -2s \\ -s & t & 2s & 0 \\ 0 & 2s & t & -s \\ - 2s & 0 &s & t\end{pmatrix},
\end{array}\eqno (7.16)$$
where $g, h, a, b, s, t$ are are expressed via $m, n, k, l$ as in Eqns (5.19), (5.20).
In particular, the coset corresponding to $T_\sq(p)$ is $p$ times scaled and rotated coset for $r=1$. 
Cosets $T_\sq(p)T_\tr(p)T_{\lz}(p)/p$ and $T_{\lz} (p)$ have the same corresponding values of $(t_j, s_j)$ 
but different signatures $(a,b)$ and thus constitute a B1 family. The same is true for cosets $T_\sq(p)$ 
and $T_\tr(p)$.

This completes the part of the argument dealing with $r$ that are the squares of the 
values discussed  in {\bf 7.1} and {\bf 7.2}. This covers several cases in Theorems 5 -- 8. 
\bigskip

{\bf 7.4.} The situation becomes more involved when the corresponding power $\vrho (p) >2$. The 
corresponding 
amount of different cosets
$$\big((\vrho (p) + 2)(2 \vrho (p)^2+ 8 \vrho (p) + 15 + 9 (-1)^{\vrho (p)} \big)/ 48\eqno (7.17)$$
is known as OEIS A053307 or spreading number $\a_3(\vrho (p))$ (see 
\cite{GGR}, \cite{BT}, \cite{OEIS}). 
All these cosets can be partitioned into a single Class B1 family and remaining Class B2 families. We 
provide details for the cases $r=p^3$ and $r=p^4$ while the generalization to higher powers of $p$ 
is straightforward. If $r=p^3$ and
$$p = \{m|n|k|l\} \{m,-n, k, -l\} \{k,l,m,n\} \{k,-l,m,-n\}$$
then the corresponding B1 family contains two cosets
$$\{m|n|k|l\}\{m|n|k|l\}\{m|n|k|l\},\;\;\{m|n|k|l\}\{k,-l,m,-n\}\{k,-l,m,-n\}\eqno (7.18)$$
while the corresponding B2 family contains 3 cosets
$$\begin{array}{c}\{m|n|k|l\}\{m,-n, k, -l\}\{m,-n, k, -l\},
\{m|n|k|l\}\{k,l,m,n\}\{k,l,m,n\},\\
\{m,-n, k, -l\} \{k,l,m,n\} \{k,-l,m,-n\}.\end{array}\eqno (7.19)$$
The first and the third cosets from the B2 family contain congruent  triangles (for the same $j$). 
This conclusion becomes transparent if we use the representation of conjugates of $p$ via commuting 
linear transformations $T_{\bullet}(p)$:
$$\begin{array}{c}\{m|n|k|l\} \to \sqrt{1/p}T_\sq^{1/2}(p)T_\tr^{1/2}(p)T_\lz^{1/2}(p), \\
\{m,-n, k, -l\}  \to \sqrt{p}\;T_\sq^{-1/2}(p)T_\tr^{1/2}(p)T_\lz^{-1/2}(p), \\
\{k,l,m,n\}  \to \sqrt{p}\;T_\sq^{1/2}(p)T_\tr^{-1/2}(p)T_\lz^{-1/2}(p), \\
\{k,-l,m,-n\}  \to \sqrt{p}\; T_\sq^{-1/2}(p)T_\tr^{-1/2}(p)T_\lz^{1/2}(p). \end{array}\eqno (7.20)$$

Similarly, if $r=p^4$ then the corresponding B1 family ${\mathfrak B}_1$ consists of 3 cosets
$$ \begin{array}{l}\{m|n|k|l\}\{m|n|k|l\}\{m|n|k|l\}\{m|n|k|l\},\\
\{m|n|k|l\}\{m|n|k|l\}\{m|n|k|l\}\{k,-l,m,-n\},\\
\{m|n|k|l\}\{m|n|k|l\}\{k,-l,m,-n\}\{k,-l,m,-n\}.\end{array}\eqno (7.21)$$
The corresponding two B2 families ${\mathfrak B}^{(1)}_2$ and ${\mathfrak B}^{(2)}_2$ are
$$ \begin{array}{c} \{m|n|k|l\}\{m|n|k|l\}\{m|n|k|l\}\{m,-n, k, -l\},\\
 \{m|n|k|l\}\{m|n|k|l\}\{m|n|k|l\}\{k,l,m,n\},\\
 \{m|n|k|l\}\{m|n|k|l\}\{m,-n, k, -l\} \{k,-l,m,-n\},\\
 \{m|n|k|l\}\{m|n|k|l\}\{k,l,m,n\} \{k,-l,m,-n\}\end{array}\eqno (7.22)$$
and
$$ \begin{array}{c} \{m|n|k|l\}\{m|n|k|l\}\{m,-n, k, -l\}\{m,-n, k, -l\},\\
 \{m|n|k|l\}\{m|n|k|l\}\{k,l,m,n\}\{k,l,m,n\},\\
 \{m|n|k|l\}\{m|n|k|l\} \{m,-n, k, -l\} \{k,l,m,n\}, \\
  \{m|n|k|l\} \{m,-n, k, -l\} \{k,l,m,n\} \{k,-l,m,-n\}.\end{array}\eqno (7.23)$$
Among them the cosets 
$$ \begin{array}{c} \{m|n|k|l\}\{m|n|k|l\}\{m|n|k|l\}\{m,-n, k, -l\},\\
 \{m|n|k|l\}\{m|n|k|l\}\{m|n|k|l\}\{k,l,m,n\}\end{array}\eqno (7.24)$$
contain congruent triangles (labeled by the same $j$). Similarly, the cosets
$$ \begin{array}{c}
\{m|n|k|l\}\{m|n|k|l\}\{m,-n, k, -l\}\{m,-n, k, -l\},\\
\{m|n|k|l\} \{m,-n, k, -l\} \{k,l,m,n\} \{k,-l,m,-n\}\end{array}\eqno (7.25)$$
also contain congruent triangles (again with the same label $j$). For larger values
$\vrho (p)$ the corresponding B1 family contains $[\vrho (p)/2]+1$ cosets
$$\begin{array}{l} T_\sq^{\vrho (p) /2}(p)\; T_\tr^{\vrho (p) /2}(p)\; T_{\lz}^{\vrho (p) / 2}(p)\; p ^{- \vrho (p) / 2},\\
pT_\sq^{(\vrho (p)-2)/ 2}(p)\; T_\tr^{(\vrho (p)-2)/ 2}(p)\; T_{\lz}^{\vrho (p) / 2}(p)\; p ^{- \vrho (p) / 2},\\
p^2T_\sq^{(\vrho (p)-4)/2}(p)\; T_\tr^{(\vrho (p)-4)/ 2}(p)\; T_{\lz}^{\vrho (p) / 2}(p)\; p ^{- {\vrho (p) / 2}}, \ldots
\end{array}\eqno (7.26)$$
which all have different signatures (as they correspond to different powers of $T_\tr$). 

Adding the next power of $p$, i.e. considering $\vrho (p)+1$ instead of $\vrho (p)$ splits this B1 family 
into Class B1 and Class B2 families. The Class B1 family corresponds to multiplication by conjugates 
$T_\sq^{1/2}(p) T_\tr^{1/2}(p) T_\lz^{1/2}(p)/\sqrt{p}$ and 
$T_\sq^{-{1/2}}(p) T_\tr^{-{1/2}}(p)  T_\lz^{{1/2}}(p)\sqrt{p}$. The Class B2 family corresponds to 
the multiplication 
by conjugates $T_\sq^{1/2}(p) T_\tr^{-{1/2}}(p) T_\lz^{-{1/2}}(p)\sqrt{p}$ and $T_\sq^{-{1/2}}(p)T_\tr^{{1/2}}(p)
T_\lz^{-{1/2}}(p)\sqrt{p}$.

The arbitrary powers of prime factors $2, 3, q,w,z$ behave in a somewhat simpler fashion 
(we continue using letters $w, z,q$ for primes having values 
$5, 7, 11$, respectively,  modulo $12$). 
The factor $2^{2\vrho (2)}$ generates a single solution coset $T_\sq(2)^{\vrho (2)}$ which is 
the $2^{\vrho (2)}$-times  
magnification of the coset corresponding to $r=1$. Similarly, the factor $3^{2\vrho (3)}$ generates 
a single coset $T_\sq(3)^{\vrho (3)}$. 

Each factor 
$w^{2\vrho (w)}$ generates $[\vrho (w)/2]+1$ cosets 
$$T_\sq(w)^{\vrho (w)},\; T_\sq(w)^{\vrho (w)-2}w^2,\;  T_\sq(w)^{\vrho (w)-4}w^4, \ldots , \eqno (7.27)$$ 
which form a Class B0 family. We count here only non-negative powers of $T_\sq(w)$ because all triangles 
in these cosets are $\bbZ^2$-symmetric (and therefore isosceles).

Next, each factor 
$z^{2\vrho (z)}$ generates $[\vrho (z)/2]+1$ cosets 
$$T_\tr(z)^{\vrho (z)},\; T_\tr(z)^{\vrho (z)-2}z^2,\;  T_\tr(z)^{\vrho (z)-4}z^4, \ldots , \eqno (7.28)$$ 
which form a Class B1 family. We count here only non-negative powers of $T_\tr(z)$ because all 
triangles in these cosets are $\bbA_2$-symmetric.

Further, each factor 
$q^{2\vrho (q)}$ generates $[\vrho (q)/2]+1$ Class A cosets 
$$T_\tr(q)^{\vrho (q)},\; T_\tr(q)^{\vrho (q)-2}q^2,\; T_\tr(q)^{\vrho (q)-4}q^4, \ldots .\eqno (7.29)$$
They again consist of $\bbZ^2$-symmetric triangles.
\bigskip

{\bf 7.5.} The arguments in {\bf 7.1} -- {\bf 7.4} cover the situation where $r$ is an arbitrary 
power of a single rational prime number. To complete the 
proof of Theorems 5 -- 10 one needs to understand what happens when powers of 
different primes are 
combined with each other. Viz., consider $r=\prod_j q_j^2$ and assume that the corresponding 
statement of Theorem 6 is already verified for $r$ (the initial case $r=1$ is 
straightforward). Then $r=q^2 \prod_j q_j^2$ 
has a doubled number of solution cosets, and each of them is of Class A. This is true because 
every partitioning of  $\prod_j q_j^2$ into the product of 4 conjugated algebraic integers can be combined 
with each of two partitioning of $q$ into the product of 4 conjugated algebraic integers. In this way we
are able to form every  
partitioning of $q^2\prod_j q_j^2$ into the product of 4 conjugated algebraic integers. It is not hard to see 
that an additional factor $p$, $w^2$ or $z^2$ also doubles the amount of cosets, and all of them still are
of Class A. In addition, the factor $w^2$ magnifies and rotates all triangles in each 
coset while the factor $z^2$ 
changes the signature of all triangles in each coset. The factor $p$ performs 3 transformations 
$T_\sq^{1/2}$, $T_\tr^{{1/2}}$ and $T_\lz^{{1/2}}$ simultaneously. Observe that the 
presence of any additional factor $p,w^2$ or $z^2$ breaks the $\bbZ^2$-symmetry of the initial coset. 
This symmetry-breaking is responsible for the difference in the coset count between $r=q^{2\vrho (q)}$ and $r=*\cdot r$ where $*=p,w^2$ or $z^2$. A general case where each $q_j$ enters 
the product at an arbitrary even power $2\vrho (q_j)$ is similar but the total number of cosets 
cannot be expressed as an elementary function of values 
$\vrho (q_j)$. This comes as a result of different symmetry properties of different cosets
generated by single factor $q_j^{2\vrho (q_j)}$. Consequently, different pairs of cosets corresponding to 
$q_{j'}^{2\vrho (q_{j'})}$ and $q_{j''}^{2\vrho (q_{j''})}$ generate different amount of cosets corresponding 
to the product $q_{j'}^{2\vrho (q_{j'})} q_{j''}^{2\vrho (q_{j''})}$, and the resulting total amount turns out
to be smaller than $[(\vrho (q_{j'})/2]+1)([\vrho (q_{j''})/2]+1)$.

Now we can take any $r$ of the  above type and start adding factors $w_k^2$. The first of them, $w_1^2$, 
doubles the number of cosets and converts each Class A coset into two cosets forming a Class B0 family. 
The first of the cosets in the pair corresponds to the magnifying rotation $T_\sq(w_1)$ while the second one
correspond to magnification $w_1^2 I$. Adding another factor $w_2^2$ doubles the number of cosets 
in each B0 family, and so on. A generalization to the case $\prod_k w_k^{2\vrho (w_k)}$ is 
straightforward but is affected by a similar inability to express the resulting amount of B0 families 
as an elementary function of $\vrho (w_k)$.

The statements of Theorem 8 regarding Class B1 families are verified in a similar way. The only 
difference is that B1 families are generated by two types of prime factors: $p_i$ and 
$z_l^{2\vrho (z_l)}$. Recall that the factors $p_i^{\vrho (p_i)}$ with $\vrho (p_i) > 1$ also generate 
Class B1 families but they are accompanied by Class B2 families. 

A simpler way to generate Class B2 families is to combine factors $p_i$ and/or $w_k^{2\vrho (w_k)}$ 
with factors $z_l^{\vrho (z_l)}$. Observe that, due to ${\bbZ}^2$-symmetry of triangles generated by a 
single $w^2$ factor, adding a single $z^2$ factor produces a single Class B1 family containing two 
cosets with different signatures. Nevertheless, an additional $p$ or $w_1^2$ factor converts this family 
into two B2 families (due to the fact that $w^2p$ or $w^2w_1^2$ generate two cosets forming a
B0 family). Any additional factor $z_l^2$ doubles the amount of cosets in each of two families. 
Finally, any additional factor $w_k^2$ doubles the amount of B2 families.

In view of earlier considerations, additional factors $2^{2\vrho (2)}$ and/or $3^{2\vrho (3)}$  
change neither the amount of cosets nor the class type of the family nor the amount of cosets 
in a given family.

\section{Proof of Theorem S}\label{Sec8}

A characteristic pattern of sliding is that we have a  $\bbZ^2$-trapeze $OABW$ where both 
$\triangle OAW$ and $\triangle OBW$ are M-triangles, and side $OW$ is the base of sliding,
with $|OA|=|BW|$. Cf. Figures \ref{Fig1}C and \ref{Fig6}.  (The property $|OA|=|BW|$ may 
require to re-label some vertices.)  
Let $z = |AB|$. Then we say that each of the M-triangles $OAW$ and $OBW$ has 
$z$-sliding (for a given $D$).

The following Lemmas 8.1 and 8.3 are essentially repetitions of Theorem 2  and Lemma 2.1 from \cite{K}.

\bl For a $D$-admissible \rM-triangle $\triangle$ with $z$-sliding the corresponding 
doubled area $\tA(\triangle )$ is bounded from below as follows
$$\tA(\triangle ) > \frac{1}{2}D\sqrt{(3 D-1) (D + z)}. \eqno{(8.1)}$$
\el

\bp The trapeze $OABW$ is a cyclic quadrilateral; therefore 
$$|OB|\cdot |WA| = |AB| \cdot |OW| + |OA| \cdot |WB|. \eqno{(8.2)}$$
As $|OA| = |WB|$ and $|OB|= |WA|$, we have
$$|WA|^2 = |AB| \cdot |OW| + |OA|^2. \eqno{(8.3)}$$ 
Observe that $|OW|, |OA| \ge D$. Hence, $|WA|$ is minimal when $|AB| = z$ and 
$|OW| = |OA|=D$. In this case $|WA|^2=D^2+zD$ and 
$$
\tA(\triangle )\ge \frac{1}{2}D\sqrt{(3 D-1) (D + z)}.\eqno{(8.4)}$$
\ep

%\bigskip\noindent
%{\bf Lemma 9.2.} {\sl A $\bbZ^2$-triangle with vertices
%$$(0,0),\quad(\llbracket   x-\sqrt{3}y\rrbracket , \llbracket   y+\sqrt{3} x\rrbracket ),\quad (2x,2y) \eqno{(9.12)}$$ 
%has all sides longer than $D$ and a double area 
%$$\tA(\triangle ) \;<\; \left(D + 2\right)  \left(\frac{\sqrt{3}}{2} D + \sqrt{2}+\sqrt{3} \right)
%\;<\; \frac{1}{ 2}D\sqrt{(3 D-1) (D + z)} \eqno{(9.13)}$$
%as soon as
%$$z \ge 12,\quad D > 128 \quad {\rm and} \quad D + \sqrt{2} \le \sqrt{4x^2 + 4y^2} \le D + 2. \eqno{(9.14)}$$
%}
%
%\medskip\noindent
%{\bf Proof.} The first inequality in (9.13) is obtained from the last two inequalities in (9.14) by formally setting 
%$z=20$ in the proof of Lemma 9.2. The only part in Lemma 9.2 which is dependent on assumptions 
%regarding $z$ and $D$ is the second inequality in (9.7). It is a direct calculation to verify that under the
% assumptions on $z$ and $D$ in (9.14) the second inequality in (9.13) holds true.~\qed 

\bl For $D > 30$ and $z^2 \ge 14$ there is no $z$-sliding.
\el

\bp According to estimate (3.1) in Lemma~3.4.4 the doubled area of M-triangle $S(D)<\frac{{\sqrt 3}D^2}{2}+{\sqrt 2}D$. It is not hard to see that 
$$\frac{{\sqrt 3}D^2}{2}+{\sqrt 2}D < \frac{1}{2}D\sqrt{(3 D-1) (D + z)}.$$
for $D > 30$ and $z^2 \ge 14$.
Now one can apply Lemma 8.1.
\ep

For $x\in \bbR$ let $\{ x \}$ be the fractional part of $x$, 
$ \llbracket x\rrbracket$ an integer closest to $x$ and 
$$ \ldblbrace x \rdblbrace = \min\Big(\{ x \},\; 1-\{ x \} \Big) = \Big|x - \llbracket x\rrbracket \Big|. \eqno{(8.5)}$$
The most important idea in \cite{K} is the way to construct admissible triangles with small enough area which is described in Lemma~8.3 below. For small $z$ the corresponding estimate (8.7) is better than (3.1).

\bl A $\bbZ^2$-triangle $\triangle$ with vertices
$$(0,0),\quad(\llbracket   x-\sqrt{3}y\rrbracket , \llbracket   y+\sqrt{3} x\rrbracket ),\quad (2x,2y) \eqno{(8.6)}$$ 
has all sides longer than $D$ and a double area 
$$\tA (\triangle ) \;<\; \left(D + \frac{z}{10}\right)  \left(\frac{\sqrt 3}{2} D +\frac{z\sqrt{2}+z\sqrt{3}}{20}\right) 
< \frac{1}{ 2}D\sqrt{(3 D-1) (D + z)} \eqno{(8.7)}$$
as soon as 
$$1 \le z < 12,\quad D > 5,\quad \ldblbrace  x\sqrt{3} \rdblbrace,  \ldblbrace  y\sqrt{3} \rdblbrace
< \frac{z}{20}$$ 
and
$$D + \frac{z\sqrt{2}}{20} \le \sqrt{4x^2 + 4y^2} \le D + \frac{z}{10}. \eqno{(8.8)}$$
\el

\bp Observe that under the assumptions on $\sqrt{4x^2 + 4y^2}$ in (8.8) the distance $\ov D$ from 
the point $(x-\sqrt{3}y$, $y+\sqrt{3}x)$ to each of two sites $(0,0)$ and $(2x,2y)$ is upper-bounded by 
$$D + \frac{z\sqrt{2}}{20} \le \ov D \le D + \frac{z}{10}, \eqno{(8.9)}$$
as the corresponding triangle is equilateral. Owing to the assumptions on $\ldblbrace\,\cdot\,\rdblbrace$  
in (8.8), the distance between site $(\llbracket x-\sqrt{3}y\rrbracket , \llbracket   y+\sqrt{3}x\rrbracket )$ and point 
$(x-\sqrt{3}y,\;y+\sqrt{3}x)$ is at most $z{\sqrt 2}/20$. Thus, by the triangle inequality all sides of $\triangle$ 
are longer than $D$.

To estimate the doubled area of  $\triangle$, observe that the base of  $\triangle$ is shorter than 
$D + z/10$ while the corresponding height can be estimated as
$$\frac{\sqrt{3}}{2}\left(D +\frac{z}{10}\right) +\frac{z\sqrt{2}}{20}= \frac{\sqrt{3}}{2} D 
+\frac{z\sqrt{2}+z\sqrt{3}}{20}. \eqno{(8.10)}$$
Therefore,
$$\tA(\triangle) \le \left(D +\frac{z}{10}\right)  
\left(\frac{\sqrt{3}}{2} D +\frac{z\sqrt{2}+z\sqrt{3}}{20}\right). \eqno{(8.11)}$$
A direct calculation verifies that the last inequality in (8.7) holds true for $1<z<12$ and 
$D > 5$.
\ep

\bl Consider a countable set of integers   
$$\left\{x \in \bbZ:\; \ldblbrace\sqrt{3} x\rdblbrace < \frac{1}{ 20}\right\} \eqno{(8.12)}$$
and numerate them as $x_n$ in increasing order (so that $x_n < x_{n+1}$). Then for $n>1$ 
$$ x_{n+1}-x_n  \in \{4,11,15\}  \eqno{(8.13)}$$%(9.19)}$$
\el

\bp Consider the following fractional parts
$$\begin{array}{ll}%%
\del_{11} = \{11\sqrt{3}\} &\approx 0.05255888325764957, \\
\del_4 = \{4 \sqrt{3}\} &\approx 0.9282032302755088,\\
\del_{15} = \{15 \sqrt{3}\} &\approx 0.9807621135331566.\ena$$%(9.20)}$$
If $0.95 <\{\sqrt{3} x_n\} < 1.9-\del_4$ then $x_{n+1}=x_n+4$. If $1.9-\del_4 <\{\sqrt{3} x_n\} < 1$ or 
$0 < \{\sqrt{3} x_n\} < \del_{11}-0.05$ then $x_{n+1}=x_n+15$. If $\del_{11}-0.05  <\{\sqrt{3} x_n\} 
< 0.05$ then $x_{n+1}=x_n+11$. It is a direct calculation to verify that 
$\ldblbrace\sqrt{3} x_n+j\rdblbrace> 0.05$ for any other positive integer $j$ smaller than 15.
\ep

\bl Consider a countable set of integers 
$$\left\{d=4x^2+4y^2:\; x,y\in \bbZ,\; \ldblbrace\sqrt{3} x\rdblbrace ,  \ldblbrace\sqrt{3} y\rdblbrace 
< \frac{1}{ 20}\right\} \eqno{(8.14)}$$%(9.21)}$$
and numerate them as $d_n$ in increasing order, as above. Then for $d_n > 810000$
$$d_{n+1}-d_n \le 477 \root{4}\of{d_n}. \eqno{(8.15)}$$%(9.22)}$$
\el

\bp Given an integer $d\geq 1$, find $n_1$ such that
$$4x_{n_1}^2 \le d \le 4x_{n_1+1}^2. \eqno{(8.16)}$$%9.23)}$$
Then, according to Lemma 8.4, 
$$4x_{n_1+1}^2-4x_{n_1}^2 < 120x_{n_1}+ 900< 60\sqrt{d} + 900 \le 61\sqrt{d},  \eqno{(8.17)}$$%9.24)}$$
where the last inequality is true for $d \ge 810000$. Now find $n_2$ such that
$$4x_{n_1}^2 + 4x_{n_2}^2 \le d \le 4x_{n_1}^2 + 4x_{n_2+1}^2. \eqno{(8.18)}$$%9.25)}$$
Then
$$4x_{n_2+1}^2 - 4x_{n_2}^2 <
120x_{n_2}+ 900 <  61\sqrt{4x_{n_1+1}^2-4x_{n_1}^2}  
<  61 \sqrt{61\sqrt{d}} < 477 \root{4}\of{d}. \eqno{(8.19)}$$%9.26)}$$
Taking $d=d_n$ finishes the proof.
\ep

\bl For $1 \le z < 12$ and any $D = \sqrt{d}$, there exists $\ov d \in \bbN$ such that 
$$\left(D + \frac{z\sqrt{2}}{20}\right)^2 \le{\ov d}\le \left(D + \frac{z}{10}\right)^2,\;\;
{\ov d} = 4x^2+4y^2, \;\; \ldblbrace\sqrt{3} x\rdblbrace ,  
\ldblbrace\sqrt{3} y\rdblbrace < \frac{1}{ 20}\le \frac{z}{20}$$
as soon as $D > 66438468\; z^{-2}$. Consequently there is no $z$-sliding for $D > 66438468\; z^{-2}$.
\el

\bp Set ${\wh d}=  \left(D + z{\sqrt 2}/20\right)^2$. According to Lemma 8.3, there exists $\ov d$ 
such that
$${\wh d}\le{\ov d}\le {\wh d} + 477 \root{4}\of{\wh d}, \quad {\ov d}= 4x^2+4y^2, \quad 
\ldblbrace\sqrt{3} x\rdblbrace , \ldblbrace\sqrt{3} y\rdblbrace < \frac{1}{ 20}. \eqno{(8.20)}$$%9.27)}$$
Therefore,
$${\ov d}\le \left(D +\frac{z\sqrt{2}}{20}\right)^2 + 477\; \sqrt{D +\frac{z\sqrt{2}}{20}} \le 
\left(D +\frac{z}{10}\right)^2 \eqno{(8.21)}$$%9.28)}$$
as soon as $D > 66438468\; z^{-2}$. Indeed, in terms of ${\ov D}:= D/477$ and 
${\ov z}:=z \sqrt{2}/(20\cdot 477)$ we need
$$\begin{array} {l}%%
({\ov D} + \sqrt{2} {\ov z})^2- ({\ov D} + {\ov z})^2 - \sqrt{\diy\frac{{\ov D}}{477} + \frac{{\ov z}}{477}} > 0,\\
\big(2 (\sqrt{2}-1){\ov D}  {\ov z} + {\ov z}^2 \big)^2 > \diy\frac{{\ov D}}{477} +\frac{{\ov z}}{477},\\
{\ov z}^2 4(\sqrt{2}-1)^2 \left({\ov D}  +\diy\frac{{\ov z}}{2(\sqrt{2}-1)} \right) 477 > 1,\\
D  > \diy\frac{477^2 \;200}{4(\sqrt{2}-1)^2 z^2}\,,\ena\eqno{(8.22)}$$%9.29)}$$
which holds for $D > 66438468 z^{-2}= 292\cdot 477^2 z^{-2}$. 
\ep

\bl Suppose $z^2 = p^2 + q^2$ where integers $p,q\geq 1$ are mutually prime,
and there is no $z$-sliding. Then there is no $kz$-sliding for any positive integer $k$.
\el

\bp If there exists $kz$-sliding then the segment $AB$ contains $k+1$ lattice sites at distance 
$z$ from each other. Among them any two consecutive sites generate $z$-sliding which is impossible.
\ep

\bl Let $z^2 = p^2+q^2$ and ${\ov z}^2 = {\ov p}^2 + {\ov q}^2$ where integers $p,q\geq 1$ are 
mutually prime and ${\ov p},{\ov q}\geq 0$ are mutually prime, with ${\ov z}^2\geq 1$. If all 4 vertices of a 
trapeze implementing
a $z\cdot{\ov z}$-sliding belong to the skewed sub-lattice $\bbZ^2 \bpma{\ov p} & {\ov q} \\ -{\ov q} & {\ov p}\epma$ 
then this trapeze is a ${\ov z}$ times magnified version of a trapeze implementing $z$-sliding.
\el

\bp If for the scaled down triangle there exists an admissible triangle with a smaller area then the 
same is true for the scaled up version of this admissible triangle.
\ep

Given integers $z^2$, ${\ov z}^2$, $p,q$ and ${\ov p},{\ov q}$ as in Lemma 8.9, we define 
a $({\ov p}, {\ov q})$-{\it constrained}
$z$-sliding as a $z$-sliding with additional requirement that all vertices of an implementing trapeze $OABW$ 
belong to $\bbZ^2 \bpma {\ov p} & {\ov q} \\ -{\ov q} & {\ov p}\epma$.

\bl Let $z^2 = p^2+q^2$ where  integers $p,q\geq 1$ are mutually prime. If there is a
$z$-sliding then it is a $(p,q$)-constrained $z^2$-sliding.
\el

\bp The lemma is a direct consequence of the observation that for lattice $\bbZ^2$ hosting a
$(p,q)$-constrained $z^2$-sliding its $\bbZ^2 \bpma p & q \\ -q & p\epma$ sub-lattice hosts a
$z$-sliding.
\ep

\bl Suppose $z^2 = p^2+q^2$ where $p,q\geq 1$ are mutually prime integers. 
Consider an instance of $(p,q$)-constrained $z^2$-sliding implemented by the trapeze with vertices
$$(0,0),\quad({\ov x}, y),\quad ({\ov x} + z^2, y),\quad (x,0),\eqno{(8.23)}$$%9.30)}$$
where $x,{\ov x}, y \in \bbN,\; x \bmod z^2 = 0$ and $x - 2z^2 < 2{\ov x} < x - z^2$. Then 
$$\sqrt{3}{\ov x} < y < \sqrt{3}{\ov x} +\frac{4}{\sqrt{3}} z^2+1 \eqno{(8.24)}$$%9.31)}$$
and
$$\frac{\sqrt{3}}{2}  x - \sqrt{3} z^2 < y < \frac{\sqrt{3}}{2}x + \frac{1}{\sqrt{3}} z^2 + 1. \eqno{(8.25)}$$%9.32)}$$
Consequently, given ${\ov x}$ or $x$, only 3 corresponding values of $y$ are possible.
\el

\bp For any trapeze with vertices as in (8.23) which is a candidate for a $(p,q$)-constrained 
$z^2$-sliding, either $d = x^2$ or $d=(x-{\ov x}-z^2)^2 + y^2$. If $(x-z^2)^2 \ge d$ then the admissible 
triangle with vertices
$$(0,0),\quad ({\ov x}, y),\quad (x-z^2,0) \eqno{(8.26)}$$%9.33)}$$
has a smaller area. Therefore, $(x-z^2)^2  < (x-{\ov x}-z^2)^2 + y^2$, implying
$$ y^2 > 2 x {\ov x} -{\ov x}^2- 2{\ov x} z^2 > 3 {\ov x}^2 > 3 \left(\frac{x}{2} -z^2\right)^2. \eqno{(8.27)}$$%9.34)}$$
As $p$ and $q$ are mutually prime there exists ${\ov x} < \hat x < {\ov x} + z^2$ such that $(\hat x, y-1)$ 
belongs to sub-lattice $\bbZ^2 \bpma p & q \\ -q & p\epma$. The condition that the triangle with vertices
$$(0,0),\quad (\hat x, y-1),\quad (x-z^2,0)\eqno{(8.28)}$$%9.35)}$$
is not admissible is $(\min(\hat x, x-\hat x))^2 + (y-1)^2 < x^2$. Consequently,
$$(y-1)^2 < \frac{3}{4} x^2 + x z^2 - z^4 < 3 {\ov x}^2 + 8 z^2{\ov x} + 4 z^4. \eqno{(8.29)}$$%9.36)}$$
\ep 

\bl All existing 39 instances of $z$-sliding are listed in Table~2.
\el
%They correspond to either 1-sliding or $\sqrt{2}$-sliding.}??

\bp Due to Lemma 8.6, we only need to examine $D < {\ov D} = 66438468$ and due to 
Lemma 8.2 only with $z^2 < 14$. This investigation is done via a computer-assisted enumeration by {\tt ZSliding.java}. 
%See Section 9 and Supplement. 
Given integers $z^2=p^2+q^2$ and $d=D^2$ we need to (i) enumerate all $d$-admissible triangles which are 
{\it candidates} for $z$-sliding and (ii) perform an exhausting search for an {\it improvement}, i.e. a 
$d$-admissible triangle with smaller area. 

To construct an efficient algorithm we consider an auxiliary unit square lattice $\bbZ^2_{p,q}$ and 
embed into it $z^2$ times magnified original unit lattice. For definiteness we assume that this embedding 
is given by sublattice $\bbZ^2_{p,q} \bpma p & q \\ -q & p\epma$. Inside $\bbZ^2_{p,q}$ any $z^2$ times 
magnified trapeze $OABW$ that is a $z$-sliding candidate has an embedding with sides $AB$ and $OW$ 
being horizontal. Thus, placing the trapeze vertex $O$ at the origin we only need to choose vertex $W$ 
from the list of sites $(z^2, 0), (2z^2, 0), (3z^2,0), \ldots, ([{\ov D} / z^2]z^2 , 0)$, which grows linearly in ${\ov D}$. 
With vertex $W$ being selected, there are at most $3[z^2/2]$ possible positions for vertex $A$ satisfying 
$|OA| \le |OB|,|AW|, |BW|$. The $3[z^2/2]$-upper bound is a consequence of Lemma 8.11; accordingly,  
we obtain $O({\ov D})$ candidates for $z$-sliding to examine.

Now, for each of $z$-sliding candidates we need to find an improvement. The candidates without an 
improvement are the desired  instances of $z$-sliding. Given a trapeze $OABW$ with $|OW|= X$,
there are $O(X^2 z^{-4})$ possibilities to select the second vertex of the improvement and $O(1)$ 
possibilities to then select the third vertex. This amounts to $O({\ov D}^2)$ iterations 
for each of $O({\ov D})$ candidates. The total order of an $O({\ov D}^3)$
enumeration is unreachable to modern computers. To this end, we use an empirical approach 
to achieve a reasonable calculation time. Namely, given $W = (X,0)$, we perform the search for 
the second vertex of an improvement in the following order:
$$\begin{array}{c} %%
(X-1, 0), (X-1, 1), \ldots, (X-1, X),\\
(X-2, 0), (X-2, 1), \ldots, (X-2, X), \\
\ldots \ena\eqno{(8.30)}$$
In the worst situation it still requires $O(X^2 z^{-4})$ iterations. Because of that we additionally limit the search 
for the second coordinate to only 10 attempts, i.e. we perform the following enumeration
$$\begin{array}{c} %%
(X-1, 0), (X-1, 1), \ldots, (X-1, 10),\\
(X-2, 0), (X-2, 1), \ldots, (X-2, 10), \\
\ldots \ena\eqno{(8.31)}$$
Now it is $O(X z^{-2})$ iterations summing up to a total $O({\ov D}^2)$ iterations for the entire algorithm, which is 
still a considerable lot. Nevertheless, it appears that typically the improvement is found already among 
$(X-1, \cdot)$, and  the algorithm rarely goes to $(X-1, \cdot)$ or further. Because of that the actual amount of 
iterations performed by the algorithm is of order $O({\ov D})$, and on a mediocre personal computer the entire 
algorithm finishes in less than 5 hours.

Due to a non-exhaustive nature of our empirical search for an improvement, any identified instance of 
$z$-sliding needs to be additionally verified to be an actual instance of $z$-sliding. In practice, all 
39 instances  of sliding  detected 
by the algorithm are those listed in Table~2.

We also stress that our algorithm performs a search only among integers $z^2 < 14$ which correspond to 
mutually prime integers $p$ and $q$. The algorithm finds no sliding for $(p,q)$ pairs different from $(1,0)$ 
and $(1,1)$. Due to Lemmas 8.7--8.9, the last step is to verify that, being scaled more that $\sqrt{2}$ times 
but less than 12 times, none of the 39 sliding instances generate an additional instance of 
$z$-sliding. This appears to be the case, which completes our computer-assisted proof.
\ep

\vskip .5cm

{\bf Acknowledgement.} \ IS and YS thank the Math Department, Penn State University, for
hospitality and support. YS thanks St John's College, Cambridge, for long-term support.

\end{document}